# Some series and integrals involving the Riemann zeta function, binomial coefficients and the harmonic numbers

## Volume III

Donal F. Connon

18 February 2008

**Abstract**

In this series of seven papers, predominantly by means of elementary analysis, we establish a number of identities related to the Riemann zeta function, including the following:

$$Li_{s+1}(x) = \sum_{n=1}^{\infty} \frac{1}{n2^n} \sum_{k=1}^{n} \binom{n}{k} \frac{x^k}{k^s}$$

$$\sum_{n=1}^{\infty} \frac{1}{n^2} \sum_{k=1}^{n} \binom{n}{k} (-1)^k \frac{x^k}{k^s} = \log x \, Li_{s+1}(x) - (s+1)Li_{s+2}(x)$$

$$\sum_{n=1}^{\infty} \frac{w^n}{n^p} \sum_{k=1}^{n} \binom{n}{k} \frac{(-1)^k x^k}{k^s} = \frac{(-1)^{s+1} x}{\Gamma(s+1)} \int_0^1 \frac{\log^s t \, Li_{p-1}[w(1-xt)]}{1-xt} dt$$

$$\sum_{n=0}^{\infty} \frac{1}{2^{n+1}} \sum_{k=0}^{n} \binom{n}{k} \frac{(-1)^k}{(k+x)^s} = \sum_{k=0}^{\infty} \frac{(-1)^k}{(k+x)^s}$$

$$\frac{x}{s-1} \sum_{n=0}^{\infty} \frac{1}{n+1} \sum_{k=0}^{n} \binom{n}{k} (-1)^k \frac{x^k}{(k+1)^{s-1}} = \frac{\log x}{s-1} Li_{s-1}(x) + Li_s(x)$$

$$\sum_{n=0}^{\infty} \frac{w^n}{n+1} \sum_{k=0}^{n} \binom{n}{k} \frac{x^k}{(k+1)^{s-1}} = \frac{(-1)^{s-1}}{w\Gamma(s-1)} \int_0^1 \frac{\log^{s-2} v \log[1-w(1+xv)]}{1+xv} dv$$

$$\int_o^1 \left\{ t - \frac{Li_p[(1-y)t]}{1-y} \right\} \frac{dy}{\log y} = \sum_{n=1}^{\infty} \frac{t^{n+1}}{(n+1)^p} \sum_{k=0}^{n} \binom{n}{k} (-1)^{k+1} \log(1+k)$$

Whilst the paper is mainly expository, some of the formulae reported in it are believed to be new, and the paper may also be of interest specifically due to the fact that most of the various identities have been derived by elementary methods.

**(iii) Theorem 4.3:**

(4.4.45) $$Li_{s+1}(x) = \sum_{n=1}^{\infty} \frac{1}{n2^n} \sum_{k=1}^{n} \binom{n}{k} \frac{x^k}{k^s}$$

With $x = \pm 1$ this becomes

(4.4.45a) $$\varsigma(s+1) = \sum_{n=1}^{\infty} \frac{1}{n2^n} \sum_{k=1}^{n} \binom{n}{k} \frac{1}{k^s}$$

(4.4.45b) $$\varsigma_a(s+1) = \sum_{n=1}^{\infty} \frac{1}{n2^n} \sum_{k=1}^{n} \binom{n}{k} \frac{(-1)^{k+1}}{k^s}$$

**Proof:**

Let us consider the series $S_p(x,s)$ defined by

(4.4.46) $$S_p(x,s) = \sum_{n=1}^{\infty} \frac{1}{n^p 2^n} \sum_{k=1}^{n} \binom{n}{k} \frac{x^k}{k^s}$$

and consider the specific case

(4.4.46a) $$S_1(x,s) = \sum_{n=1}^{\infty} \frac{1}{n2^n} \sum_{k=1}^{n} \binom{n}{k} \frac{x^k}{k^s}$$

Using (4.4.33)

$$S_n(x) = \sum_{k=1}^{n} \binom{n}{k} \frac{x^k}{k^s} = \frac{1}{\Gamma(s)} \int_0^{\infty} u^{s-1} \left\{ (1 + xe^{-u})^n - 1 \right\} du$$

and the integral identity

(4.4.47) $$\frac{1}{n} = \int_0^{\infty} e^{-nv} dv$$

we obtain

(4.4.48a) $$S_1(x,s) = \sum_{n=1}^{\infty} \int_0^{\infty} e^{-nv} \frac{1}{2^n} \frac{1}{\Gamma(s)} \int_0^{\infty} u^{s-1} \left\{ (1 + xe^{-u})^n - 1 \right\} du\, dv$$

(4.4.48b)
$$= \frac{1}{\Gamma(s)} \int_0^\infty \int_0^\infty u^{s-1} \left\{ \frac{e^{-v}(1+xe^{-u})}{2-e^{-v}(1+xe^{-u})} - \frac{e^{-v}}{2-e^{-v}} \right\} du\, dv$$

(4.4.48c)
$$= -\frac{1}{\Gamma(s)} \int_0^\infty u^{s-1} \log(1-xe^{-u})\, du$$

If we let $s=1$ and $x=\pm 1$ in (4.4.48c), then G&R [74, p.525] tells us that

(4.4.49a)
$$\int_0^\infty \log(1-e^{-u})\, du = -\frac{\pi^2}{6}$$

(4.4.49b)
$$\int_0^\infty \log(1+e^{-u})\, du = \frac{\pi^2}{12}$$

and this at least provided me with some comfort that I was on the right track! It should be noted that the above two integrals are easily evaluated by making the substitution $y = e^{-u}$.

Now substitute $u = -\log t$ in (4.4.48c) and let $x=1$ to obtain

(4.4.50)
$$S_1(1,s) = \frac{(-1)^{s-1}}{\Gamma(s)} \int_0^1 \frac{\log^{s-1} t}{t} \log(1-t)\, dt$$

(see also (4.4.91j) re Nielsen's generalised polylogarithms).

With integration by parts we have

(4.4.51)
$$\int_0^1 \frac{\log^{s-1} t}{t} \log(1-t)\, dt = \log(1-t) \frac{\log^s t}{s} \Big|_0^1 + \frac{1}{s} \int_0^1 \frac{\log^s t}{1-t}\, dt$$

(4.4.51a)
$$= \frac{1}{s} \int_0^1 \frac{\log^s t}{1-t}\, dt$$

The integral in (4.4.51a) has already been calculated in (4.4.42) above and we have

(4.4.51b)
$$\int_0^1 \frac{\log^s t}{1-t}\, dt = \int_0^1 \frac{\log^s(1-t)}{t}\, dt = (-1)^s \varsigma(s+1)\Gamma(s+1)$$

Hence we obtain

(4.4.52) $$S_1(1,s) = \sum_{n=1}^{\infty} \frac{1}{n2^n} \sum_{k=1}^{n} \binom{n}{k} \frac{1}{k^s} = \varsigma(s+1) = \frac{(-1)^{s-1}}{\Gamma(s)} \int_0^1 \frac{\log^{s-1} t}{t} \log(1-t) dt$$

The excellent book, "Irresistible Integrals", by G. Boros and V.H. Moll [25, p.240] shows us another way to evaluate (4.4.51a) by expanding the denominator in a power series as follows:

(4.4.53a) $$\int_0^1 \frac{\log^s t}{1-t} dt = \sum_{k=0}^{\infty} \int_0^1 t^k \log^s t \, dt$$

(4.4.53b) $$= \sum_{k=0}^{\infty} (-1)^s \int_0^{\infty} u^s e^{-(k+1)u} du$$

(4.4.53c) $$= s! \sum_{k=0}^{\infty} \frac{(-1)^k}{(k+1)^{s+1}}$$

Therefore we get

(4.4.53d) $$\int_0^1 \frac{\log^s t}{1-t} dt = (-1)^s s! \varsigma(s+1)$$

Heretofore we have restricted our evaluation of $S_1(x,s)$ to the case where $x = 1$: we can in fact determine a more general result as set out below.

From (4.4.48c) we have

(4.4.54) $$S_1(x,s) = \frac{1}{\Gamma(s)} \int_0^{\infty} u^{s-1} \log(1 - xe^{-u}) du$$

and differentiating with respect to $x$ we obtain

(4.4.54a) $$S_1'(x,s) = \frac{1}{\Gamma(s)} \int_0^{\infty} \frac{u^{s-1} e^{-u}}{1 - xe^{-u}} du$$

(4.4.54b) $$= \frac{1}{\Gamma(s)} \int_0^{\infty} \frac{u^{s-1}}{e^u - x} du$$

(4.4.54c) $$= \frac{Li_s(x)}{x}$$

where the last step employs Appell's formula which was derived independently in (4.4.37). Therefore, using the series definition of the polylogarithm, we have

$$(4.4.55) \qquad S_1(t,s) = \int_0^t \frac{Li_s(x)}{x} dx = Li_{s+1}(t)$$

and hence

$$(4.4.56) \qquad Li_{s+1}(x) = \sum_{n=1}^{\infty} \frac{1}{n2^n} \sum_{k=1}^{n} \binom{n}{k} \frac{x^k}{k^s}$$

Yet another derivation is given below. Using (4.4.54) and the Maclaurin expansion of $\log(1-y)$ we obtain

$$(4.4.57) \qquad S_1(x,s) = -\frac{1}{\Gamma(s)} \int_0^{\infty} u^{s-1} \sum_{n=1}^{\infty} \frac{x^n e^{-nu}}{n} du$$

$$(4.4.57a) \qquad = -\frac{1}{\Gamma(s)} \sum_{n=1}^{\infty} \frac{x^n}{n} \int_0^{\infty} u^{s-1} e^{-nu} du$$

From [25, p.103] we have

$$(4.4.57b) \qquad \int_0^y u^m e^{-\mu u} du = \frac{m!}{\mu^{m+1}} - e^{-\mu y} \sum_{k=0}^{m} \frac{m! y^k}{k! \mu^{m-k+1}}$$

and hence as $y \to \infty$ we obtain

$$(4.4.57c) \qquad \int_0^{\infty} u^{s-1} e^{-nu} du = \frac{(s-1)!}{n^s}$$

The latter result could also be obtained rather more directly by the parametric differentiation of

$$(4.4.57d) \qquad \int_0^{\infty} e^{-nu} du = \frac{1}{n}$$

or, alternatively, by reference to the definition of the $\Gamma(x)$ function.

Therefore (4.4.57a) becomes

(4.4.57e) $$S_1(x,s) = \frac{1}{\Gamma(s)} \sum_{n=1}^{\infty} \frac{x^n}{n} \frac{(s-1)!}{n^s} = \sum_{n=1}^{\infty} \frac{x^n}{n^{s+1}} = Li_{s+1}(x)$$

An alternative proof is given below. From (4.4.25) we have

$$Li_s(x) = \frac{x}{\Gamma(s)} \int_0^{\infty} \frac{u^{s-1}}{e^u - x} du$$

and, after dividing this equation by $x$ and integrating, we obtain

(4.4.57f) $$\int_0^t \frac{Li_s(x)}{x} dx = \frac{1}{\Gamma(s)} \int_0^t dx \int_0^{\infty} \frac{u^{s-1}}{e^u - x} du$$

Reversing the order of integration gives us the elementary integral

$$\int_0^t \frac{1}{e^u - x} dx = -\log(e^u - x)\Big|_0^t$$

$$= -\log(e^u - t) + \log e^u$$

$$= -\log(1 - te^{-u})$$

Substituting this in (4.4.57f) we have

(4.4.57g) $$Li_{s+1}(t) = -\frac{1}{\Gamma(s)} \int_0^t u^{s-1} \log(1 - te^{-u}) du$$

and from (4.4.48c) we see that the right-hand side is equal to $S_1(x, s)$. Equations (4.4.49a) and (4.4.49b) are immediate consequences of (4.4.57g) with $s = 1$ and $t = \pm 1$.

**(iv) Theorem 4.4:**

(4.4.58) $$P_2(x,s) = \sum_{n=1}^{\infty} \frac{1}{n^2} \sum_{k=1}^{n} \binom{n}{k} \frac{(-1)^k x^k}{k^s} = -(s+1) Li_{s+2}(x) + \log x \, Li_{s+1}(x)$$

**Proof:**

Following the same approach as the previous theorem we have

$$(4.4.59) \qquad P_2(x,s) = \sum_{n=1}^{\infty} \int_0^{\infty} e^{-nv} \int_0^{\infty} e^{-nw} \frac{1}{\Gamma(s)} \int_0^{\infty} u^{s-1} \left\{ \left(1 + xe^{-u}\right)^n - 1 \right\} du\, dv\, dw$$

$$(4.4.59a) \qquad = \frac{1}{\Gamma(s)} \int_0^{\infty}\int_0^{\infty}\int_0^{\infty} u^{s-1} \left\{ \frac{e^{-(v+w)}\left(1 + xe^{-u}\right)}{1 - e^{-(v+w)}\left(1 + xe^{-u}\right)} - \frac{e^{-(v+w)}}{1 - e^{-(v+w)}} \right\} du\, dv\, dw$$

$$(4.4.59b) \qquad = -\frac{1}{\Gamma(s)} \int_0^{\infty}\int_0^{\infty} u^{s-1} \left\{ \log\left[1 - e^{-v}\left(1 + xe^{-u}\right)\right] - \log\left[1 - e^{-v}\right] \right\} du\, dv$$

I abandoned hope at this stage and reverted to the following simpler method.

Using (4.4.33) we have

$$S_n(x) = \sum_{k=1}^{n} \binom{n}{k} \frac{x^k}{k^s} = \frac{1}{\Gamma(s)} \int_0^{\infty} u^{s-1} \left\{ \left(1 + xe^{-u}\right)^n - 1 \right\} du$$

we have letting $x \to -x$

$$(4.4.60) \qquad S_n^a(x) = \sum_{k=1}^{n} \binom{n}{k} \frac{(-1)^k x^k}{k^s} = \frac{1}{\Gamma(s)} \int_0^{\infty} u^{s-1} \left\{ \left(1 - xe^{-u}\right)^n - 1 \right\} du$$

Making the substitution $u = -\log t$ we have

$$(4.4.61) \qquad S_n^a(x) = \frac{(-1)^{s-1}}{\Gamma(s)} \int_0^1 \frac{\log^{s-1} t}{t} \left\{ (1 - xt)^n - 1 \right\} dt$$

Integration by parts results in

$$(4.4.61a) \qquad S_n^a(x) = \frac{(-1)^{s-1}}{\Gamma(s)} \left\{ \frac{\log^s t}{s} \left[(1-xt)^n - 1\right] \Big|_0^1 + \frac{1}{s}\int_0^1 nx(1-xt)^{n-1} \log^s t\, dt \right\}$$

$$(4.4.61b) \qquad = \frac{(-1)^{s-1} nx}{s\Gamma(s)} \int_0^1 (1-xt)^{n-1} \log^s t\, dt$$

where we note that the very convenient factor $n$ "magically" appears in the numerator. Therefore we deduce

(4.4.61c) $$P_2(x,s) = \sum_{n=1}^{\infty} \frac{1}{n} \int_0^1 u^{n-1} du \cdot \frac{(-1)^{s-1} nx}{s\Gamma(s)} \int_0^1 (1-xt)^{n-1} \log^s t \, dt$$

(4.4.61d) $$= \frac{(-1)^{s-1} x}{\Gamma(s+1)} \int_0^1 du \int_0^1 \log^s t \sum_{n=1}^{\infty} u^{n-1} (1-xt)^{n-1} dt$$

and note that the $n$'s cancel out nicely. Completing the geometric series we have

(4.4.61e) $$= \frac{(-1)^{s-1} x}{\Gamma(s+1)} \int_0^1 du \int_0^1 \frac{\log^s t}{1-u(1-xt)} dt$$

(4.4.61f) $$= \frac{(-1)^{s-1} x}{\Gamma(s+1)} \int_0^1 \log^s t \, dt \int_0^1 \frac{1}{1-u(1-xt)} du$$

We have

$$\int_0^1 \frac{1}{1-u(1-xt)} du = -\frac{1}{1-xt} \log[1-u(1-xt)] \Big|_0^1$$

and hence we get

(4.4.61g) $$= \frac{(-1)^s x}{\Gamma(s+1)} \int_0^1 \frac{\log^s t \log xt}{1-xt} dt$$

(4.4.61h) $$= \frac{(-1)^s x}{\Gamma(s+1)} \int_0^1 \frac{\log^{s+1} t}{1-xt} dt + \frac{(-1)^s x \log x}{\Gamma(s+1)} \int_0^1 \frac{\log^s t}{1-xt} dt$$

(4.4.61i) $$= -(s+1) Li_{s+2}(x) + \log x \, Li_{s+1}(x)$$

where we have used (4.4.38b) in the last step.

Therefore we have shown that

(4.4.62) $$P_2(x,s) = \sum_{n=1}^{\infty} \frac{1}{n^2} \sum_{k=1}^{n} \binom{n}{k} \frac{(-1)^k x^k}{k^s} = -(s+1) Li_{s+2}(x) + \log x \, Li_{s+1}(x)$$

and letting $x=1$ confirms the validity of my conjecture in (4.2.51) that for $s \geq 0$

$$\varsigma(s+2) = \frac{1}{s+1} \sum_{n=1}^{\infty} \frac{1}{n^2} \sum_{k=1}^{n} \binom{n}{k} \frac{(-1)^{k+1}}{k^s}$$

When $s = 0$ we get

$$\sum_{n=1}^{\infty} \frac{1}{n^2} \sum_{k=1}^{n} \binom{n}{k} (-1)^k x^k = -Li_2(x) + \log x \, Li_1(x)$$

and we note that $\sum_{k=1}^{n} \binom{n}{k} (-1)^k x^k = \sum_{k=0}^{n} \binom{n}{k} (-1)^k x^k - 1 = (1-x)^n - 1$ and hence

$$\sum_{n=1}^{\infty} \frac{1}{n^2} \sum_{k=1}^{n} \binom{n}{k} (-1)^k x^k = \sum_{n=1}^{\infty} \frac{(1-x)^n - 1}{n^2} = Li_2(1-x) - \varsigma(2)$$

This then reduces to the Euler identity

$$Li_2(1-x) - \varsigma(2) = -Li_2(x) - \log x \log(1-x)$$

We also note that

$$\sum_{k=0}^{n} \binom{n}{k} (-1)^k k x^{k-1} = n(1-x)^{n-1}$$

and therefore

$$\sum_{n=1}^{\infty} \frac{1}{n^2} \sum_{k=1}^{n} \binom{n}{k} (-1)^k k x^k = x \sum_{n=1}^{\infty} \frac{1}{n^2} n(1-x)^{n-1} = x \sum_{n=1}^{\infty} \frac{(1-x)^{n-1}}{n}$$

With $s = -1$ in (4.4.62) we have

$$\sum_{n=1}^{\infty} \frac{1}{n^2} \sum_{k=1}^{n} \binom{n}{k} (-1)^k k \, x^k = \log x \, Li_0(x) = \frac{x \log x}{1-x}$$

in agreement with the above.

With $x = 1/2$ we have

(4.4.62a) $$P_2(1/2, s) = \sum_{n=1}^{\infty} \frac{1}{n^2} \sum_{k=1}^{n} \binom{n}{k} \frac{(-1)^{k+1}}{2^k k^s} = (s+1) Li_{s+2}(1/2) + \log 2 \, Li_{s+1}(1/2)$$

and hence we can find a specific result for $s = 1$

(4.4.62b) $$P_2(1/2, 1) = \sum_{n=1}^{\infty} \frac{1}{n^2} \sum_{k=1}^{n} \binom{n}{k} \frac{(-1)^k}{k 2^k} = \frac{7}{4} \varsigma(3) - \frac{\pi^2}{12} \log 2 - \frac{1}{6} \log^3 2$$

where we have used the Euler/Landen identities (3.43).

Having regard to (4.4.62b), and using (4.1.6), we can write

(4.4.63) $$\sum_{k=1}^{n}\binom{n}{k}\frac{(-1)^k}{k2^k} = \sum_{k=1}^{n}\frac{1}{k2^k} - H_n$$

Substituting this in (4.4.62b) we have

(4.4.63a) $$\sum_{n=1}^{\infty}\frac{1}{n^2}\sum_{k=1}^{n}\binom{n}{k}\frac{(-1)^k}{k2^k} = \sum_{n=1}^{\infty}\frac{1}{n^2}\left\{\sum_{k=1}^{n}\frac{1}{k2^k} - H_n\right\}$$

(4.4.63b) $$= \sum_{n=1}^{\infty}\frac{1}{n^2}\sum_{k=1}^{n}\frac{1}{k2^k} - \sum_{n=1}^{\infty}\frac{H_n}{n^2}$$

(4.4.63c) $$= \sum_{n=1}^{\infty}\frac{1}{n2^n}\sum_{k=n}^{\infty}\frac{1}{k^2} - 2\varsigma(3)$$

where we used (3.23) and (4.2.33) in the previous step

(4.4.63d) $$= \sum_{n=1}^{\infty}\frac{1}{n2^n}\left\{\varsigma(2) - \sum_{k=1}^{n-1}\frac{1}{k^2}\right\} - 2\varsigma(3)$$

(4.4.63e) $$= \log 2\varsigma(2) - \sum_{n=1}^{\infty}\frac{H_{n-1}^{(2)}}{n2^n} - 2\varsigma(3)$$

(4.4.63f) $$= \log 2\varsigma(2) - \sum_{n=1}^{\infty}\frac{H_n^{(2)}}{n2^n} + \sum_{n=1}^{\infty}\frac{1}{n^3 2^n} - 2\varsigma(3)$$

(4.4.63g) $$= \log 2\varsigma(2) - \sum_{n=1}^{\infty}\frac{H_n^{(2)}}{n2^n} + Li_3(1/2) - 2\varsigma(3)$$

Using (4.4.62b) and (4.4.63g) we can therefore show that

(4.4.64) $$\sum_{n=1}^{\infty}\frac{H_n^{(2)}}{n2^n} = \frac{5}{8}\varsigma(3)$$

(which in fact was previously derived in (3.49)). It may be possible to generalise the results by employing, for example, (4.1.10) and (4.1.17) but the output would appear to have a fairly complex structure. Reference should also be made to (4.4.156) where (4.4.64) is derived by a different method.

With $s = 1$ in (4.4.62) and using (4.1.6) we obtain

$$\sum_{n=1}^{\infty}\frac{1}{n^2}\sum_{k=1}^{n}\binom{n}{k}\frac{(-1)^k x^k}{k} = -\sum_{n=1}^{\infty}\frac{1}{n^2}\sum_{k=1}^{n}\frac{1-(1-x)^k}{k} = -2Li_3(x) + \log x\, Li_2(x)$$

and therefore we have

(4.4.64a) $$\sum_{n=1}^{\infty}\frac{1}{n^2}\sum_{k=1}^{n}\frac{(1-x)^k}{k} = 2\varsigma(3) - 2Li_3(x) + \log x\, Li_2(x)$$

This formula is used in (3.141) in Volume I.

With $s = 2$ in (4.4.62) we obtain

$$\sum_{n=1}^{\infty}\frac{1}{n^2}\sum_{k=1}^{n}\binom{n}{k}\frac{(-1)^k x^k}{k^2} = -3Li_4(x) + \log x\, Li_3(x)$$

and using (4.1.10) in Volume II(a) we may write this as

$$\sum_{n=1}^{\infty}\frac{1}{n^2}\sum_{k=1}^{n}\frac{1}{k}\left\{\sum_{l=1}^{k}\frac{1-(1-x)^l}{l}\right\} = 3Li_4(x) - \log x\, Li_3(x)$$

We have

$$\sum_{n=1}^{\infty}\frac{1}{n^2}\sum_{k=1}^{n}\frac{1}{k}\left\{\sum_{l=1}^{k}\frac{1-(1-x)^l}{l}\right\} = \sum_{n=1}^{\infty}\frac{1}{n^2}\sum_{k=1}^{n}\frac{H_k}{k} - \sum_{n=1}^{\infty}\frac{1}{n^2}\sum_{k=1}^{n}\frac{1}{k}\left\{\sum_{l=1}^{k}\frac{(1-x)^l}{l}\right\}$$

and from (4.1.14) we have

$$\sum_{k=1}^{n}\frac{H_k}{k} = \frac{1}{2}\left(H_n^{(1)}\right)^2 + \frac{1}{2}H_n^{(2)}$$

This then gives us

(4.4.64ai) $$\frac{1}{2}\sum_{n=1}^{\infty}\frac{\left(H_n^{(1)}\right)^2 + H_n^{(2)}}{n^2} - \sum_{n=1}^{\infty}\frac{1}{n^2}\sum_{k=1}^{n}\frac{1}{k}\left\{\sum_{l=1}^{k}\frac{(1-x)^l}{l}\right\} = 3Li_4(x) - \log x\, Li_3(x)$$

and with $x = 1$ we immediately see that

(4.4.64aii)
$$\frac{1}{2}\sum_{n=1}^{\infty}\frac{\left(H_n^{(1)}\right)^2+H_n^{(2)}}{n^2}=3\varsigma(4)$$

and therefore we see that

(4.4.64aiii)
$$3\varsigma(4)-\sum_{n=1}^{\infty}\frac{1}{n^2}\sum_{k=1}^{n}\frac{1}{k}\left\{\sum_{l=1}^{k}\frac{(1-x)^l}{l}\right\}=3Li_4(x)-\log x\,Li_3(x)$$

Letting $f(x)=\sum_{l=1}^{k}\frac{(1-x)^l}{l}$ we see that $f'(x)=-\sum_{l=1}^{k}(1-x)^{l-1}=\frac{(1-x)^k-1}{x}$. Therefore, differentiating (4.4.64aiii) results in

(4.4.64aiv)
$$\sum_{n=1}^{\infty}\frac{1}{n^2}\sum_{k=1}^{n}\frac{(1-x)^k}{k}=2\varsigma(3)-2Li_3(x)+\log x\,Li_2(x)$$

□

Dividing (4.4.64a) by $1-x$ and integrating results in

$$\sum_{n=1}^{\infty}\frac{1}{n^2}\sum_{k=1}^{n}\frac{(1-t)^k}{k^2}-\sum_{n=1}^{\infty}\frac{1}{n^2}\sum_{k=1}^{n}\frac{1}{k^2}=2\int_0^t\frac{\varsigma(3)-Li_3(x)}{1-x}dx+\int_0^t\frac{\log x\,Li_2(x)}{1-x}dx$$

See also (3.226b). The integral $\int_0^t\frac{\varsigma(3)-Li_3(1-u)}{1-u}du$ is evaluated in (4.4.168j) in Volume IV and hence we may determine the integral $\int_0^t\frac{\log x\,Li_2(x)}{1-x}dx$ in terms of polylogarithms and the series $\sum_{n=1}^{\infty}\frac{1}{n^2}\sum_{k=1}^{n}\frac{(1-t)^k}{k^2}$ (the Wolfram Integrator was unable to evaluate this integral).

Letting $t=1$ we immediately see that

$$2\int_0^1\frac{\varsigma(3)-Li_3(x)}{1-x}dx=-\int_0^1\frac{\log x\,Li_2(x)}{1-x}dx$$

Reference to (4.4.59b) shows that

$$P_2(x,s)=\sum_{n=1}^{\infty}\frac{1}{n^2}\sum_{k=1}^{n}\binom{n}{k}\frac{(-1)^k x^k}{k^s}=-(s+1)Li_{s+2}(x)+\log x\,Li_{s+1}(x)$$

(4.4.64b)
$$= -\frac{1}{\Gamma(s)} \int_0^\infty \int_0^\infty u^{s-1} \left\{ \log\left[1 - e^{-v}\left(1 + xe^{-u}\right)\right] - \log\left[1 - e^{-v}\right] \right\} du\, dv$$

Reference to (4.4.61c) shows that

$$P_3(x, s) = \sum_{n=1}^\infty \frac{1}{n^3} \sum_{k=1}^n \binom{n}{k} \frac{(-1)^k x^k}{k^s}$$

$$= \sum_{n=1}^\infty \frac{1}{n} \int_0^1 u^{n-1} du \int_0^1 v^{n-1} dv \cdot \frac{(-1)^{s-1} nx}{s\Gamma(s)} \int_0^1 (1-xt)^{n-1} \log^s t\, dt$$

$$= \frac{(-1)^{s-1} x}{\Gamma(s+1)} \int_0^1 du \int_0^1 dv \int_0^1 \frac{\log^s t}{1 - uv(1-xt)} dt$$

We see that

$$\int_0^1 \frac{dv}{1 - uv(1-xt)} = -\frac{\log[1 - uv(1-xt)]}{u(1-xt)} \bigg|_0^1 = -\frac{\log[1 - u(1-xt)]}{u(1-xt)}$$

$$\int_0^1 \frac{\log[1 - u(1-xt)]}{u(1-xt)} du = -\frac{Li_2[u(1-xt)]}{1-xt} \bigg|_0^1 = -\frac{Li_2(1-xt)}{1-xt}$$

Therefore we obtain

(4.4.64bi) $$P_3(x, s) = \sum_{n=1}^\infty \frac{1}{n^3} \sum_{k=1}^n \binom{n}{k} \frac{(-1)^k x^k}{k^s} = \frac{(-1)^{s-1} x}{\Gamma(s+1)} \int_0^1 \frac{\log^s t\, Li_2(1-xt)}{1-xt} dt$$

We have with $x = s = 1$

$$\int \frac{\log t\, Li_2(1-t)}{1-t} dt = \frac{1}{2}[Li_2(1-t)]^2$$

and therefore

$$\int_0^1 \frac{\log t\, Li_2(1-t)}{1-t} dt = -\frac{1}{2} \varsigma^2(2)$$

Hence we get

$$(4.4.64c) \quad P_3(1,1) = \sum_{n=1}^{\infty} \frac{1}{n^3} \sum_{k=1}^{n} \binom{n}{k} \frac{(-1)^{k+1}}{k} = \frac{1}{2}\varsigma^2(2)$$

and, using (4.1.7), this is equivalent to

$$\sum_{n=1}^{\infty} \frac{H_n^{(1)}}{n^3} = \frac{1}{2}\varsigma^2(2) = \frac{5}{4}\varsigma(4)$$

With $x = 1$ and $s = 2$ we get

$$\sum_{n=1}^{\infty} \frac{1}{n^3} \sum_{k=1}^{n} \binom{n}{k} \frac{(-1)^k}{k^2} = -\frac{1}{2} \int_0^1 \frac{\log^2 t \, Li_2(1-t)}{1-t} dt$$

Having regard to (4.1.15) this becomes

$$(4.4.64ci) \quad \int_0^1 \frac{\log^2 t \, Li_2(1-t)}{1-t} dt = \sum_{n=1}^{\infty} \frac{\left(H_n^{(1)}\right)^2}{n^3} + \sum_{n=1}^{\infty} \frac{H_n^{(2)}}{n^3}$$

Employing integration by parts we have

$$\int_a^1 \frac{\log^2 t \, Li_2(1-t)}{1-t} dt = -\frac{1}{2}[Li_2(1-a)]^2 \log a - \frac{1}{2}\int_a^1 \frac{[Li_2(1-t)]^2}{t} dt$$

$$= \frac{1}{2}\left(\varsigma^2(2) - [Li_2(1-a)]^2\right) \log a - \frac{1}{2}\int_a^1 \frac{[Li_2(1-t)]^2 - \varsigma^2(2)}{t} dt$$

Therefore we have as $a \to 0$

$$\int_0^1 \frac{\log^2 t \, Li_2(1-t)}{1-t} dt = -\frac{1}{2}\int_0^1 \frac{[Li_2(1-t)]^2 - \varsigma^2(2)}{t} dt$$

Using Euler's dilogarithm identity (1.6c) we see that

$$[Li_2(t)]^2 + [Li_2(1-t)]^2 + 2Li_2(t)Li_2(1-t) = \log^2 t \log^2(1-t) + \varsigma^2(2) - 2\varsigma(2)\log t \log(1-t)$$

We then get after dividing by $t$ and integrating

$$\int_0^1 \frac{[Li_2(t)]^2}{t}dt + \int_0^1 \frac{[Li_2(1-t)]^2 - \varsigma^2(2)}{t}dt + 2\int_0^1 \frac{Li_2(t)Li_2(1-t)}{t}dt =$$

$$\int_0^1 \frac{\log^2 t \log^2(1-t)}{t}dt - 2\varsigma(2)\int_0^1 \frac{\log t \log(1-t)}{t}dt$$

From (3.211f) in Volume I we have

$$\int_0^1 \frac{[Li_2(t)]^2}{t}dt = 2\varsigma(2)\varsigma(3) - 3\varsigma(5)$$

From (4.4.24zi) in Volume II(b) we note that

$$\int_0^1 \frac{[Li_2(1-t)]^2 - \varsigma^2(2)}{t}dt = -4\sum_{n=0}^{\infty}\frac{1}{(n+1)^2}\sum_{k=0}^{n}\binom{n}{k}\frac{(-1)^k}{(k+1)^3}$$

and we have from (4.2.28)

$$\sum_{k=0}^{n}\binom{n}{k}\frac{(-1)^k}{(1+k)^3} = \frac{1}{n+1}\left\{\frac{1}{2}(H_{n+1}^{(1)})^2 + \frac{1}{2}H_{n+1}^{(2)}\right\}$$

We therefore see that

$$\int_0^1 \frac{[Li_2(1-t)]^2 - \varsigma^2(2)}{t}dt = -2\sum_{n=0}^{\infty}\frac{(H_{n+1}^{(1)})^2 + H_{n+1}^{(2)}}{(n+1)^3} = -2\sum_{n=1}^{\infty}\frac{(H_n^{(1)})^2 + H_n^{(2)}}{n^3}$$

and this concurs with (4.4.64ci).

From (4.4.100gi) we have

$$\int_0^1 \frac{\log t \log(1-t)}{t}dt = \varsigma(3)$$

Integration by parts gives us

$$\int_0^1 \frac{Li_2(t)Li_2(1-t)}{t}dt = Li_3(t)Li_2(1-t)\Big|_0^1 - \int_0^1 \frac{Li_3(t)\log t}{t}dt$$

$$\int_0^1 \frac{Li_3(t)\log t}{t}\,dt = Li_4(t)\log t\Big|_0^1 - \int_0^1 \frac{Li_4(t)}{t}\,dt$$

Hence we obtain

(4.4.64cii) $$\int_0^1 \frac{Li_2(t)Li_2(1-t)}{t}\,dt = \varsigma(5)$$

Therefore we deduce that

(4.4.64ciii) $$\int_0^1 \frac{\log^2 t \log^2(1-t)}{t}\,dt = 4\varsigma(2)\varsigma(3) - \varsigma(5) - \frac{1}{2}\sum_{n=1}^{\infty}\frac{\left(H_n^{(1)}\right)^2}{n^3} - \frac{1}{2}\sum_{n=1}^{\infty}\frac{H_n^{(2)}}{n^3}$$

The Wolfram Integrator was not able to compute the above integral.

□

Let us now consider the expression

$$P_2(y,x,s) = \sum_{n=1}^{\infty}\frac{y^n}{n^2}\sum_{k=1}^{n}\binom{n}{k}\frac{(-1)^k x^k}{k^s}$$

Using (4.4.61d) this becomes

$$= \frac{(-1)^{s-1}xy}{\Gamma(s+1)}\int_0^1 du \int_0^1 \log^s t \sum_{n=1}^{\infty} y^{n-1}u^{n-1}(1-xt)^{n-1}\,dt$$

$$= \frac{(-1)^{s-1}xy}{\Gamma(s+1)}\int_0^1 du \int_0^1 \frac{\log^s t}{1-uy(1-xt)}\,dt$$

$$= \frac{(-1)^{s-1}xy}{\Gamma(s+1)}\int_0^1 \log^s t\,dt \int_0^1 \frac{1}{1-uy(1-xt)}\,du$$

$$= \frac{(-1)^{s-1}x}{\Gamma(s+1)}\int_0^1 \log^s t\,dt \left\{-\frac{1}{1-xt}\log[1-uy(1-xt)]\Big|_0^1\right\}$$

and hence we have

(4.4.64d) $$\sum_{n=1}^{\infty}\frac{y^n}{n^2}\sum_{k=1}^{n}\binom{n}{k}\frac{(-1)^k x^k}{k^s} = \frac{(-1)^s x}{\Gamma(s+1)}\int_0^1 \frac{\log^s t \log[1-y(1-xt)]}{1-xt}\,dt$$

When $y=1$ this reduces to (4.4.58) and we then obtain

$$\frac{(-1)^s x}{\Gamma(s+1)} \int_0^1 \frac{\log^s t \log(xt)}{1-xt} dt = -(s+1)Li_{s+2}(x) + \log x \, Li_{s+1}(x)$$

which may be seen from $\log(xt) = \log x + \log t$.

Letting $s = 2$, and reference to (3.115c) in Volume I shows that

(4.4.64e) $$\sum_{n=1}^{\infty} \frac{y^n}{n^2} \sum_{k=1}^{n} \binom{n}{k} \frac{(-1)^k}{k^2} = \frac{1}{2} \int_0^1 \frac{\log^2 t \log[1-y(1-t)]}{1-t} dt$$

$$= -\frac{1}{2} \sum_{n=1}^{\infty} \frac{y^n}{n^2} \left[ \left(H_n^{(1)}\right)^2 + H_n^{(2)} \right] = \int_0^y Li_3\left(\frac{-u}{1-u}\right) \frac{du}{u}$$

A division of (4.4.64d) by $y$ and an integration results in

$$\sum_{n=1}^{\infty} \frac{w^n}{n^3} \sum_{k=1}^{n} \binom{n}{k} \frac{(-1)^k x^k}{k^s} = \frac{(-1)^s x}{\Gamma(s+1)} \int_0^1 \frac{\log^s t}{1-xt} dt \int_0^w \frac{\log[1-y(1-xt)]}{y} dy$$

and since $\int_0^w \frac{\log[1-y(1-xt)]}{y} dy = -Li_2[y(1-xt)]\Big|_0^w = -Li_2[w(1-xt)]$

we have

(4.4.64f) $$\sum_{n=1}^{\infty} \frac{w^n}{n^3} \sum_{k=1}^{n} \binom{n}{k} \frac{(-1)^k x^k}{k^s} = \frac{(-1)^{s+1} x}{\Gamma(s+1)} \int_0^1 \frac{\log^s t \, Li_2[w(1-xt)]}{1-xt} dt$$

In the simplest case where $w = x = s = 1$ we get (4.4.64c) above. A further integration of (4.4.64f) results in

$$\sum_{n=1}^{\infty} \frac{w^n}{n^4} \sum_{k=1}^{n} \binom{n}{k} \frac{(-1)^k x^k}{k^s} = \frac{(-1)^{s+1} x}{\Gamma(s+1)} \int_0^1 \frac{\log^s t}{1-xt} dt \int_0^w \frac{Li_2[w(1-xt)]}{w} dw$$

and using

$$\int_0^w \frac{Li_2[w(1-xt)]}{w} dw = Li_3[w(1-xt)]$$

we obtain

(4.4.64g) $$\sum_{n=1}^{\infty}\frac{w^n}{n^4}\sum_{k=1}^{n}\binom{n}{k}\frac{(-1)^k x^k}{k^s} = \frac{(-1)^{s+1} x}{\Gamma(s+1)}\int_0^1 \frac{\log^s t\, Li_3[w(1-xt)]}{1-xt}\,dt$$

More generally, it is easily seen that

(4.4.64h) $$\sum_{n=1}^{\infty}\frac{w^n}{n^p}\sum_{k=1}^{n}\binom{n}{k}\frac{(-1)^k x^k}{k^s} = \frac{(-1)^{s+1} x}{\Gamma(s+1)}\int_0^1 \frac{\log^s t\, Li_{p-1}[w(1-xt)]}{1-xt}\,dt$$

With $p=1$, $w=1/2$ and $x \to -x$ we obtain

$$\sum_{n=1}^{\infty}\frac{1}{n 2^n}\sum_{k=1}^{n}\binom{n}{k}\frac{x^k}{k^s} = \frac{(-1)^{s+1} x}{\Gamma(s+1)}\int_0^1 \frac{\log^s t\, Li_0[(1+xt)/2]}{1+xt}\,dt$$

and reference to (4.4.45) shows that

$$Li_{s+1}(x) = \frac{(-1)^{s+1} x}{\Gamma(s+1)}\int_0^1 \frac{\log^s t\, Li_0[(1+xt)/2]}{1+xt}\,dt$$

or alternatively

$$Li_{s+1}(x) = \frac{(-1)^{s+1} x}{\Gamma(s+1)}\int_0^1 \frac{\log^s t\,[(1+xt)/2]}{(1+xt)\{1-[(1+xt)/2]\}}\,dt$$

and this simplifies to

(4.4.64i) $$Li_{s+1}(x) = \frac{(-1)^{s+1} x}{\Gamma(s+1)}\int_0^1 \frac{\log^s t}{1-xt}\,dt$$

Integration by parts shows that this is equivalent to (4.4.38e) in Volume II(b).

**(v)** In (4.4.25) we showed that

(4.4.65) $$Li_s(x) = \frac{x}{\Gamma(s)}\int_0^{\infty}\frac{u^{s-1}}{e^u - x}\,du$$

Using the expression

$$\frac{1}{A-B} - \frac{1}{A+B} = \frac{2B}{A^2 - B^2}$$

it is easily seen that

$$(4.4.66) \quad \frac{x}{\Gamma(s)}\int_0^\infty \frac{u^{s-1}}{e^u - x}du - \frac{x}{\Gamma(s)}\int_0^\infty \frac{u^{s-1}}{e^u + x}du = \frac{2x^2}{\Gamma(s)}\int_0^\infty \frac{u^{s-1}}{e^{2u} - x^2}du$$

$$= \frac{2x^2}{\Gamma(s)2^s}\int_0^\infty \frac{u^{s-1}}{e^u - x^2}du$$

Therefore we have

$$(4.4.67) \quad Li_s(x) + Li_s(-x) = 2^{1-s} Li_s(x^2)$$

(and this identity can also be derived directly from the series definition of $Li_s(x)$). From (4.4.67) it is easily seen that

$$(4.4.67a) \quad Li_s(-1) = (2^{1-s} - 1)Li_s(1) = (2^{1-s} - 1)\varsigma(s)$$

$$(4.4.67b) \quad Li_s(i) + Li_s(-i) = 2^{1-s} Li_s(-1)$$

where $i = \sqrt{-1}$.

As the following note shows, such simple identities can give rise to some marvellous formulae. In 1997 Bailey et al. [17] discovered the following fast converging series for $\pi$ by "a combination of inspired guessing and extensive searching using the PSLQ integer relation algorithm".

$$(4.4.68) \quad \pi = \sum_{n=0}^\infty \frac{1}{16^n}\left\{\frac{4}{8n+1} - \frac{2}{8n+4} - \frac{1}{8n+5} - \frac{1}{8n+6}\right\}$$

Having empirically discovered the identity, they were subsequently able to prove that it was equivalent to

$$(4.4.69a) \quad \pi = \int_0^{1/\sqrt{2}} \frac{4\sqrt{2} - 8x^3 - 4\sqrt{2}x^4 - 8x^5}{1 - x^8}dx$$

which on substituting $t = \sqrt{2}x$ becomes

$$(4.4.69b) \quad \pi = \int_0^1 \frac{16(1-t)}{t^4 - 2t^3 + 4t - 4}dt$$

The equivalence of (4.4.68) and (4.4.69a) is straightforward: it follows from the identity

(4.4.69c)
$$\int_0^{1/\sqrt{2}} \frac{x^{k-1}}{1-x^8}dx = \int_0^{1/\sqrt{2}} \sum_{j=0}^{\infty} x^{k-1+8j}dx$$

$$= \frac{1}{\sqrt{2}^k} \sum_{j=0}^{\infty} \frac{1}{16^j(8j+k)}$$

The proof that the integral (4.4.69a) evaluates to $\pi$ is an exercise in partial fractions most easily carried out in Maple or Mathematica.

In 1998, Broadhurst [36] considerably simplified the analysis by making the following observations:

(4.4.70a) $$Li_1(z) = -\log(1-z) = \sum_{n=1}^{\infty} \frac{z^n}{n}$$

With $w = \frac{1+i}{2}$ we have

(4.4.70b) $$Li_1(w) - \frac{1}{2}Li_1(1/2) = \frac{i\pi}{4}$$

Broadhurst then multiplies (4.4.70b) by $8(1-w) = 4(1-i)$ and takes the real part to obtain

(4.4.70c) $$\pi = 8\,\mathrm{Re}\{(1-w)Li_1(w)\} - 2Li_1(1/2)$$

$$= \sum_{n=0}^{\infty} \frac{1}{16^n}\left\{\frac{4}{8n+1} - \frac{2}{8n+4} - \frac{1}{8n+5} - \frac{1}{8n+6}\right\}$$

As a postscript, it may be mentioned that in 1998 by calculating the integrals $\int_0^x \frac{u^n}{1-u^8}du$ for $n = 0,1,...7$, Hirschhorn [79a] showed that

(4.4.70d) $$\log\left[\frac{(1+\sqrt{2}x+x^2)(1-x^2)^2}{(1-\sqrt{2}x+x^2)(1+x^4)}\right] + 2\tan^{-1}\left(\frac{\sqrt{2}x}{1-x^2}\right) + 2\tan^{-1}(x^2)$$

$$= 4\sqrt{2}\left[\sum_{n=1}^{\infty}\frac{x^{8n+1}}{8n+1} - \sum_{n=1}^{\infty}\frac{x^{8n+5}}{8n+5}\right] - 8\left[\sum_{n=1}^{\infty}\frac{x^{8n+4}}{8n+4} + \sum_{n=1}^{\infty}\frac{x^{8n+6}}{8n+6}\right]$$

and formula (4.4.70c) is a particular case where $x = 1/\sqrt{2}$. □

Letting $x = 1$ in (4.4.66) we have

(4.4.71a) $$\varsigma(s) = \frac{1}{(1-2^{1-s})\Gamma(s)} \int_0^\infty \frac{u^{s-1}}{e^u + 1} du$$

or alternatively

(4.4.71b) $$\varsigma_a(s) = \frac{1}{\Gamma(s)} \int_0^\infty \frac{u^{s-1}}{e^u + 1} du$$

and this formula is contained in [126, p.103]. This identity may also be obtained as follows.

Letting $t = 2u$ in $\varsigma(s) = \frac{1}{\Gamma(s)} \int_0^\infty \frac{t^{s-1}}{e^t - 1} dt$ we get

$$\varsigma(s) = 2 \frac{2^{s-1}}{\Gamma(s)} \int_0^\infty \frac{u^{s-1}}{e^{2u} - 1} du$$

and therefore subtraction results in

$$(1 - 2^{1-s})\varsigma(s) = \frac{1}{\Gamma(s)} \int_0^\infty \frac{u^{s-1}}{e^u - 1} du - \frac{1}{\Gamma(s)} \int_0^\infty \frac{2u^{s-1}}{e^{2u} - 1} du = \frac{1}{\Gamma(s)} \int_0^\infty \frac{u^{s-1}}{e^u + 1}$$

This method is employed in [80a, p.262] to provide an analytic continuation of the zeta function.

Now let us differentiate (4.4.65) with respect to $x$ to obtain

(4.4.72) $$Li'_s(x) = \frac{1}{\Gamma(s)} \int_0^\infty \frac{u^{s-1}}{e^u - x} du + \frac{x}{\Gamma(s)} \int_0^\infty \frac{u^{s-1}}{(e^u - x)^2} du$$

and simple algebra gives us

(4.4.72a) $$Li'_s(x) = \frac{1}{\Gamma(s)} \int_0^\infty \frac{u^{s-1} e^u}{(e^u - x)^2} du$$

A variant of (4.4.72a) is given in Gradshteyn and Ryzhik ("G&R") [74, p.354]. In (4.4.72) let $x = 1$ to obtain

(4.4.72b) $$Li'_s(1) = Li_s(1) + \frac{1}{\Gamma(s)} \int_0^\infty \frac{u^{s-1}}{(e^u - 1)^2} du$$

and hence we have

(4.4.72c) $$\int_0^\infty \frac{u^{s-1}}{(e^u-1)^2}\,du = \Gamma(s)\{Li_s'(1) - Li_s(1)\} = \Gamma(s)\{\varsigma(s-1) - \varsigma(s)\}$$

(which is also to be found in G&R [74, p.354, 3.423.1])

It is shown in [43] that

(4.4.72d) $$\int_0^\infty \frac{u^s}{(e^u+1)^2}\,du = C(s+1)\varsigma(s+1) - sC(s)\varsigma(s) \quad, s > 0$$

where $$C(s) = \Gamma(s)(1 - 2^{1-s})$$

Therefore we have

(4.4.72e) $$\int_0^\infty \frac{u^{s-1}}{(e^u+1)^2}\,du = C(s)\varsigma(s) - (s-1)C(s-1)\varsigma(s-1)$$

(4.4.72f) $$= \Gamma(s)(\varsigma_a(s) - \varsigma_a(s-1))$$

In their book [126, p.103] Srivastava and Choi state the following formula for $\varsigma(s)$

(4.4.72g) $$2\Gamma(s+1)(1-2^{1-s})\varsigma(s) = \int_0^\infty \frac{u^s e^u}{(e^u+1)^2}\,du \text{ (see the correction in (4.4.72j))}$$

Using integration by parts we have

(4.4.72h) $$\int_0^\infty \frac{u^s e^u}{(e^u+1)^2}\,du = \left.\frac{-u^s}{e^u+1}\right|_0^\infty + s\int_0^\infty \frac{u^{s-1}}{e^u+1}\,du$$

(4.4.72i) $$= s\int_0^\infty \frac{u^{s-1}}{e^u+1}\,du$$

Using (4.4.71a) this becomes

(4.4.72j) $$\int_0^\infty \frac{u^s e^u}{\left(e^u+1\right)^2}\,du = \Gamma(s+1)\left(1-2^{1-s}\right)\varsigma(s)$$

and it is therefore apparent that (4.4.72g) should not contain the factor 2 (this typographical error was also spotted by Apostol in his review of [126]).

Since $Li'_s(x) = \dfrac{Li_{s-1}(x)}{x} = \dfrac{1}{\Gamma(s-1)}\int_0^\infty \dfrac{u^{s-2}}{e^u - x}\,du$ we have

(4.4.72k) $$\frac{1}{\Gamma(s-1)}\int_0^\infty \frac{u^{s-2}}{e^u - x}\,du = \frac{1}{\Gamma(s)}\int_0^\infty \frac{u^{s-1}e^u}{\left(e^u - x\right)^2}\,du$$

Letting $u = \alpha t$ in (4.4.25) we get

$$Li_s(x) = \frac{x\alpha}{\Gamma(s)}\int_0^\infty \frac{\alpha^{s-1}t^{s-1}}{e^{\alpha t} - x}\,dt$$

and differentiating with respect to $\alpha$ we see that

$$\frac{\partial}{\partial \alpha} Li_s(x) = \frac{x}{\Gamma(s)}\int_0^\infty \frac{t^{s-1}\left[\left(e^{\alpha t} - x\right)s\alpha^{s-1} - \alpha^s t e^{\alpha t}\right]}{\left(e^{\alpha t} - x\right)^2}\,dt$$

Since $\dfrac{\partial}{\partial \alpha} Li_s(x) = 0$ we have with $\alpha = 1$

$$s\int_0^\infty \frac{t^{s-1}}{e^t - x}\,dt = \int_0^\infty \frac{t^s e^t}{(e^t - x)^2}\,dt$$

The above integral identity is simply another way of expressing the obvious relationship

$$s\Gamma(s)\frac{Li_s(x)}{x} = \Gamma(s+1)Li'_{s+1}(x)$$

From (4.4.15) we have

(4.4.73) $$g^{(p)}(x) = (-1)^p p!\sum_{k=0}^n \binom{n}{k}\frac{(-1)^k}{(k+x)^{p+1}} = \int_0^1 t^{x-1}(1-t)^n \log^p t\,dt$$

Similar integrals are reported in G&R [74, pp.535, 538] for values of $p$ equal to 1 and 2 (and their evaluation is considerably simplified by just computing successive derivatives of $g(x)$ as in (4.4.73) above). These integrals are shown below (and they can easily be transformed into (4.4.73) with the substitution $t = x^r$).

(4.4.73a) $$\int_0^1 x^{\mu-1}(1-x^r)^{\upsilon-1} \log x \, dx = \frac{1}{r^2} B\left(\frac{\mu}{r}, \upsilon\right)\left[\psi\left(\frac{\mu}{r}\right) - \psi\left(\frac{\mu}{r} + \upsilon\right)\right]$$

(4.4.73b)
$$\int_0^1 x^{\mu-1}(1-x^r)^{\upsilon-1} \log^2 x \, dx = \frac{1}{r^3} B\left(\frac{\mu}{r}, \upsilon\right)\left[\psi'\left(\frac{\mu}{r}\right) - \psi'\left(\frac{\mu}{r} + \upsilon\right) + \left\{\psi\left(\frac{\mu}{r}\right) - \psi\left(\frac{\mu}{r} + \upsilon\right)\right\}^2\right]$$

where $\text{Re}(\mu) > 0$, $\text{Re}(\upsilon) > 0$, $r > 0$ in both cases.

With $r = 1$ and $\nu = n+1$ in (4.4.73a) we have

(4.4.73c) $$\int_0^1 t^{x-1}(1-t)^n \log t \, dt = B(x, n+1)\left[\psi(x) - \psi(x+n+1)\right]$$

and comparing this with (4.4.73) results in

(4.4.73d) $$\sum_{k=0}^{n} \binom{n}{k} \frac{(-1)^k}{(k+x)^2} = -B(x, n+1)\left[\psi(x) - \psi(x+n+1)\right]$$

Similarly, with (4.4.73b) we have

(4.4.73e) $$2\sum_{k=0}^{n} \binom{n}{k} \frac{(-1)^k}{(k+x)^3} = B(x, n+1)\left[\psi'(x) - \psi'(x+n+1) + \{\psi(x) - \psi(x+n+1)\}^2\right]$$

which of course may be directly obtained by differentiating (4.4.73c).

**(vi)** We showed previously that

(4.4.74) $$Li_s(x) = \frac{x}{\Gamma(s)} \int_0^\infty \frac{u^{s-1}}{e^u - x} du$$

and letting $x \to -x$ we have

$$-Li_s(-x)\Gamma(s) = x\int_0^\infty \frac{u^{s-1}}{e^u + x}du$$

$$= \int_0^\infty \frac{u^{s-1}}{(e^u/x)+1}du$$

Now make the substitution $e^t = e^u/x$, $(u = t + \log x)$, to obtain

(4.4.74a) $\quad -Li_s(-x)\Gamma(s) = \int_{-\log x}^\infty \frac{(t+\log x)^{s-1}}{e^t+1}dt$

(4.4.74b) $\quad = \int_0^\infty \frac{(t+\log x)^{s-1}}{e^t+1}dt - \int_0^{-\log x} \frac{(t+\log x)^{s-1}}{e^t+1}dt$

First of all consider the case where $s = 2$.

(4.4.74c) $\quad -Li_2(-x)\Gamma(2) = \int_0^\infty \frac{(t+\log x)}{e^t+1}dt - \int_0^{-\log x} \frac{(t+\log x)}{e^t+1}dt$

(4.4.74d) $\quad = \int_0^\infty \frac{t}{e^t+1}dt + \log x \int_0^\infty \frac{1}{e^t+1}dt - \int_0^{-\log x} \frac{t}{e^t+1}dt - \log x \int_0^{-\log x} \frac{1}{e^t+1}dt$

(4.4.74e) $\quad = -Li_2(-1)\Gamma(2) - \log x\, Li_1(-1)\Gamma(1) - \int_0^{-\log x} \frac{t}{e^t+1}dt - \log x \int_0^{-\log x} \frac{1}{e^t+1}dt$

where we have used (4.4.74) to evaluate the first two integrals. The other integrals may be determined using the excellent resource of the Wolfram Integrator (http://integrals.wolfram.com/index.cgi) (and then making a human proof by differentiating the machine-produced result): however I prefer the following direct evaluation

(4.4.75a) $\quad I_1 = \int_0^{-\log x} \frac{1}{e^t+1}dt = \int_0^{-\log x} \frac{e^{-t}}{e^{-t}+1}dt = -\log(1+e^{-t})\Big|_0^{-\log x} = t - \log(1+e^t)\Big|_0^{-\log x}$

(4.4.75b) $\quad = -\log x - \log(1+1/x) + \log 2$

Integration by parts is used for the following integrals.

(4.4.76a) $$I_2 = \int_0^{-\log x} \frac{t}{e^t+1} dt = t\{t - t\log(1+e^t)\}\Big|_0^{-\log x} - \int_0^{-\log x} \{t - t\log(1+e^t)\} dt$$

We have by substituting $y = -e^t$

(4.4.76b) $$\int \log(1+e^t) dt = -\int \frac{\log(1-y)}{y} dy = -Li_2(y) = -Li_2(-e^t)$$

and hence we obtain

(4.4.76c) $$\int_0^{-\log x} \frac{t}{e^t+1} dt = \frac{t^2}{2} - t\log(1+e^t) - Li_2(-e^t)\Big|_0^{-\log x}$$

(4.4.76d) $$= \frac{1}{2}\log^2 x + \log x \log(1+1/x) - Li_2(-1/x) + Li_2(-1)$$

Similarly we have

(4.4.77) $$\int_0^{-\log x} \frac{t^2}{e^t+1} dt = \frac{t^3}{3} - t^2 \log(1+e^t) - 2tLi_2(-e^t) + 2Li_3(-e^t)\Big|_0^{-\log x}$$

(4.4.77a) $$= -\frac{1}{3}\log^3 x - \log^2 x \log(1+1/x) + 2\log x Li_2(-1/x) + 2Li_3(-1/x) - 2Li_3(-1)$$

Therefore, combining (4.4.74e), (4.4.75b) and (4.4.76d) we have

$$-Li_2(-x) = -Li_2(-1) + \log 2 \log x - \frac{1}{2}\log^2 x - \log x \log(1+1/x) + Li_2(-1/x) - Li_2(-1)$$

$$- \log x\{-\log x - \log(1+1/x) + \log 2\}$$

Therefore we obtain

(4.4.78) $$Li_2(-x) = -\varsigma(2) - \frac{1}{2}\log^2 x - Li_2(-1/x)$$

This formula is cited in [126, p.106] where it is derived in a much more elementary manner.

Letting $x = i$ in (4.4.78) we find

(4.4.78a) $$Li_2(-i) + Li_2(i) = -\varsigma(2) - \frac{1}{2}\log^2 i$$

and, since $\log i = \log 1 + i\frac{\pi}{2} = i\frac{\pi}{2}$, we have

(4.4.78b) $$Li_2(-i) + Li_2(i) = -\frac{1}{4}\varsigma(2)$$

(this was also shown in (4.4.67b)).

Letting $x = e^{i\pi}$ in (4.4.78) we find that $\varsigma(2) = \frac{\pi^2}{6}$.

With (4.4.74b) and (4.4.77a) it is easily shown that

(4.4.78c) $$Li_3(-x) - Li_3(-1/x) = -\frac{1}{3!}\log^3 x + 2\log x\, Li_2(-1)$$

Srivastava and Choi [126, p116] show how this formula can be generalised to

(4.4.78d) $$Li_n(-x) + (-1)^n Li_n(-1/x) = -\frac{1}{n!}\log^n x + 2\sum_{k=1}^{\left[\frac{n}{2}\right]} \frac{\log^{n-2k} x}{(n-2k)!} Li_{2k}(-1) \quad \forall n \geq 2$$

and they demonstrate how this can be used to derive an elementary proof of Euler's formula (1.7)

$$\varsigma(2n) = \sum_{k=1}^{\infty} \frac{1}{k^{2n}} = (-1)^{n+1} \frac{2^{2n-1}\pi^{2n} B_{2n}}{(2n)!}$$

Substituting $x = e^{it}$ in (4.4.78) we obtain

$$Li_2(-e^{it}) = -\varsigma(2) - \frac{1}{2}t^2 - Li_2(-e^{it})$$

and thereby results the familiar Fourier series [130, p.148]

(4.4.78e) $$\sum_{n=1}^{\infty}(-1)^n \frac{\cos nt}{n^2} = \frac{1}{4}t^2 - \frac{1}{2}\varsigma(2)$$

Similarly, with (4.4.78c) we have

(4.4.78f)
$$\sum_{n=1}^{\infty}(-1)^n \frac{\sin nt}{n^3} = \frac{1}{12}t^3 - \frac{t}{2}\varsigma(2)$$

Further Fourier series may be easily obtained by differentiating and integrating (4.4.78e) and (4.4.78f), and more general Fourier series may be obtained by making the same substitution in (4.4.78d).

Integrating (4.4.78) we obtain

$$\int_1^x \frac{Li_2(-t)+\varsigma(2)}{t}dt = -\frac{1}{6}\log^3 x - \int_1^x \frac{Li_2(-1/t)}{t}dt$$

We also see that

$$\int_1^x \frac{Li_2(-t)+\varsigma(2)}{t}dt = Li_3(-x) - Li_3(-1) + \varsigma(2)\log x$$

With the assistance of the Wolfram Integrator we get

(4.4.78g) $\int_1^x \frac{Li_2(-1/t)}{t}dt = \frac{1}{3}\log^3 x + Li_2(-1/x)\log x + Li_2(-x)\log x - Li_3(-x) + Li_3(-1)$

and therefore we simply recover (4.4.78)

$$\varsigma(2) = -\frac{1}{2}\log^2 x - Li_2(-1/x) - Li_2(-x)$$

With the substitution $u = 1/t$ we get

$$\int_1^x \frac{Li_2(-1/t)}{t}dt = -\int_1^{1/x} \frac{Li_2(-u)}{u}du = -Li_3(-u)\Big|_1^{1/x} = Li_3(-1) - Li_3(-1/x)$$

and hence we see that

(4.4.78h) $\frac{1}{3}\log^3 x + Li_2(-1/x)\log x + Li_2(-x)\log x - Li_3(-x) = -Li_3(-1/x)$

The next two theorems extend the Hasse/Sondow formulae (3.11) and (3.12).

**(vii) Theorem 4.5:**

(4.4.79) $$\sum_{n=0}^{\infty}\frac{1}{2^{n+1}}\sum_{k=0}^{n}\binom{n}{k}\frac{(-1)^k}{(k+x)^s} = \sum_{n=0}^{\infty}\frac{(-1)^n}{(n+x)^s} = \Phi(-1,s,x)$$

where $\Phi(z,s,x)$ is the Hurwitz-Lerch zeta function. With $x=1$ this becomes the Hasse/Sondow identity (3.11)

$$\varsigma_a(s) = \sum_{n=0}^{\infty}\frac{1}{2^{n+1}}\sum_{k=0}^{n}\binom{n}{k}\frac{(-1)^k}{(k+1)^s}$$

It is clear from (4.4.79) that $x$ can neither be zero nor a negative integer. I subsequently noted that equation (4.4.79) is a special case of Theorem 2.1 presented by Guillera and Sondow in [75aa].

**Proof:**

From (4.4.15) and (4.4.16) we have

(4.4.80) $$g^{(s-1)}(x) = (-1)^{s-1}(s-1)!\sum_{k=0}^{n}\binom{n}{k}\frac{(-1)^k}{(k+x)^s}$$

(4.4.80a) $$= \int_{0}^{1} t^{x-1}(1-t)^n \log^{s-1} t\, dt$$

Hence we get

(4.4.81) $$\sum_{n=0}^{\infty}\frac{1}{2^{n+1}}\sum_{k=0}^{n}\binom{n}{k}\frac{(-1)^k}{(k+x)^s} = \frac{(-1)^{s-1}}{2(s-1)!}\sum_{n=0}^{\infty}\frac{1}{2^n}\int_{0}^{1}t^{x-1}(1-t)^n\log^{s-1}t\,dt$$

(4.4.81a) $$= \frac{(-1)^{s-1}}{(s-1)!}\int_{0}^{1}\frac{t^{x-1}\log^{s-1}t}{1+t}dt$$

With the substitution $y = -\log t$ this becomes

(4.4.81b) $$= \frac{1}{(s-1)!}\int_{0}^{\infty}\frac{y^{s-1}e^{-y(x-1)}}{e^y+1}dy$$

and reference to [126, p.121] shows that (4.4.81b) is equal to $\Phi(-1,s,x)$ where the Hurwitz-Lerch zeta function $\Phi(z,s,x)$ is defined by

(4.4.82) $$\Phi(z, s, x) = \sum_{n=0}^{\infty} \frac{z^n}{(n+x)^s}$$

and has the integral representation [126, p.121]

(4.4.82a) $$\Phi(z, s, x) = \frac{1}{\Gamma(s)} \int_0^{\infty} \frac{y^{s-1} e^{-y(x-1)}}{e^y - z} dy$$

which is valid for $\text{Re}(x) > 0$; $|z| \leq 1, z \neq 1, \text{Re}(s) > 0; z = 1, \text{Re}(s) > 1$.

Therefore we have

(4.4.83) $$\sum_{n=0}^{\infty} \frac{1}{2^{n+1}} \sum_{k=0}^{n} \binom{n}{k} \frac{(-1)^k}{(k+x)^s} = \sum_{n=0}^{\infty} \frac{(-1)^n}{(n+x)^s}$$

and for $x = 1$ we have

(4.4.84) $$\sum_{n=0}^{\infty} \frac{1}{2^{n+1}} \sum_{k=0}^{n} \binom{n}{k} \frac{(-1)^k}{(k+1)^s} = \sum_{n=0}^{\infty} \frac{(-1)^n}{(n+1)^s} = \sum_{n=1}^{\infty} \frac{(-1)^{n+1}}{n^s} = \varsigma_a(s)$$

and (4.4.84) is the Hasse/Sondow identity (3.11).

**(viii) Theorem 4.6:**

(4.4.85) $$\frac{y}{s-1} \sum_{n=0}^{\infty} \frac{1}{n+1} \sum_{k=0}^{n} \binom{n}{k} \frac{(-1)^k y^k}{(k+1)^{s-1}} = \frac{\log y}{s-1} Li_{s-1}(y) + Li_s(y)$$

With $y = 1$ this becomes the formula discovered by Hasse [77] in 1930 (see (3.12))

(4.4.85a) $$\frac{1}{s-1} \sum_{n=0}^{\infty} \frac{1}{n+1} \sum_{k=0}^{n} \binom{n}{k} \frac{(-1)^k}{(k+1)^{s-1}} = Li_s(1) = \varsigma(s)$$

**Proof:**

Consider the sum

(4.4.86) $$S_n = \sum_{k=0}^{n} \binom{n}{k} \frac{x^{k+1}}{(k+1)^{s-1}} = \frac{1}{\Gamma(s-1)} \int_0^{\infty} u^{s-2} \sum_{k=0}^{n} \left\{ \binom{n}{k} \left[ e^{-u} x \right]^{k+1} \right\} du$$

(where we have used (4.4.32) above and note that the summation starts at $k = 0$ this time).

$$= \frac{1}{\Gamma(s-1)} \int_0^\infty u^{s-2} e^{-u} x (1+xe^{-u})^n \, du$$

We have the integral identity

(4.4.87) $$\frac{1}{n+1} = \int_0^\infty e^{-(n+1)t} dt = \int_0^\infty e^{-nt} e^{-t} dt$$

and hence we have

(4.4.88a) $$\sum_{n=0}^\infty \frac{S_n}{n+1} = \frac{1}{\Gamma(s-1)} \sum_{n=0}^\infty \int_0^\infty e^{-nt} e^{-t} dt \int_0^\infty u^{s-2} e^{-u} x (1+xe^{-u})^n \, du$$

(4.4.88b) $$= \frac{1}{\Gamma(s-1)} \int_0^\infty e^{-t} dt \int_0^\infty u^{s-2} e^{-u} x \sum_{n=0}^\infty e^{-nt} (1+xe^{-u})^n \, du$$

(4.4.88c) $$= \frac{1}{\Gamma(s-1)} \int_0^\infty e^{-t} dt \int_0^\infty u^{s-2} e^{-u} x \frac{1}{1-e^{-t}(1+xe^{-u})} \, du$$

Let us first of all integrate (4.4.88c) with respect to $t$. We have

(4.4.89) $$\int_0^\infty \frac{e^{-t} dt}{1-Ae^{-t}} = \frac{1}{A} \log(1-Ae^{-t}) \Big|_0^\infty$$

$$= -\frac{1}{A} \log(1-A)$$

$$= -\frac{1}{1+xe^{-u}} \log(-xe^{-u}) \quad, \text{ where } A = 1+xe^{-u}.$$

At this stage we now specify that $x$ is a strictly negative number and substitute $y = -x$ ($y > 0$). We therefore have

(4.4.90a) $$= \frac{1}{\Gamma(s-1)} \int_0^\infty \frac{-y(\log y + u) u^{s-2} e^{-u} du}{1-ye^{-u}}$$

(4.4.90b) $$= -\frac{y \log y}{\Gamma(s-1)} \int_0^\infty \frac{u^{s-2} du}{e^u - y} - \frac{y}{\Gamma(s-1)} \int_0^\infty \frac{u^{s-1} du}{e^u - y}$$

(4.4.90c) $$= -\log y\, Li_{s-1}(y) - \frac{\Gamma(s)}{\Gamma(s-1)} Li_s(y)$$

where we have used the integral identity in (4.4.25). This therefore proves that

(4.4.91) $$\frac{y}{s-1}\sum_{n=0}^{\infty}\frac{1}{n+1}\sum_{k=0}^{n}\binom{n}{k}\frac{(-1)^k y^k}{(k+1)^{s-1}} = \frac{\log y}{s-1} Li_{s-1}(y) + Li_s(y)$$

and with $y = 1$ this reverts to the second Hasse identity referred to in (3.12)

$$\varsigma(s) = \frac{1}{s-1}\sum_{n=0}^{\infty}\frac{1}{n+1}\sum_{k=0}^{n}\binom{n}{k}\frac{(-1)^k}{(k+1)^{s-1}}$$

With $y = 1/2$ in (4.4.91) we have

(4.4.91a) $$\frac{1}{2(s-1)}\sum_{n=0}^{\infty}\frac{1}{n+1}\sum_{k=0}^{n}\binom{n}{k}\frac{(-1)^k}{2^k(k+1)^{s-1}} = \frac{\log(1/2)}{s-1} Li_{s-1}(1/2) + Li_s(1/2)$$

and, using the Euler/Landen identities (3.43) for $Li_s(1/2)$ for $s = 2$ and $s = 3$ we obtain respectively

(4.4.91b) $$\frac{1}{2}\sum_{n=0}^{\infty}\frac{1}{n+1}\sum_{k=0}^{n}\binom{n}{k}\frac{(-1)^k}{2^k(k+1)} = \frac{\pi^2}{12} - \frac{3}{2}\log^2 2$$

(4.4.91c) $$\frac{1}{4}\sum_{n=0}^{\infty}\frac{1}{n+1}\sum_{k=0}^{n}\binom{n}{k}\frac{(-1)^k}{2^k(k+1)^2} = \frac{7}{8}\varsigma(3) - \frac{\pi^2}{8}\log 2 + \frac{5}{12}\log^3 2$$

Dividing (4.4.91) by $y$ and integrating over the interval $[0, x]$ results in

(4.4.91d) $$\frac{x}{s-1}\sum_{n=0}^{\infty}\frac{1}{n+1}\sum_{k=0}^{n}\binom{n}{k}\frac{(-1)^k x^k}{(k+1)^s} = \int_0^x\left(\frac{\log y}{(s-1)y} Li_{s-1}(y) + \frac{Li_s(y)}{y}\right)dy$$

Using integration by parts we get

$$\int_0^x \frac{\log y}{y} Li_{s-1}(y)dy = \frac{1}{2}\log^2 x\, Li_{s-1}(x) - \frac{1}{2}\int_0^x \frac{\log^2 y}{y} Li_{s-2}(y)dy$$

Continuing the integration by parts we obtain

$$\int_0^x \frac{\log y}{y} Li_{s-1}(y)dy = \frac{1}{2!}\log^2 x \, Li_{s-1}(x) - \frac{1}{3!}\log^3 x \, Li_{s-2}(x) + \frac{1}{3!}\int_0^x \frac{\log^3 y}{y} Li_{s-3}(y)dy$$

and this becomes

$$\int_0^x \frac{\log y}{y} Li_{s-1}(y)dy = \sum_{n=1}^{s-3} \frac{(-1)^{n+1}}{(n+1)!}\log^{n+1} x \, Li_{s-n}(x) + \frac{(-1)^{s-1}}{(s-2)!}\int_0^x \frac{\log^{s-2} y}{y} Li_2(y)dy$$

By a final integration by parts we have

$$\int_0^x \frac{\log^{s-2} y}{y} Li_2(y)dy = \frac{1}{s-1}\log^{s-1} x \, Li_2(x) + \frac{1}{s-1}\int_0^x \frac{\log^{s-1} y \log(1-y)}{y}dy$$

and hence we get

(4.4.91e)

$$\int_0^x \frac{\log y}{y} Li_{s-1}(y)dy = \sum_{n=1}^{s-2} \frac{(-1)^{n+1}}{(n+1)!}\log^{n+1} x \, Li_{s-n}(x) + \frac{(-1)^{s-1}}{(s-1)!}\int_0^x \frac{\log^{s-1} y \log(1-y)}{y}dy$$

An alternative integration by parts provides us with the simpler expression

(4.4.91f) $$\int_0^x \log y \frac{Li_{s-1}(y)}{y} dy = \log x \, Li_s(x) - \int_0^x \frac{Li_s(y)}{y}dy$$

$$= \log x \, Li_s(x) - Li_{s+1}(x)$$

We therefore obtain

(4.4.91g) $$\frac{(-1)^{s-1}}{(s-1)!}\int_0^x \frac{\log^{s-1} y \log(1-y)}{y}dy = -Li_{s+1}(x) + \sum_{n=0}^{s-2} \frac{(-1)^n}{(n+1)!}\log^{n+1} x \, Li_{s-n}(x)$$

$$= \sum_{m=0}^{s-1} \frac{(-1)^{m+1}}{m!}\log^m x \, Li_{s-m+1}(x)$$

This formula was also derived by Lewin [100, p.199] by using repeated integration by parts of $Li_n(x) = \int_0^x \frac{Li_{n-1}(t)}{t}dt$, but his result does not appear to be quite correct.

With $x=1$ we get

$$(4.4.91\text{ga}) \qquad \int_0^1 \frac{\log^{s-1} y \log(1-y)}{y} dy = (-1)^s (s-1)!\varsigma(s+1)$$

By letting $s \to s+1$, reference to (4.4.91) readily shows that the left-hand side of (4.4.91d) is equal to

$$\frac{\log x}{s-1} Li_s(x) + \frac{s}{s-1} Li_{s+1}(x)$$

Since $\int_0^x \frac{Li_s(y)}{y} dy = Li_{s+1}(x)$ we obtain from (4.4.91d) and (4.4.91e)

$$\frac{\log x}{s-1} Li_s(x) + \frac{s}{s-1} Li_{s+1}(x) = \frac{1}{s-1} \sum_{n=1}^{s-2} \frac{(-1)^{n+1}}{(n+1)!} \log^{n+1} x \, Li_{s-n}(x) + \frac{1}{s-1} \frac{(-1)^s}{(s-1)!} \int_0^x \frac{\log^{s-1} y \log(1-y)}{y} dy + Li_{s+1}(x)$$

and hence we get another proof of (4.4.91g).

The above result is related to Nielsen's generalised polylogarithms which are defined by

$$(4.4.91\text{h}) \qquad S_{n,p}(z) = \frac{(-1)^{n+p-1}}{(n-1)!p!} \int_0^1 \frac{\log^{n-1} t \log^p(1-zt)}{t} dt \qquad n, p = 1, 2, 3...$$

and were briefly referred to by Maximon in [101b].

Maximon [101b] also reports that

$$(4.4.91\text{i}) \qquad S_{n+1,1}(z) = Li_{n+2}(z) = \frac{(-1)^n}{n!} \int_0^1 \frac{\log^n t \log(1-zt)}{t} dt \qquad n,= 0,1,2...$$

where in particular we have with $z=1$

$$(4.4.91\text{j}) \qquad \varsigma(n+2) = \frac{(-1)^n}{n!} \int_0^1 \frac{\log^n t \log(1-t)}{t} dt \qquad n,= 0,1,2...$$

It should be noted that (4.4.91j) is a direct consequence of (4.4.91g) when we put $x=1$ and $s=n+1$. See also (4.4.50).

Nielsen's generalised polylogarithms were considered in some detail by Kölbig in 1986 in [91a].

There also appears to be a relationship with the generalised Rogers' function $L_n(z)$ defined in Lewin's survey "Structural Properties of Polylogarithms" [101, p.31] by

(4.4.91k) $$Li_n(z) = L_n(z) + \sum_{r=1}^{n-2} \frac{\log^r x \, Li_{n-r}(z)}{r!} - \frac{1}{n!} \log^{n-1} z \log(1-z)$$

but I shall leave further analysis of this aspect to the experts.

With a change of variables $t = y/x$ we have

$$\int_0^x \frac{\log^{s-1} y \log(1-y)}{y} dy = \int_0^1 \frac{\log^{s-1}(xt) \log(1-xt)}{t} dt$$

$$= \int_0^1 \frac{(\log x + \log t)^{s-1} \log(1-xt)}{t} dt$$

Using the binomial theorem this becomes

$$= \sum_{k=0}^{s-1} \binom{s-1}{k} \log^k x \int_0^1 \frac{\log^{s-1-k} t \log(1-xt)}{t} dt$$

and hence we have a connection with Nielsen's generalised polylogarithms.

**(viii) Theorem 4.6(a):**

(4.4.91l) $$\sum_{n=0}^{\infty} \frac{w^n}{n+1} \sum_{k=0}^{n} \binom{n}{k} \frac{x^k}{(k+1)^{s-1}} = \frac{(-1)^{s-1}}{w \Gamma(s-1)} \int_0^1 \frac{\log^{s-2} v \log[1-w(1+xv)]}{1+xv} dv$$

**Proof:**

Generalising (4.4.88a) we have

$$\sum_{n=0}^{\infty} \frac{w^n S_n}{n+1} = \frac{1}{\Gamma(s-1)} \sum_{n=0}^{\infty} \int_0^{\infty} e^{-nt} e^{-t} dt \int_0^{\infty} u^{s-2} e^{-u} xw^n (1+xe^{-u})^n du$$

$$= \frac{x}{\Gamma(s-1)} \int_0^{\infty} e^{-t} dt \int_0^{\infty} u^{s-2} e^{-u} \frac{1}{1 - e^{-t} w(1+xe^{-u})} du$$

$$= -\frac{x}{\Gamma(s-1)} \int_0^{\infty} \frac{u^{s-2} e^{-u} \log[1-w(1+xe^{-u})]}{w(1+xe^{-u})} du$$

$$= -\frac{x}{w\Gamma(s-1)} \int_0^\infty \frac{u^{s-2} \log\left[1-w(1+xe^{-u})\right] du}{e^u + x}$$

With the substitution $v = e^{-u}$ this becomes (4.4.91l)

$$\sum_{n=0}^\infty \frac{w^n}{n+1} \sum_{k=0}^n \binom{n}{k} \frac{x^k}{(k+1)^{s-1}} = \frac{(-1)^{s-1} x}{w\Gamma(s-1)} \int_0^1 \frac{\log^{s-2} v \log[1-w(1+xv)]}{1+xv} dv$$

When $w = 1, x = -y$ we get

(4.4.91m) $\sum_{n=0}^\infty \frac{1}{n+1} \sum_{k=0}^n \binom{n}{k} \frac{(-1)^k y^k}{(k+1)^{s-1}} = \frac{(-1)^s}{\Gamma(s-1)} \int_0^1 \frac{\log^{s-2} v \log(yv)}{1-yv} dv$

$$= \frac{(-1)^s \log y}{\Gamma(s-1)} \int_0^1 \frac{\log^{s-2} v}{1-yv} dv + \frac{(-1)^s}{\Gamma(s-1)} \int_0^1 \frac{\log^{s-1} v}{1-yv} dv$$

Then referring to (4.4.64i)

$$Li_{s+1}(y) = \frac{(-1)^{s+1} y}{\Gamma(s+1)} \int_0^1 \frac{\log^s v}{1-yv} dv$$

we see that

(4.4.91n)

$$\frac{(-1)^s}{\Gamma(s-1)} \int_0^1 \frac{\log^{s-2} v \log(yv)}{1-yv} dv = \sum_{n=0}^\infty \frac{1}{n+1} \sum_{k=0}^n \binom{n}{k} \frac{(-1)^k y^k}{(k+1)^{s-1}} = \frac{\log y}{y} Li_{s-1}(y) + \frac{s-1}{y} Li_s(y)$$

and hence we have another proof of (4.4.85).

Integrating (4.4.91l) with respect to $w$ we get

(4.4.91o) $\sum_{n=0}^\infty \frac{w^{n+1}}{(n+1)^2} \sum_{k=0}^n \binom{n}{k} \frac{x^k}{(k+1)^{s-1}} = \frac{(-1)^s}{\Gamma(s-1)} \int_0^1 \frac{\log^{s-2} v\, Li_2[w(1+xv)]}{1+xv} dv$

We have

(4.4.91p) $\int_0^t \frac{Li_2[w(1+xv)]}{1+xv} dv = \frac{Li_3[w(1+xt)]}{x} - \frac{Li_3(w)}{x}$

and hence with $s = 2$ we obtain from (4.4.91o)

(4.4.91q) $$\sum_{n=0}^{\infty} \frac{w^{n+1}}{(n+1)^2} \sum_{k=0}^{n} \binom{n}{k} \frac{x^k}{k+1} = \int_0^1 \frac{Li_2[w(1+xv)]}{1+xv} dv = \frac{Li_3[w(1+x)]}{x} - \frac{Li_3(w)}{x}$$

Letting $x = 1$ in the above equation results in

(4.4.91r) $$\sum_{n=0}^{\infty} \frac{w^{n+1}}{(n+1)^2} \sum_{k=0}^{n} \binom{n}{k} \frac{1}{k+1} = \int_0^1 \frac{Li_2[w(1+v)]}{1+v} dv = Li_3(2w) - Li_3(w)$$

and putting $w = 1/2$ gives us

(4.4.91s) $$\sum_{n=0}^{\infty} \frac{1}{2^{n+1}(n+1)^2} \sum_{k=0}^{n} \binom{n}{k} \frac{1}{k+1} = \int_0^1 \frac{Li_2[(1+v)/2]}{1+v} dv = \varsigma(3) - Li_3(1/2)$$

Referring to (4.4.123a) in Volume IV

(4.4.91si) $$\sum_{k=0}^{n} \binom{n}{k} \frac{1}{k+1} = \frac{2^{n+1} - 1}{n+1}$$

we note that the series in (4.4.91r) simplifies to

$$\sum_{n=0}^{\infty} \frac{w^{n+1}}{(n+1)^2} \sum_{k=0}^{n} \binom{n}{k} \frac{1}{k+1} = \sum_{n=0}^{\infty} \frac{w^{n+1}}{(n+1)^2} \frac{2^{n+1} - 1}{n+1} = \sum_{n=0}^{\infty} \frac{(2w)^{n+1}}{(n+1)^3} - \sum_{n=0}^{\infty} \frac{w^{n+1}}{(n+1)^3}$$

and the rest of equation (4.4.91r) becomes obvious.

With $x = -1$ in (4.4.91l) we have

(4.4.91t) $$\sum_{n=0}^{\infty} \frac{w^n}{n+1} \sum_{k=0}^{n} \binom{n}{k} \frac{(-1)^k}{(k+1)^{s-1}} = \frac{(-1)^{s-1}}{w\Gamma(s-1)} \int_0^1 \frac{\log^{s-2} v \log[1 - w(1-v)]}{1-v} dv$$

and $s = 2, 3$ and $4$ respectively gives us in conjunction with (4.4.123), (4.4.127) and (4.4.130)

(4.4.91u) $$\sum_{n=0}^{\infty} \frac{w^n}{n+1} \sum_{k=0}^{n} \binom{n}{k} \frac{(-1)^k}{(k+1)} = \sum_{n=0}^{\infty} \frac{w^n}{(n+1)^2} = \frac{(-1)^{s-1}}{w\Gamma(s-1)} \int_0^1 \frac{\log[1 - w(1-v)]}{1-v} dv$$

(4.4.91v) $$\sum_{n=0}^{\infty} \frac{w^n}{n+1} \sum_{k=0}^{n} \binom{n}{k} \frac{(-1)^k}{(k+1)^2} = \sum_{n=0}^{\infty} \frac{H_{n+1}^{(1)}}{(n+1)^2} w^n = \frac{(-1)^{s-1}}{w\Gamma(s-1)} \int_0^1 \frac{\log v \log[1 - w(1-v)]}{1-v} dv$$

(4.4.91w)

$$\sum_{n=0}^{\infty}\frac{w^n}{n+1}\sum_{k=0}^{n}\binom{n}{k}\frac{(-1)^k}{(k+1)^3} = \frac{1}{2}\sum_{n=0}^{\infty}\frac{\left(H_{n+1}^{(1)}\right)^2+H_{n+1}^{(2)}}{(n+1)^3}w^n = \frac{(-1)^{s-1}}{w\Gamma(s-1)}\int_0^1\frac{\log^2 v\log[1-w(1-v)]}{1-v}dv$$

Carrying out the same operation on (4.4.91o) results in

(4.4.91x) $$\sum_{n=0}^{\infty}\frac{w^n}{(n+1)^2}\sum_{k=0}^{n}\binom{n}{k}\frac{(-1)^k}{(k+1)} = \sum_{n=0}^{\infty}\frac{w^n}{(n+1)^3} = \frac{(-1)^{s-1}}{w\Gamma(s-1)}\int_0^1\frac{Li_2[w(1-v)]}{1-v}dv$$

(4.4.91y)

$$\sum_{n=0}^{\infty}\frac{w^n}{(n+1)^2}\sum_{k=0}^{n}\binom{n}{k}\frac{(-1)^k}{(k+1)^2} = \sum_{n=0}^{\infty}\frac{H_{n+1}^{(1)}}{(n+1)^3}w^n = \frac{(-1)^{s-1}}{w\Gamma(s-1)}\int_0^1\frac{\log v\, Li_2[w(1-v)]}{1-v}dv$$

(4.4.91z)

$$\sum_{n=0}^{\infty}\frac{w^n}{(n+1)^2}\sum_{k=0}^{n}\binom{n}{k}\frac{(-1)^k}{(k+1)^3} = \frac{1}{2}\sum_{n=0}^{\infty}\frac{\left(H_{n+1}^{(1)}\right)^2+H_{n+1}^{(2)}}{(n+1)^4}w^n = \frac{(-1)^{s-1}}{w\Gamma(s-1)}\int_0^1\frac{\log^2 v\, Li_2[w(1-v)]}{1-v}dv$$

We now multiply (4.4.91l) by $w$ and differentiate to obtain

(4.4.91zi) $$\sum_{n=0}^{\infty}w^n\sum_{k=0}^{n}\binom{n}{k}\frac{x^k}{(k+1)^{s-1}} = \frac{(-1)^s}{\Gamma(s-1)}\int_0^1\frac{\log^{s-2}v}{1-w(1+xv)}dv$$

and with $w=1/2$ we get

(4.4.91zii) $$\sum_{n=0}^{\infty}\frac{1}{2^{n+1}}\sum_{k=0}^{n}\binom{n}{k}\frac{x^k}{(k+1)^{s-1}} = \frac{(-1)^s}{\Gamma(s-1)}\int_0^1\frac{\log^{s-2}v}{1-xv}dv$$

With $s \to s+1$ this may be written as

(4.4.91ziii) $$\sum_{n=0}^{\infty}\frac{1}{2^{n+1}}\sum_{k=0}^{n}\binom{n}{k}\frac{x^k}{(k+1)^s} = \frac{(-1)^{s+1}}{\Gamma(s)}\int_0^1\frac{\log^{s-1}v}{1-xv}dv$$

**(ix)** Sondow [123] has proved the interesting result

(4.4.92) $$\gamma = \int_1^\infty \sum_{n=1}^\infty \frac{n!}{(n+1)x(x+1)...(x+n)} dx = \int_1^\infty \sum_{n=1}^\infty \frac{g(x)}{n+1} dx$$

$$= \int_1^\infty \sum_{n=1}^\infty \frac{\Gamma(n+1)\Gamma(x)}{(n+1)\Gamma(n+x+1)} dx$$

$$= \int_1^\infty \sum_{n=1}^\infty \frac{B(n+1,x)}{(n+1)} dx$$

(4.4.92a) $$= \sum_{n=1}^\infty \frac{1}{n+1} \sum_{k=0}^n \binom{n}{k} (-1)^{k+1} \log(k+1)$$

where $g(x)$ is defined in (4.2.1) and $\gamma$ is Euler's constant

(4.4.92b) $$\gamma = \lim_{n \to \infty} (H_n - \log n)$$

An alternative proof of (4.4.92) is summarised below. The following identity was given by Anglesio in [9].

(4.4.93) $$\int_0^\infty \frac{e^{-ax}(1-e^{-x})^n}{x^r} dx = \frac{(-1)^r}{(r-1)!} \sum_{k=0}^n \binom{n}{k} (-1)^k (a+k)^{r-1} \log(a+k)$$

where $a \geq 0$ and $1 \leq r \leq n$ (except for $a = 0, r = 1$).

In fact, Anglesio's proof contains the identity for $a, p > 0$

$$\int_0^\infty e^{-ax}(1-e^{-x})^n x^{p-1} dx = \Gamma(p) \sum_{k=0}^n \binom{n}{k} \frac{(-1)^k}{(a+k)^p}$$

As $p \to 0$ the right-hand side is indeterminate but we may use L'Hôpital's rule when considering the ratio (since $\sum_{k=0}^n \binom{n}{k} (-1)^k = \delta_{n,0}$ we require $n \neq 0$ to ensure that the numerator is equal to zero as $p \to 0$)

$$\frac{\sum_{k=0}^{n}\binom{n}{k}\frac{(-1)^k}{(a+k)^p}}{1/\Gamma(p)}$$

We have using Euler's reflection formula

$$\frac{d}{dp}[1/\Gamma(p)] = \frac{1}{\pi}\frac{d}{dp}[\Gamma(1-p)\sin\pi p]$$

$$= \Gamma(1-p)\cos\pi p - \frac{1}{\pi}\Gamma'(1-p)\sin\pi p$$

and therefore we see that where $r$ is an integer

$$\lim_{p\to 1-r}\frac{d}{dp}[1/\Gamma(p)] = (-1)^{r+1}\Gamma(r)$$

We now let $a=1$ and $r=1$ to obtain

(4.4.94) $$\int_0^\infty \frac{e^{-x}(1-e^{-x})^n}{x}\,dx = \sum_{k=0}^{n}\binom{n}{k}(-1)^{k+1}\log(1+k) = f(n)$$

We may also note that for $p < n$ we have

$$\sum_{k=0}^{n}\binom{n}{k}(-1)^k(a+k)^p = 0$$

Making the substitution $y = e^{-x}$ (4.4.94) becomes

(4.4.94a) $$\int_0^1 \frac{(1-y)^n}{\log y}\,dy = \sum_{k=0}^{n}\binom{n}{k}(-1)^k\log(1+k)$$

and this is in agreement with G&R [74, p.541, Eqn. 4.267 1] as corrected by Boros and Moll (see the Errata to G&R [74] dated 26 April 2005).

Equation (4.4.94a) may be derived more directly by considering (4.4.24c)

$$\int_0^1 \frac{(t^{\beta-1} - t^{\alpha-1})}{\log t}\,dt = \log\frac{\beta}{\alpha}$$

and hence we get with $\alpha = 1$ (see also (4.4.99k))

(4.4.94b) $$\int_0^1 \frac{(t^k-1)}{\log t} dt = \log(1+k)$$

Using the binomial theorem we have

$$\int_0^1 \frac{(1-y)^n}{\log y} dy = \int_0^1 \sum_{k=0}^n \binom{n}{k}(-1)^k y^k \frac{dy}{\log y}$$

$$= \int_0^1 \sum_{k=0}^n \binom{n}{k} \frac{(-1)^k [(y^k-1)+1]}{\log y} dy$$

A simple application of the binomial theorem then shows that $\sum_{k=0}^n \binom{n}{k}(-1)^k = 0$ provided $n \neq 0$ and hence for $n \geq 1$ we get

$$\int_0^1 \frac{(1-y)^n}{\log y} dy = \sum_{k=0}^n \binom{n}{k}(-1)^k \int_0^1 \frac{(y^k-1)}{\log y} dy$$

Reference to (4.4.94b) then results in (4.4.94a).

Let us now return to the proof of Sondow's formula.

Applying the integral identity

(4.4.95) $$\frac{1}{n+1} = \int_0^\infty e^{-nt} e^{-t} dt$$

we have

(4.4.96) $$\sum_{n=1}^\infty \frac{f(n)}{n+1} = \sum_{n=1}^\infty \int_0^\infty e^{-nt} e^{-t} dt \int_0^\infty \frac{e^{-x}(1-e^{-x})^n}{x} dx$$

Using the geometric series we have

$$\sum_{n=1}^\infty \left(e^{-t}(1-e^{-x})\right)^n = \frac{e^{-t}(1-e^{-x})}{1-e^{-t}(1-e^{-x})} = \frac{Ae^{-t}}{1-Ae^{-t}}$$

and we may therefore write the right hand side of (4.4.96) as
44

(4.4.97) $$I = \int_0^\infty \frac{e^{-x}}{x} dx \int_0^\infty \frac{Ae^{-2t}}{1-Ae^{-t}} dt$$

With the substitution $y = e^{-t}$ we obtain

$$\int \frac{Ae^{-2t}}{1-Ae^{-t}} dt = -\int \frac{Ay}{1-Ay} dy = \int \left(1 - \frac{1}{1-Ay}\right) dy$$

$$= y + \frac{1}{A}\log(1-Ay)$$

$$= e^{-t} + \frac{1}{A}\log(1-Ae^{-t})$$

Therefore we have

(4.4.98) $$\int_0^\infty \frac{Ae^{-2t}}{1-Ae^{-t}} dt = -1 - \frac{1}{A}\log(1-A) = -1 + \frac{x}{1-e^{-x}}$$

Substituting this in (4.4.97) we obtain

(4.4.99) $$I = \int_0^\infty \left(\frac{1}{1-e^{-x}} - \frac{1}{x}\right) e^{-x} dx$$

$$= \int_0^1 \left[\frac{1}{\log y} + \frac{1}{1-y}\right] dy$$

and this is one of the many integral forms of Euler's constant $\gamma$ (see for example [25, p.177]: a formal proof is also contained in Appendix E of Volume VI ). Accordingly, we have an alternative proof of Sondow's formula for $\gamma$.

$$\gamma = \sum_{n=1}^\infty \frac{1}{n+1} \sum_{k=0}^n \binom{n}{k} (-1)^{k+1} \log(1+k)$$

Making the substitution $x = \alpha t$ in (4.4.99) we get

$$I = \alpha \int_0^\infty \left(\frac{1}{e^{\alpha t}-1} - \frac{1}{\alpha t e^{\alpha t}}\right) dt$$

Differentiation with respect to $\alpha$ results in

$$\frac{dI}{d\alpha} = \int_0^\infty \left( \frac{1}{e^{\alpha t}-1} - \frac{1}{\alpha t e^{\alpha t}} \right) dt + \alpha \int_0^\infty \left( -\frac{t e^{\alpha t}}{(e^{\alpha t}-1)^2} + \frac{1+\alpha t}{\alpha^2 t e^{\alpha t}} \right) dt$$

Letting $\alpha = 1$ and noting that $\dfrac{dI}{d\alpha} = 0$ we obtain a (new?) integral identity for $\gamma$

(4.4.99a) $\qquad \gamma + 1 = \int_0^\infty \left( \frac{t e^t}{(e^t - 1)^2} - \frac{1}{t e^t} \right) dt = -\int_0^1 \left[ \frac{y \log y}{(1-y)^2} + \frac{y}{\log y} \right] dy$

Other integrals may be evaluated by considering $\dfrac{d^p I}{d\alpha^p} = 0$ etc.

Sondow's other representation for $\gamma$ is

$$\gamma = \int_1^\infty \sum_{n=1}^\infty \frac{n!}{(n+1)x(x+1)\ldots(x+n)}\, dx$$

and using (4.4.13) this is equivalent to

$$= \int_1^\infty \sum_{n=1}^\infty \frac{1}{n+1} \int_0^1 t^{x-1}(1-t)^n \, dx\, dt$$

Now using $\dfrac{1}{n+1} = \int_0^1 y^n dy$ we may write this as

$$= \int_1^\infty \sum_{n=1}^\infty \int_0^1 y^n dy \int_0^1 t^{x-1}(1-t)^n \, dx\, dt$$

$$= \int_1^\infty dx \int_0^1 dy \int_0^1 \frac{t^{x-1} y(1-t)}{1 - y(1-t)} dt$$

$$= \int_1^\infty dx \int_0^1 dy \int_0^1 t^{x-1}\left( \frac{1}{1 - y(1-t)} - 1 \right) dt$$

$$= \int_1^\infty dx \int_0^1 t^{x-1}\left( -\frac{\log t}{1-t} - 1 \right) dt$$

Sondow then reverses the order of integration (the integrand being non-negative) and obtains the integral (4.4.99) for $\gamma$.

In (4.3.48) of Volume II(a) we showed that

$$\psi^{(p)}(x) = \sum_{k=1}^{\infty} \frac{(-1)^{k+1}}{k} \frac{p!\, s(k,p)}{x(x+1)...(x+k-1)}$$

and with $p=1$ we have

$$\psi'(x) = \sum_{k=1}^{\infty} \frac{(-1)^{k+1}}{k} \frac{s(k,1)}{x(x+1)...(x+k-1)} = \sum_{k=1}^{\infty} \frac{k!}{k^2 x(x+1)...(x+k-1)}$$

where we have used $s(k,1) = (-1)^{k+1}(k-1)!$ in (3.105i). This is similar in structure to Sondow's representation for $\gamma$.

**Theorem:**

(4.4.99ai) $$\sum_{n=1}^{\infty} t^n \sum_{k=0}^{n} \binom{n}{k} \frac{x^k}{(k+y)^s} = \frac{1}{\Gamma(s)} \int_0^{\infty} \frac{u^{s-1} e^{-yu}\left(1+xe^{-u}\right)t}{\left[1-\left(1+xe^{-u}\right)t\right]} du$$

(4.4.99aii) $$\sum_{n=1}^{\infty} \frac{1}{n+1} \sum_{k=0}^{n} \binom{n}{k} \frac{(-1)^k}{(k+y)^s} = -\frac{1}{y^s} + \frac{(-1)^s}{\Gamma(s)} \int_0^1 \frac{(1-w)^{y-1} \log^s(1-w)}{w} dw$$

(4.4.99aiii) $$\sum_{n=1}^{\infty} \frac{1}{n+1} \sum_{k=0}^{n} \binom{n}{k} (-1)^k \log(k+z) = \psi(z+1) - \log z - \frac{1}{z}$$

**Proof:**

Referring back to the modus operandi employed in (4.4.44) we have for $s > 0$ and $k+y > 0$

$$\frac{1}{(k+y)^s} = \frac{1}{\Gamma(s)} \int_0^{\infty} u^{s-1} e^{-u(k+y)} du$$

We now consider the finite sum set out below (where the summation now starts at $k=0$ and we specify that $y$ is neither zero nor a negative integer)

$$S_n^0(x,y) = \sum_{k=0}^{n} \binom{n}{k} \frac{x^k}{(k+y)^s}$$

(where the superscript highlights the fact that the summation starts this time at $k = 0$)

Now combine the above two equations to obtain

$$S_n^0(x, y) = \sum_{k=0}^{n} \binom{n}{k} \frac{x^k}{(k+y)^s} = \sum_{k=0}^{n} \binom{n}{k} x^k \cdot \frac{1}{\Gamma(s)} \int_0^\infty u^{s-1} e^{-u(k+y)} du$$

$$= \frac{1}{\Gamma(s)} \int_0^\infty u^{s-1} \sum_{k=0}^{n} \left\{ \binom{n}{k} e^{-u(k+y)} x^k \right\} du$$

$$= \frac{1}{\Gamma(s)} \int_0^\infty u^{s-1} \sum_{k=0}^{n} \left\{ \binom{n}{k} \left[ e^{-u} x \right]^k \right\} e^{-yu} du$$

Applying the binomial theorem we have

$$S_n^0(x, y) = \frac{1}{\Gamma(s)} \int_0^\infty u^{s-1} \left(1 + x e^{-u}\right)^n e^{-yu} du$$

We now consider the following function $Y_s^0(x, y, t)$ which, for $|t| < 1$, is defined by the series

$$Y_s^0(x, y, t) = \sum_{n=1}^{\infty} t^n \sum_{k=0}^{n} \binom{n}{k} \frac{x^k}{(k+y)^s}$$

(where we note that the summation starts this time at $n = 1$). We then substitute the relation obtained for $S_n^0(x, y)$ to deduce

$$Y_s^0(x, y, t) = \sum_{n=1}^{\infty} t^n \frac{1}{\Gamma(s)} \int_0^\infty u^{s-1} \left(1 + x e^{-u}\right)^n e^{-yu} du$$

$$= \frac{1}{\Gamma(s)} \int_0^\infty u^{s-1} \sum_{n=1}^{\infty} \left(1 + x e^{-u}\right)^n t^n e^{-yu} du$$

(where we have again assumed that interchanging the order of summation and integration is permitted). Use of the geometric series simplifies this to

(4.4.99aiv) $\quad Y_s^0(x, y, t) = \sum_{n=1}^{\infty} t^n \sum_{k=0}^{n} \binom{n}{k} \frac{x^k}{(k+y)^s} = \frac{1}{\Gamma(s)} \int_0^\infty \frac{u^{s-1} e^{-yu} \left(1 + x e^{-u}\right) t}{\left[1 - \left(1 + x e^{-u}\right) t\right]} du$

We have for $x = -1$

(4.4.99av) $$Y_s^0(-1,y,t) = \sum_{n=1}^{\infty} t^n \sum_{k=0}^{n} \binom{n}{k} \frac{(-1)^k}{(k+y)^s} = \frac{1}{\Gamma(s)} \int_0^{\infty} \frac{u^{s-1} e^{-yu}(1-e^{-u})t}{\left[1-(1-e^{-u})t\right]} du$$

Also, we have for $s=1$ and $x=-1$

$$Y_1^0(-1,y,t) = \sum_{n=1}^{\infty} t^n \sum_{k=0}^{n} \binom{n}{k} \frac{(-1)^k}{k+y} = \int_0^{\infty} \frac{e^{-yu}(1-e^{-u})t}{\left[1-(1-e^{-u})t\right]} du$$

Letting $t = 1/2$ we get

$$\sum_{n=1}^{\infty} \frac{1}{2^n} \sum_{k=0}^{n} \binom{n}{k} \frac{(-1)^k}{(k+y)^s} = \frac{1}{\Gamma(s)} \int_0^{\infty} \frac{u^{s-1} e^{-yu}(1-e^{-u})}{1+e^{-u}} du$$

$$= \frac{1}{\Gamma(s)} \int_0^{\infty} \frac{u^{s-1}}{e^u+1}\left[e^{-(y-1)u} - e^{-yu}\right] du$$

and using (4.4.44fi) this is equal to

$$= \Phi(-1,s,y) - \Phi(-1,s,y+1)$$

Then (4.4.82) gives us

$$= \sum_{n=0}^{\infty} \frac{(-1)^n}{(n+y)^s} - \sum_{n=0}^{\infty} \frac{(-1)^n}{(n+y+1)^s}$$

and reference to (4.4.24) gives us

$$= \sum_{n=0}^{\infty} \frac{(-1)^n}{(n+y)^s}$$

We then obtain another derivation of (4.4.24a)

$$\sum_{n=0}^{\infty} \frac{1}{2^{n+1}} \sum_{k=0}^{n} \binom{n}{k} \frac{(-1)^k}{(k+y)^s} = \sum_{n=0}^{\infty} \frac{(-1)^n}{(n+y)^s}$$

Integrating (4.4.99av) with respect to $t$ gives us

$$\sum_{n=1}^{\infty}\frac{v^{n+1}}{n+1}\sum_{k=0}^{n}\binom{n}{k}\frac{(-1)^k}{(k+y)^s} = \frac{1}{\Gamma(s)}\int_0^v dt\int_0^{\infty}\frac{u^{s-1}e^{-yu}(1-e^{-u})t}{[1-(1-e^{-u})t]}du$$

We have the elementary integral

$$\int_0^v \frac{(1-e^{-u})t}{[1-(1-e^{-u})t]}dt = \int_0^v\left(-1+\frac{1}{1-e^{-u}}\frac{(1-e^{-u})t}{[1-(1-e^{-u})t]}\right)dt$$

$$= -t - \frac{1}{1-e^{-u}}\log[1-(1-e^{-u})t]\Big|_0^v$$

$$= -v - \frac{1}{1-e^{-u}}\log[1-(1-e^{-u})v]$$

and therefore we have

$$\sum_{n=1}^{\infty}\frac{v^{n+1}}{n+1}\sum_{k=0}^{n}\binom{n}{k}\frac{(-1)^k}{(k+y)^s} = \frac{1}{\Gamma(s)}\int_0^{\infty}u^{s-1}e^{-yu}\left(-v-\frac{1}{1-e^{-u}}\log[1-(1-e^{-u})v]\right)du$$

Using the integral definition of the gamma function this becomes

(4.4.99avi)
$$= -\frac{v}{y^s} - \frac{1}{\Gamma(s)}\int_0^{\infty}\frac{u^{s-1}e^{-yu}}{1-e^{-u}}\log[1-(1-e^{-u})v]du$$

Making the substitution $w = 1-e^{-u}$ this becomes

$$= -\frac{v}{y^s} + \frac{(-1)^s}{\Gamma(s)}\int_0^1\frac{(1-w)^{y-1}\log^{s-1}(1-w)\log(1-wv)}{w}dw$$

With $v = 1$ we have

$$\sum_{n=1}^{\infty}\frac{1}{n+1}\sum_{k=0}^{n}\binom{n}{k}\frac{(-1)^k}{(k+y)^s} = -\frac{1}{y^s} + \frac{(-1)^s}{\Gamma(s)}\int_0^1\frac{(1-w)^{y-1}\log^s(1-w)}{w}dw$$

$$= -\frac{1}{y^s} + \frac{(-1)^s}{\Gamma(s)}\frac{d^s}{dy^s}\int_0^1\frac{(1-w)^{y-1}}{w}dw$$

Therefore we have (the interchange of operations may not be valid here)

$$\sum_{n=0}^{\infty}\frac{1}{n+1}\sum_{k=0}^{n}\binom{n}{k}\frac{(-1)^k}{(k+y)^s} = \frac{(-1)^s}{\Gamma(s)}\int_0^1 \frac{(1-w)^{y-1}\log^s(1-w)}{w}dw = \frac{(-1)^s}{\Gamma(s)}\frac{d^s}{dy^s}\int_0^1 \frac{(1-w)^{y-1}}{w}dw$$

The Wolfram Integrator evaluates the last integral in terms of the Gauss hypergeometric function (convergent?)

$$\int \frac{(1-w)^{y-1}}{w}dw = -\frac{1}{y}(1-w)^y\left[1-\left(\frac{w-1}{w}\right)^{-y}{}_2F_1(-y,-y;1-y;1/w)\right]$$

Let us now try something rather less formidable! In (4.4.99avi) let $s = v = 1$ to obtain

$$\sum_{n=1}^{\infty}\frac{1}{n+1}\sum_{k=0}^{n}\binom{n}{k}\frac{(-1)^k}{k+y} = -\frac{1}{y} + \int_0^\infty \frac{ue^{-yu}}{1-e^{-u}}du = -\frac{1}{y} + \int_0^\infty \frac{ue^{-(y-1)u}}{e^u-1}du$$

$$= -\frac{1}{y} + \Phi(1,2,y)$$

We now integrate this with respect to $y$ and get

$$\sum_{n=1}^{\infty}\frac{1}{n+1}\sum_{k=0}^{n}\binom{n}{k}(-1)^k\log\frac{k+z}{k+1} = -\log z + \int_0^\infty \frac{u}{1-e^{-u}}\frac{\left(e^{-u}-e^{-zu}\right)}{u}du$$

We see that

$$\int_0^\infty \frac{u}{1-e^{-u}}\frac{\left(e^{-u}-e^{-zu}\right)}{u}du = \int_0^\infty \frac{e^u}{e^u-1}\left(e^{-u}-e^{-zu}\right)du$$

and integration by parts, and two applications of L'Hôpital's rule, gives us

$$\int_0^\infty \frac{e^u}{e^u-1}\left(e^{-u}-e^{-zu}\right)du = \int_0^\infty \log(e^u-1)\left(e^{-u}-ze^{-zu}\right)du$$

With the substitution $q = e^{-u}$ this becomes

$$= \int_0^1 \log\left(\frac{1-q}{q}\right)(1-zq^{z-1})dq$$

(4.4.99avii) $$= \int_0^1 \log\left(\frac{1-q}{q}\right)dq - z\int_0^1 q^{z-1}\log(1-q)dq + z\int_0^1 q^{z-1}\log q\, dq$$

It is easily shown that $\int_0^1 \log\left(\frac{1-q}{q}\right)dq = 0$ and we now refer to the beta function

$$B(z,y) = \int_0^1 q^{z-1}(1-q)^{y-1}\, dq$$

where differentiation results in

$$\frac{\partial}{\partial y}B(z,y) = \int_0^1 q^{z-1}(1-q)^{y-1}\log(1-q)dq$$

Hence we obtain

$$\left.\frac{\partial}{\partial y}B(z,y)\right|_{y=1} = \int_0^1 q^{z-1}\log(1-q)dq$$

From (4.4.8) we know that

$$B(z,y) = \frac{\Gamma(z)\Gamma(y)}{\Gamma(z+y)}$$

and hence as we have seen previously

$$\frac{\partial}{\partial y}B(z,y) = B(z,y)[\psi(y)-\psi(z+y)]$$

Therefore we obtain the well-known result

$$\left.\frac{\partial}{\partial y}B(z,y)\right|_{y=1} = B(z,1)[\psi(1)-\psi(z+1)] = \frac{1}{z}[\psi(1)-\psi(z+1)]$$

Hence we have

$$\int_0^1 q^{z-1}\log(1-q)dq = \frac{1}{z}[\psi(1)-\psi(z+1)]$$

Since $B(z,1) = \int_0^1 q^{z-1} \, dq = \dfrac{1}{z}$ we have upon differentiation

$$\int_0^1 q^{z-1} \log q \, dq = -\frac{1}{z^2}$$

Therefore (4.4.99avii) becomes

$$= -[\psi(1) - \psi(z+1)] - \frac{1}{z}$$

and we then deduce that

(4.4.99aviii) $\displaystyle\sum_{n=1}^{\infty} \frac{1}{n+1} \sum_{k=0}^{n} \binom{n}{k} (-1)^k \log\frac{k+z}{k+1} = -\log z - [\psi(1) - \psi(z+1)] - \frac{1}{z}$

At this stage we note Sondow's result (4.4.92a)

$$\gamma = -\sum_{n=1}^{\infty} \frac{1}{n+1} \sum_{k=0}^{n} \binom{n}{k} (-1)^k \log(k+1)$$

and then obtain for $z > 0$

$$\sum_{n=1}^{\infty} \frac{1}{n+1} \sum_{k=0}^{n} \binom{n}{k} (-1)^k \log(k+z) = \psi(z+1) - \log z - \frac{1}{z}$$

Since $\psi(z+1) = \psi(z) + \dfrac{1}{z}$ we may write this as

(4.4.99aix) $\displaystyle\sum_{n=1}^{\infty} \frac{1}{n+1} \sum_{k=0}^{n} \binom{n}{k} (-1)^k \log(k+z) = \psi(z) - \log z$

Starting the summation at $n = 0$ we see that

(4.4.99ax) $\displaystyle\sum_{n=0}^{\infty} \frac{1}{n+1} \sum_{k=0}^{n} \binom{n}{k} (-1)^k \log(k+z) = \psi(z)$

which we have seen before in (4.4.74) in Volume II(a).

My first attempt in this area is worth reproducing. It proceeded in the following manner:

Basic integration results in

$$\int_1^z \int_0^1 Y_1^0(-1,y,t)\,dt\,dy = \sum_{n=1}^{\infty} \frac{1}{n+1} \sum_{k=0}^{n} \binom{n}{k}(-1)^k \log\frac{k+z}{k+1} = \frac{1}{2}\int_1^z \int_0^1 \int_0^{\infty} \frac{e^{-yu}(1-e^{-u})t}{\left[1-(1-e^{-u})t\right]}\,du\,dt\,dy$$

We have

$$\int_0^1 \frac{(1-e^{-u})t}{\left[1-(1-e^{-u})t\right]}\,dt = \frac{1-e^{-u}-u}{1-e^{-u}} = 1 - \frac{ue^u}{e^u-1}$$

and hence

$$\int_1^z \int_0^1 \int_0^{\infty} \frac{e^{-yu}(1-e^{-u})t}{\left[1-(1-e^{-u})t\right]}\,du\,dt\,dy = \int_0^{\infty}\left(1 - \frac{ue^u}{e^u-1}\right)\frac{e^{-u}-e^{-zu}}{u}\,du$$

$$= \int_0^{\infty} \frac{e^{-u}-e^{-zu}}{u}\,du - \int_0^{\infty}\left(\frac{e^u}{e^u-1}\right)(e^{-u}-e^{-zu})\,du$$

The first integral is Frullani's integral (see (E.22aa) in Volume VI) and is equal to $\log z$. The second integral is

$$\int_0^{\infty}\left(\frac{e^u}{e^u-1}\right)(e^{-u}-e^{-zu})\,du = \int_0^{\infty} \frac{1-e^{-(z-1)u}}{e^u-1}\,du$$

and the Wolfram Integrator reports the indefinite integral in terms of the hypergeometric function.

$$\int \frac{1-e^{-(z-1)u}}{1-e^u}\,du = \frac{e^{-(z-1)u}}{z-1}\,{}_2F_1(1,1-z;2-z;e^u) - \log(1-e^{-u})$$

This hypergeometric function is however only convergent for $e^u < 1$ and hence I have not been able to pursue this particular line of approach any further. However, as seen above, a "human" integration produces an interesting result.

Ser [119b] produced the following analysis in a short paper published in 1926. He wrote

(4.4.99b)    $$\frac{1}{(1+x)^s} = S_0 - S_1 x - S_2 \frac{x(1-x)}{2!} - S_3 \frac{x(1-x)(2-x)}{3!} - \ldots$$

Letting $x=0$ we see that $S_0 = 1$. Similarly, letting $x=1$ we see that $S_1 = 1 - \frac{1}{2^s}$ and, more generally, Ser notes that (?)

(4.4.99c) $$S_n = 1 - \binom{1}{n}\frac{1}{2^s} + \binom{2}{n}\frac{1}{3^s} - \ldots + (-1)^n \frac{1}{n^s}$$

where the summation is over the upper index of the binomial coefficient.

Furthermore, we have

$$\varsigma(s) = \sum_{n=1}^{\infty} \frac{1}{(n+x)^s} + x + S_1 \frac{x(1-x)}{2!} + S_2 \frac{x(1-x)(2-x)}{3!} + \ldots$$

Then, integrating over the interval $[0,1]$ with respect to $x$, we obtain

$$\varsigma(s) = \frac{1}{s-1} + p_2 + S_1 p_3 + S_2 p_4 + \ldots$$

where

$$p_{n+1} = \frac{1}{n!} \int_0^1 x(1-x)\ldots(n-1-x)\,dx$$

Since [135, p.271] $\lim_{s \to 1}\left[\varsigma(s) - \frac{1}{s-1}\right] = \gamma$ (and also see (4.4.99m)), Ser concludes that

$$p_2 + \frac{p_3}{2} + \ldots = \gamma$$

From the February 2006 edition of The American Mathematical Monthly (p.168) I noted that

$$\int_0^1 x(x+1)\ldots(x+n-1)\,dx = (-1)^n B_n^{(n)}$$

where $B_n^{(n)}$ is the $n$th Bernoulli number of order $n$. The exponential generating function for these generalised Bernoulli numbers is

$$\sum_{n=0}^{\infty} (-1)^n B_n^{(n)} \frac{t^n}{n!} = -\frac{t}{(1-t)\log(1-t)}$$

We also have the generating function [126, p.61]

$$\left(\frac{t}{e^t-1}\right)^k = \sum_{n=0}^{\infty} B_n^{(k)} \frac{t^n}{n!} \quad \text{for } |t| < 2\pi$$

and therefore

$$(-1)^k \left(\frac{t}{e^{-t}-1}\right)^k = \sum_{n=0}^{\infty} (-1)^n B_n^{(k)} \frac{t^n}{n!}$$

We note from (3.104) that the Sterling numbers of the first kind may be defined by the generating function

$$x(x-1)\ldots(x-n+1) = \sum_{k=0}^{n} s(n,k) x^k$$

or equivalently

$$x(1-x)\ldots(n-1-x) = (-1)^{n-1} \sum_{k=0}^{n} s(n,k) x^k$$

Therefore we have

$$p_{n+1} = \frac{1}{n!} \int_0^1 x(1-x)\ldots(n-1-x) \, dx = \frac{1}{n!} \int_0^1 (-1)^{n-1} \sum_{k=0}^{n} s(n,k) x^k \, dx$$

$$= \frac{(-1)^{n-1}}{n!} \sum_{k=0}^{n} \frac{s(n,k)}{k+1}$$

I subsequently ascertained that Sondow's formula [123] for $e^\gamma$, i.e. the exponentiated form of (4.4.92a), was also reported by Ser [119b] in 1926 (this paper is available on the Comptes Rendus website mentioned in the References).

According to [114aa] Ser also showed that

$$\gamma = \frac{1}{n} \sum_{m=0}^{\infty} \binom{m+n}{m}^{-1} t_m + H_{n-1}^{(1)} - \log n$$

where $t_m = \dfrac{1}{(m+1)!} \int_0^1 x(1-x)\ldots(m-x) \, dx$. In this regard see also Yingying's paper [139a].

From (4.4.15) we have

(4.4.99d) $$\int_0^1 t^{x-1} (1-t)^n \log^{s-1} t \, dt = (-1)^{s-1} (s-1)! \sum_{k=0}^{n} \binom{n}{k} \frac{(-1)^k}{(k+x)^s}$$

Making the substitution $t = e^{-u}$ in (4.4.99d) we obtain

(4.4.99e) $$\int_0^1 t^{x-1}(1-t)^n \log^{s-1} t \, dt = (-1)^{s-1}\int_0^\infty e^{-xu}\left(1-e^{-u}\right)^n u^{s-1} du$$

and deduce the identity

(4.4.99f) $$\int_0^\infty e^{-xu}\left(1-e^{-u}\right)^n u^{s-1} du = (s-1)!\sum_{k=0}^n \binom{n}{k}\frac{(-1)^k}{(k+x)^s}$$

Integrating (4.4.99f) with respect to $x$ we have

$$\int_a^y dx \int_0^\infty e^{-xu}\left(1-e^{-u}\right)^n u^{s-1} du = (s-1)!\int_a^y \sum_{k=0}^n \binom{n}{k}\frac{(-1)^k}{(k+x)^s} dx$$

and, provided $s \neq 1$, this results in

(4.4.99g) $$\frac{1}{(s-2)!}\int_0^\infty \frac{\left(e^{-au}-e^{-yu}\right)}{u}\left(1-e^{-u}\right)^n u^{s-1} du = \sum_{k=0}^n \binom{n}{k}\frac{(-1)^k}{(k+y)^{s+1}} - \sum_{k=0}^n \binom{n}{k}\frac{(-1)^k}{(k+a)^{s+1}}$$

For $s = 1$ we have

(4.4.99h)

$$\int_0^\infty \frac{\left(e^{-au}-e^{-yu}\right)}{u}\left(1-e^{-u}\right)^n du = \sum_{k=0}^n \binom{n}{k}(-1)^{k+1}\log(y+k) - \sum_{k=0}^n \binom{n}{k}(-1)^{k+1}\log(a+k)$$

Letting $y = 1$ in (4.4.99h) we have

(4.4.99i) $$\int_0^\infty \frac{\left(e^{-au}-e^{-u}\right)}{u}\left(1-e^{-u}\right)^n du = \sum_{k=0}^n \binom{n}{k}(-1)^{k+1}\log(1+k) - \sum_{k=0}^n \binom{n}{k}(-1)^{k+1}\log(a+k)$$

and this is equivalent to Anglesio's formula (4.4.93)

$$\int_0^\infty \frac{e^{-au}(1-e^{-u})^n}{u^r} du = \frac{(-r)^r}{(r-1)!}\sum_{k=0}^n \binom{n}{k}(-1)^k (a+k)^{r-1} \log(a+k)$$

where $a \geq 0$ and $1 \leq r \leq n$ (except for $a = 0, r = 1$). With $r = 1$ we have as above

(4.4.99j) $$\int_0^\infty \frac{e^{-au}(1-e^{-u})^n}{u}du = \sum_{k=0}^n \binom{n}{k}(-1)^{k+1}\log(a+k)$$

We may also note that

$$\int_0^1 x^{nt}dt = \int_0^1 (x^n)^t dt = \frac{(x^n)^t}{\log(x^n)}\bigg|_0^1 = \frac{x^n-1}{n\log x}$$

Therefore we have

$$\int_0^1 \frac{x^n-1}{\log x}dx = n\int_0^1 dx \int_0^1 x^{nt}dt = n\int_0^1 \frac{dt}{nt+1} = \log(n+1)$$

Hence by summation we obtain for $s > 1$

$$\sum_{n=1}^\infty \frac{1}{n^s}\int_0^1 \frac{x^n-1}{\log x}dx = \int_0^1 \frac{Li_s(x)-\varsigma(s)}{\log x}dx = \sum_{n=1}^\infty \frac{\log(n+1)}{n^s}$$

It is clear that the above analysis may be extended to show that

(4.4.99k) $$\int_0^1 \frac{x^\alpha-1}{\log x}dx = \log(\alpha+1)$$

In [114aa] Rivoal recently gave an elementary proof of the following lemma:

(4.4.99l) $$\int_0^1 x^n \Omega(x)dx = \gamma - [H_n - \log(n+1)]$$

where, for convenience $\Omega(x)$, is defined as $\Omega(x) = \frac{1}{1-x} + \frac{1}{\log x}$. Using L'Hôpital's theorem it is easily seen that $\Omega(x)$ is continuous on $[0,1]$. For $x \in (0,1)$ we also have by a straightforward integration

$$\Omega(x) = \int_0^1 \frac{1-x^t}{1-x}dt$$

We have

$$\int_0^1 \Omega(x)dx = \int_0^1 \left(\frac{1}{1-x} + \frac{1}{\log x}\right)dx = \int_0^\infty \left[\frac{1}{1-e^{-u}} - \frac{1}{u}\right]e^{-u}du = \gamma$$

From (4.4.38) we have

$$\varsigma(s)\Gamma(s) = \int_0^\infty \frac{u^{s-1}}{e^u - 1} du$$

and we also see that $\Gamma(s) = (s-1)\Gamma(s-1) = (s-1)\int_0^\infty e^{-u} u^{s-2} du$. Therefore we get

$$\left[\varsigma(s) - \frac{1}{s-1}\right]\Gamma(s) = \int_0^\infty u^{s-1} \left[\frac{1}{e^u - 1} - \frac{1}{ue^u}\right] du$$

Hence, in the limit as $s \to 1$, we obtain

(4.4.99m) $\quad \lim_{s \to 1}\left[\varsigma(s) - \frac{1}{s-1}\right]\Gamma(s) = \lim_{s \to 1}\left[\varsigma(s) - \frac{1}{s-1}\right] = \int_0^\infty \left[\frac{1}{1-e^{-u}} - \frac{1}{u}\right] e^{-u} du = \gamma$

where, in the final part, we have employed (4.4.99). We now continue with Rivoal's lemma.

We have by simple algebra

$$\int_0^1 x^n \Omega(x) dx = \int_0^1 \Omega(x) dx - \int_0^1 \frac{x^n - 1}{x - 1} dx + \int_0^1 \frac{x^n - 1}{\log x} dx$$

$$= \gamma - H_n + \int_0^1 \frac{x^n - 1}{\log x} dx = \gamma - [H_n - \log(n+1)]$$

Hence by summation we obtain for $s > 1$ (see also (E.22bii) in Volume VI)

(4.4.99n) $\quad \gamma\varsigma(s) - \sum_{n=1}^\infty \frac{H_n}{n^s} + \sum_{n=1}^\infty \frac{\log(n+1)}{n^s} = \int_0^1 Li_s(x)\left(\frac{1}{1-x} + \frac{1}{\log x}\right) dx$

In 1997 Candelpergher et al. [38a] produced a somewhat similar result in the case where $s = 2$ by reference to "Ramanujan summation" but unfortunately I am not au fait with the underlying analysis contained in that paper.

Using the Georghiou and Philippou formula [69c] for the case where $s = m \in \mathbf{N}$

(4.4.99o) $\quad \sum_{n=1}^\infty \frac{H_n}{n^m} = \left(1 + \frac{m}{2}\right)\varsigma(m+1) - \frac{1}{2}\sum_{j=2}^{m-1}\varsigma(j)\varsigma(m-j+1)$

this may be written as

$$\gamma\varsigma(m) - \left(1 + \frac{m}{2}\right)\varsigma(m+1) + \frac{1}{2}\sum_{j=2}^{m-1}\varsigma(j)\varsigma(m-j+1) = \int_0^1 Li_m(x)\left(\frac{1}{1-x} + \frac{1}{\log x}\right)dx - \int_0^1 \frac{Li_m(x) - \varsigma(m)}{\log x}dx$$

We therefore have

(4.4.99p) $$\gamma\varsigma(m) - \left(1 + \frac{m}{2}\right)\varsigma(m+1) + \frac{1}{2}\sum_{j=2}^{m-1}\varsigma(j)\varsigma(m-j+1) = \int_0^1 \left[\frac{Li_m(x)}{1-x} + \frac{\varsigma(m)}{\log x}\right]dx$$

Rivoal's method may be slightly extended to show that

(4.4.99q) $$\int_0^1 x^\alpha \left(\frac{1}{1-x} + \frac{1}{\log x}\right)dx = -\psi(\alpha+1) + \log(\alpha+1)$$

and differentiating the above we get

$$\int_0^1 x^\alpha \log x \left(\frac{1}{1-x} + \frac{1}{\log x}\right)dx = -\psi'(\alpha+1) + \frac{1}{\alpha+1}$$

(in passing, it may be noted that both identities vanish as $\alpha \to \infty$).

Accordingly we have

$$\int_0^1 \frac{x^\alpha \log x}{1-x}dx = -\psi'(\alpha+1)$$

and this may be extended to

(4.4.99r) $$\int_0^1 \frac{x^\alpha \log^p x}{1-x}dx = -\psi^{(p)}(\alpha+1)$$

Integration of (4.4.99q) results in

$$\int_a^b d\alpha \int_0^1 x^\alpha \left(\frac{1}{1-x} + \frac{1}{\log x}\right)dx =$$

$$-\log\Gamma(b+1) + \log\Gamma(a+1) + a - b + (1+b)\log(b+1) - (1+a)\log(a+1)$$

Reversing the order of integration we obtain

$$\int_a^b d\alpha \int_0^1 x^\alpha \left( \frac{1}{1-x} + \frac{1}{\log x} \right) dx = \int_0^1 \frac{x^b - x^a}{\log x} \left( \frac{1}{1-x} + \frac{1}{\log x} \right) dx$$

and we therefore have

$$\int_0^1 \frac{x^b - x^a}{\log x} \left( \frac{1}{1-x} + \frac{1}{\log x} \right) dx =$$

$$-\log \Gamma(b+1) + \log \Gamma(a+1) + a - b + (1+b)\log(b+1) - (1+a)\log(a+1)$$

With $a = 0$, $b = 1$ we get

(4.4.99s) $\quad \displaystyle\int_0^1 \left( -\frac{1}{\log x} + \frac{x-1}{\log^2 x} \right) dx = 1 + 2\log 2$

Instead of the zeta summation in (4.4.99n) we may proceed as follows

$$\gamma \sum_{n=1}^\infty \frac{1}{2^n} - \sum_{n=1}^\infty \frac{H_n}{2^n} + \sum_{n=1}^\infty \frac{\log(n+1)}{2^n} = \int_0^1 \frac{x}{2-x}\left( \frac{1}{1-x} + \frac{1}{\log x} \right) dx$$

and using (3.7) of Volume I we have

(4.4.99t) $\quad \displaystyle \gamma - 2\log 2 + \sum_{n=1}^\infty \frac{\log(n+1)}{2^n} = \int_0^1 \frac{x}{2-x}\left( \frac{1}{1-x} + \frac{1}{\log x} \right) dx$

## EVALUATION OF VARIOUS LOGARITHMIC INTEGRALS

Lewin [100, p.20] showed with integration by parts that

$$\int_0^t x^\alpha Li_2(x)\, dx = \frac{t^{\alpha+1}}{\alpha+1} Li_2(t) + \int_0^t \frac{x^{\alpha+1}}{\alpha+1} \frac{\log(1-x)}{x} dx$$

$$\int_0^t \frac{x^\alpha}{\alpha+1} \log(1-x)\, dx = \frac{t^{\alpha+1}}{(\alpha+1)^2} \log(1-t) + \frac{1}{(\alpha+1)^2} \int_0^t \frac{x^{\alpha+1}}{1-x} dx$$

and noting that

$$\int_0^t \frac{x^{\alpha+1}}{1-x}\,dx = -\int_0^t \frac{1-x^{\alpha+1}}{1-x}\,dx + \log(1-t)$$

we may write this as

(4.4.100) $\quad \int_0^t x^\alpha Li_2(x)\,dx = \frac{t^{\alpha+1}}{\alpha+1} Li_2(t) + \frac{t^{\alpha+1}-1}{(\alpha+1)^2}\log(1-t) - \frac{1}{(\alpha+1)^2}\int_0^t \frac{1-x^{\alpha+1}}{1-x}\,dx$

With $t=1$ we have

(4.4.100a) $\quad \int_0^1 x^\alpha Li_2(x)\,dx = \frac{\varsigma(2)}{\alpha+1} - \frac{\psi(\alpha+2)+\gamma}{(\alpha+1)^2}$

where in the final part we have used

$$\psi(\alpha)+\gamma = \int_0^1 \frac{1-x^{\alpha-1}}{1-x}\,dx$$

The Wolfram Integrator provides the following result involving the Gauss hypergeometric function (which only converges for $|t|<1$)

$$\int_0^t \frac{1-x^{\alpha+1}}{1-x}\,dx = \frac{t^{\alpha+1}}{\alpha+1}\left[{}_2F_1(\alpha+1,1;\alpha+2;t)-1\right]$$

With $\alpha = n-1$ we have

(4.4.100b) $\quad \int_0^1 x^n \frac{Li_2(x)}{x}\,dx = \frac{\varsigma(2)}{n} - \frac{H_n}{n^2}$

Completing the summation we see that

$$\sum_{n=1}^\infty \int_0^1 (-1)^n x^n \frac{Li_2(x)}{x}\,dx = \int_0^1 \sum_{n=1}^\infty (-1)^n x^n \frac{Li_2(x)}{x}\,dx = \int_0^1 \frac{Li_2(x)}{(1+x)x}\,dx = \varsigma(2)\sum_{n=1}^\infty \frac{(-1)^n}{n} - \sum_{n=1}^\infty (-1)^n \frac{H_n}{n^2}$$

Hence we obtain

(4.4.100c) $\quad \int_0^1 \frac{Li_2(x)}{(1+x)x}\,dx = \sum_{n=1}^\infty (-1)^{n+1} \frac{H_n}{n^2} - \varsigma(2)\log 2$

Partial fractions give us

$$\int_0^1 \frac{Li_2(x)}{(1+x)x}dx = \int_0^1 \frac{Li_2(x)}{x}dx - \int_0^1 \frac{Li_2(x)}{(1+x)}dx$$

and integration by parts results in

$$\int_0^1 \frac{Li_2(x)}{(1+x)}dx = \log(1+x)Li_2(x)\Big|_0^1 + \int_0^1 \frac{\log(1+x)\log(1-x)}{x}dx$$

$$= \log 2\varsigma(2) + \int_0^1 \frac{\log(1+x)\log(1-x)}{x}dx$$

As shown in Lewin's monograph [100, p.159] we have

(4.4.100d) $$\int_0^1 \frac{\log(1+x)\log(1-x)}{x}dx = -\frac{5}{8}\varsigma(3)$$

and hence we obtain

(4.4.100e) $$\int_0^1 \frac{Li_2(x)}{(1+x)}dx = \log 2\varsigma(2) - \frac{5}{8}\varsigma(3)$$

$$\int_0^1 \frac{Li_2(x)}{(1+x)x}dx = \frac{13}{8}\varsigma(3) - \log 2\varsigma(2)$$

Therefore we obtain

(4.4.100ei) $$\sum_{n=1}^{\infty}(-1)^{n+1}\frac{H_n}{n^2} = \frac{13}{8}\varsigma(3)$$

See also (3.138) in Volume I which, unfortunately, is slightly different: your mission, if you are prepared to accept it, is to fix the arithmetic and determine the correct result!

Lewin's evaluation of various logarithmic integrals is shown below. We have

$$\int_0^t \frac{\log^2(1-x^2)}{x}dx = \int_0^t \frac{[\log(1+x)+\log(1-x)]^2}{x}dx$$

$$= \int_0^t \frac{\log^2(1+x)}{x}dx + \int_0^t \frac{\log^2(1-x)}{x}dx + 2\int_0^t \frac{\log(1+x)\log(1-x)}{x}dx$$

Integration by parts produces

$$\int_0^t \frac{\log^2(1+x)}{x}dx = \log t \log^2(1+t) - 2\int_0^t \frac{\log x \log(1+x)}{1+x}dx$$

With the substitution $y = 1+x$ we get

$$\int \frac{\log x \log(1+x)}{1+x}dx = \int \frac{\log(y-1)\log y}{y}dy = \int \frac{\log y}{y}[\log y + \log(1-1/y)]dy$$

$$= \frac{\log^3 y}{3} + \int \log y \frac{\log(1-1/y)}{y}dy$$

Integration by parts gives us

$$\int \log y \frac{\log(1-1/y)}{y}dy = \log y \, Li_2(1/y) - \int \frac{Li_2(1/y)}{y}dy$$

because $\int \frac{\log(1-1/y)}{y}dy = -\int \frac{\log(1-t)}{t}dt = Li_2(t) = Li_2(1/y)$

It is easily seen that

$$\int \frac{Li_2(1/y)}{y}dy = -\int \frac{Li_2(t)}{t}dt = -Li_3(t) = -Li_3(1/y)$$

and therefore we have

(4.4.100f)

$$\int_0^t \frac{\log x \log(1+x)}{1+x}dx = \frac{1}{3}\log^3(1+t) + \log(1+t)Li_2\left(\frac{1}{1+t}\right) + Li_3\left(\frac{1}{1+t}\right) - \varsigma(3)$$

Accordingly we get

(4.4.100g)

$$\int_0^t \frac{\log^2(1+x)}{x}dx = \log t \log^2(1+t) - \frac{2}{3}\log^3(1+t) - 2\log(1+t)Li_2\left(\frac{1}{1+t}\right) - 2Li_3\left(\frac{1}{1+t}\right) + 2\varsigma(3)$$

Similarly, for the other integral we obtain

$$\int_0^t \frac{\log^2(1-x)}{x}\,dx = \log t \log^2(1-t) + 2\int_0^t \frac{\log x \log(1-x)}{1-x}\,dx$$

With the substitution $y = 1-x$ we get

$$\int \frac{\log x \log(1-x)}{1-x}\,dx = -\int \frac{\log(1-y)}{y} \log y\,dy = \log y\, Li_2(y) - \int \frac{Li_2(y)}{y}\,dy$$

$$= \log y\, Li_2(y) - Li_3(y)$$

and hence we have

(4.4.100gi) $$\int_0^t \frac{\log x \log(1-x)}{1-x}\,dx = \log(1-t) Li_2(1-t) - Li_3(1-t) + \varsigma(3)$$

Hence we obtain

(4.4.100gii) $$\int_0^t \frac{\log^2(1-x)}{x}\,dx = \log t \log^2(1-t) + 2\log(1-t) Li_2(1-t) - 2Li_3(1-t) + 2\varsigma(3)$$

Making the substitution $x = y^2$ in the above integral we obtain

$$\int_0^{t^2} \frac{\log^2(1-x)}{x}\,dx = 2\int_0^t \frac{\log^2(1-y^2)}{y}\,dy$$

$$= \log t^2 \log^2(1-t^2) + 2\log(1-t^2) Li_2(1-t^2) - 2Li_3(1-t^2) + 2\varsigma(3)$$

and therefore we have

(4.4.100giii)

$$\int_0^t \frac{\log^2(1-x^2)}{x}\,dx = \log t \log^2(1-t^2) + \log(1-t^2) Li_2(1-t^2) - Li_3(1-t^2) + \varsigma(3)$$

We now return to our original identity

$$2\int_0^t \frac{\log(1+x)\log(1-x)}{x}\,dx = \int_0^t \frac{\log^2(1-x^2)}{x}\,dx - \int_0^t \frac{\log^2(1+x)}{x}\,dx - \int_0^t \frac{\log^2(1-x)}{x}\,dx$$

which then results in

(4.4.100giv) $$2\int_0^t \frac{\log(1+x)\log(1-x)}{x}dx$$

$$= \log t \log^2(1-t^2) + \log(1-t^2) Li_2(1-t^2) - Li_3(1-t^2) + \varsigma(3)$$

$$- \log t \log^2(1+t) + \frac{2}{3}\log^3(1+t) + 2\log(1+t) Li_2\left(\frac{1}{1+t}\right) + 2Li_3\left(\frac{1}{1+t}\right) - 2\varsigma(3)$$

$$- \log t \log^2(1-t) - 2\log(1-t) Li_2(1-t) + 2Li_3(1-t) - 2\varsigma(3)$$

In particular, as reported in de Doelder's paper [55], we have (as used above in (4.4.100d))

(4.4.100gv) $$\int_0^1 \frac{\log(1+x)\log(1-x)}{x}dx = -\frac{5}{8}\varsigma(3)$$

## SOME INTEGRALS INVOLVING POLYLOGARITHMS

Referring to (4.4.100b)

$$\int_0^1 x^n \frac{Li_2(x)}{x}dx = \frac{\varsigma(2)}{n} - \frac{H_n}{n^2}$$

we also have the summation

(4.4.100hi)

$$\sum_{n=1}^{\infty}\frac{1}{n}\int_0^1 x^n \frac{Li_2(x)}{x}dx = \int_0^1 \sum_{n=1}^{\infty}\frac{x^n}{n}\frac{Li_2(x)}{x}dx = -\int_0^1 \log(1-x)\frac{Li_2(x)}{x}dx = \varsigma^2(2) - \sum_{n=1}^{\infty}\frac{H_n}{n^3}$$

It is easily seen that

$$\int \log(1-x)\frac{Li_2(x)}{x}dx = -\frac{1}{2}[Li_2(x)]^2$$

and we then obtain the familiar result

$$\sum_{n=1}^{\infty}\frac{H_n}{n^3}=\frac{1}{2}\varsigma^2(2)=\frac{5}{4}\varsigma(4)$$

We also have the dilogarithm summation

(4.4.100hii) $\displaystyle\sum_{n=1}^{\infty}\frac{1}{n^2}\int_0^1 x^n\frac{Li_2(x)}{x}dx=\int_0^1\sum_{n=1}^{\infty}\frac{x^n}{n^2}\frac{Li_2(x)}{x}dx=\int_0^1\frac{[Li_2(x)]^2}{x}dx=\varsigma(2)\varsigma(3)-\sum_{n=1}^{\infty}\frac{H_n}{n^4}$

and similarly we obtain

(4.4.100hiii) $\displaystyle\int_0^1\frac{Li_3(x)Li_2(x)}{x}dx=\varsigma(2)\varsigma(4)-\sum_{n=1}^{\infty}\frac{H_n}{n^5}$

Since $\displaystyle\int\frac{Li_3(x)Li_2(x)}{x}dx=\frac{1}{2}[Li_3(x)]^2$ we see that

(4.4.100hiv) $\displaystyle\sum_{n=1}^{\infty}\frac{H_n}{n^5}=\varsigma(2)\varsigma(4)-\frac{1}{2}\varsigma^2(3)$

More generally we have

$$\sum_{n=1}^{\infty}\frac{1}{n^p}\int_0^1 x^n\frac{Li_2(x)}{x}dx=\int_0^1\sum_{n=1}^{\infty}\frac{x^n}{n^p}\frac{Li_2(x)}{x}dx=\int_0^1\frac{Li_p(x)Li_2(x)}{x}dx=\varsigma(2)\varsigma(p+1)-\sum_{n=1}^{\infty}\frac{H_n}{n^{p+2}}$$

A different integration by parts results in

$$\int_0^t x^n\frac{Li_2(x)}{x}dx=t^n Li_3(t)-n\int_0^t x^n\frac{Li_3(x)}{x}dx$$

and we therefore see that

(4.4.100i) $\displaystyle\int_0^1 x^n\frac{Li_3(x)}{x}dx=\frac{\varsigma(3)}{n}-\frac{\varsigma(2)}{n^2}+\frac{H_n}{n^3}$

Completing the summation we see that

$$\sum_{n=1}^{\infty}\int_0^1\frac{x^n}{n}\frac{Li_3(x)}{x}dx=\int_0^1\sum_{n=1}^{\infty}\frac{x^n}{n}\frac{Li_3(x)}{x}dx=-\int_0^1\log(1-x)\frac{Li_3(x)}{x}dx=\sum_{n=1}^{\infty}\frac{H_n}{n^4}$$

We also have

$$-\int_0^1 \frac{\log(1-x)}{x} Li_3(x)\,dx = Li_2(x)Li_3(x)\Big|_0^1 - \int_0^1 \frac{[Li_2(x)]^2}{x}\,dx$$

Freitas [69a] records the following integral

$$\int_0^1 \frac{[Li_2(x)]^2}{x}\,dx = 2\varsigma(2)\varsigma(3) - 3\varsigma(5)$$

and therefore we get (as noted by Flajolet and Salvy [69] )

(4.4.100ji) $$\sum_{n=1}^\infty \frac{H_n}{n^4} = 3\varsigma(5) - \varsigma(2)\varsigma(3)$$

Continuing in the same way we see that

$$\int_0^1 Li_4(x) \frac{Li_3(x)}{x}\,dx = \varsigma(3)\varsigma(5) - \varsigma(2)\varsigma(6) + \sum_{n=1}^\infty \frac{H_n}{n^7}$$

Since $\int \frac{Li_4(x)Li_3(x)}{x}\,dx = \frac{1}{2}[Li_4(x)]^2$ we conclude that

(4.4.100jii) $$\sum_{n=1}^\infty \frac{H_n}{n^7} = \varsigma(2)\varsigma(6) - \varsigma(3)\varsigma(5) - \frac{1}{2}\varsigma^2(4)$$

We also have an alternative summation of (4.4.100b)

(4.4.100jiv) $$\sum_{n=1}^\infty \frac{(-1)^n}{n} \int_0^1 x^n \frac{Li_2(x)}{x}\,dx = \int_0^1 \sum_{n=1}^\infty (-1)^{n+1} \frac{x^n}{n} \frac{Li_2(x)}{x}\,dx$$

$$= \int_0^1 \log(1+x) \frac{Li_2(x)}{x}\,dx = \varsigma(2)\log 2 - \sum_{n=1}^\infty (-1)^{n+1} \frac{H_n}{n^3}$$

The Wolfram Integrator was not able to evaluate the above integral.

More generally we have

$$\sum_{n=1}^\infty \frac{1}{n^p} \int_0^1 x^n \frac{Li_3(x)}{x}\,dx = \int_0^1 \sum_{n=1}^\infty \frac{x^n}{n^p} \frac{Li_3(x)}{x}\,dx$$

and therefore

$$\int_0^1 \frac{Li_p(x)Li_3(x)}{x}dx = \varsigma(3)\varsigma(p+1) - \varsigma(2)\varsigma(p+2) + \sum_{n=1}^{\infty}\frac{H_n}{n^{p+3}}$$

From (4.4.100i) we may also deduce

$$\int_0^1 x^n \frac{Li_3(x)}{x}dx = \varsigma(4) - n\int_0^1 x^n \frac{Li_4(x)}{x}dx$$

and accordingly we get

(4.4.100jv) $$\int_0^1 x^n \frac{Li_4(x)}{x}dx = \frac{\varsigma(4)}{n} - \frac{\varsigma(3)}{n^2} + \frac{\varsigma(2)}{n^3} - \frac{H_n}{n^4}$$

and therefore

$$\int_0^1 \frac{Li_p(x)Li_4(x)}{x}dx = \varsigma(4)\varsigma(p+1) - \varsigma(3)\varsigma(p+2) + \varsigma(2)\varsigma(p+3) - \sum_{n=1}^{\infty}\frac{H_n}{n^{p+4}}$$

Differentiating (4.4.100a)

$$\int_0^1 x^\alpha Li_2(x)\,dx = \frac{\varsigma(2)}{\alpha+1} - \frac{\psi(\alpha+2)+\gamma}{(\alpha+1)^2}$$

with respect to $\alpha$ we obtain

(4.4.100k) $$\int_0^1 x^\alpha \log x\, Li_2(x)\,dx = -\frac{\varsigma(2)}{(\alpha+1)^2} - \frac{\psi'(\alpha+2)}{(\alpha+1)^2} + 2\frac{\psi(\alpha+2)+\gamma}{(\alpha+1)^3}$$

With $\alpha = n-1$ we have

(4.4.100ki) $$\int_0^1 x^n \log x \frac{Li_2(x)}{x}dx = -\frac{\varsigma(2)}{n^2} - \frac{\psi'(n+1)}{n^2} + 2\frac{H_n}{n^3}$$

Using (E.16) in Volume VI we see that

$$\psi'(n+1) = \sum_{k=0}^{\infty}\frac{1}{(n+1+k)^2} = \varsigma(2) - H_n^{(2)}$$

and hence we obtain

(4.4.100kiii) $$\int_0^1 x^n \log x \frac{Li_2(x)}{x} dx = -2\frac{\varsigma(2)}{n^2} + \frac{H_n^{(2)}}{n^2} + 2\frac{H_n}{n^3}$$

Completing the summation we see that

(4.4.100kiv) $$\int_0^1 \frac{\log x \, Li_2(x)}{x(1+x)} dx = -2\varsigma^2(2) + \sum_{n=1}^{\infty} \frac{H_n^{(2)}}{n^2} + 2\sum_{n=1}^{\infty} \frac{H_n}{n^3}$$

and also

(4.4.100kv) $$\int_0^1 Li_p(x) \log x \frac{Li_2(x)}{x} dx = -2\varsigma(2)\varsigma(p+2) + \sum_{n=1}^{\infty} \frac{H_n^{(2)}}{n^{p+2}} + 2\sum_{n=1}^{\infty} \frac{H_n}{n^{p+3}}$$

A further differentiation results in

(4.4.100li) $$\int_0^1 x^\alpha \log^2 x \, Li_2(x) dx = 2\frac{\varsigma(2)}{(\alpha+1)^3} - \frac{\psi''(\alpha+2)}{(\alpha+1)^2} + 4\frac{\psi'(\alpha+2)}{(\alpha+1)^3} - 6\frac{\psi(\alpha+2)+\gamma}{(\alpha+1)^4}$$

With $\alpha = n-1$ we have

(4.4.100lii) $$\int_0^1 x^n \log^2 x \frac{Li_2(x)}{x} dx = 2\frac{\varsigma(2)}{n^3} - \frac{\psi''(n+1)}{n^2} + 4\frac{\psi'(n+1)}{n^3} - 6\frac{H_n}{n^4}$$

Using (E.16a) we see that

$$\psi''(n+1) = -2\sum_{k=0}^{\infty} \frac{1}{(n+1+k)^3} = 2H_n^{(3)} - 2\varsigma(3)$$

and accordingly we get

(4.4.100liii) $$\int_0^1 x^n \log^2 x \frac{Li_2(x)}{x} dx = 2\frac{\varsigma(2)}{n^3} - \frac{2H_n^{(3)} - 2\varsigma(3)}{n^2} + 4\frac{\varsigma(2) - H_n^{(2)}}{n^3} - 6\frac{H_n}{n^4}$$

As before, the summation process gives us

(4.4.100liv) $$\int_0^1 \frac{\log^2 x \, Li_2(x)}{x(1+x)} dx = 8\varsigma(2)\varsigma(3) - 4\sum_{n=1}^{\infty} \frac{H_n^{(2)}}{n^3} - 2\sum_{n=1}^{\infty} \frac{H_n^{(3)}}{n^2} - 6\sum_{n=1}^{\infty} \frac{H_n}{n^4}$$

and generally

(4.4.100lv)

$$\int_0^1 Li_p(x)\log^2 x \frac{Li_2(x)}{x}dx = 6\varsigma(2)\varsigma(p+3)+2\varsigma(3)\varsigma(p+2)-2\sum_{n=1}^{\infty}\frac{H_n^{(3)}}{n^{p+2}}-4\sum_{n=1}^{\infty}\frac{H_n^{(2)}}{n^{p+3}}-6\sum_{n=1}^{\infty}\frac{H_n}{n^{p+4}}$$

In particular we have

(4.4.100lvi)

$$\int_0^1 Li_3(x)\log^2 x \frac{Li_2(x)}{x}dx = 6\varsigma(2)\varsigma(6)+2\varsigma(3)\varsigma(5)-2\sum_{n=1}^{\infty}\frac{H_n^{(3)}}{n^5}-4\sum_{n=1}^{\infty}\frac{H_n^{(2)}}{n^6}-6\sum_{n=1}^{\infty}\frac{H_n}{n^7}$$

We have

$$\int_0^1 Li_3(x)\log^2 x \frac{Li_2(x)}{x}dx = \frac{1}{2}\log^2 x [Li_3(x)]^2 \Big|_0^1 - \int_0^1 \frac{\log x}{x}[Li_3(x)]^2 dx$$

$$= -\int_0^1 \frac{\log x}{x}[Li_3(x)]^2 dx$$

$$\int_0^1 \frac{\log x}{x}[Li_3(x)]^2 dx = \int_0^1 \log x \, Li_3(x) \frac{Li_3(x)}{x}dx$$

$$= \log x \, Li_3(x)Li_4(x)\Big|_0^1 - \int_0^1 Li_4(x)\left[\frac{Li_3(x)}{x}+\log x \frac{Li_2(x)}{x}\right]dx$$

$$= -\int_0^1 Li_4(x)\left[\frac{Li_3(x)}{x}+\log x \frac{Li_2(x)}{x}\right]dx$$

$$= -\frac{1}{2}\varsigma^2(4) - \int_0^1 Li_4(x)\log x \frac{Li_2(x)}{x}dx$$

From (4.4.100kv) we know that

$$\int_0^1 Li_4(x)\log x \frac{Li_2(x)}{x}dx = -2\varsigma(2)\varsigma(6)+\sum_{n=1}^{\infty}\frac{H_n^{(2)}}{n^6}+2\sum_{n=1}^{\infty}\frac{H_n}{n^7}$$

Hence we have an alternative expression for the integral (see also (4.4.100lv))

$$\int_0^1 Li_3(x)\log^2 x \frac{Li_2(x)}{x}dx = \frac{1}{2}\varsigma^2(4) - 2\varsigma(2)\varsigma(6) + \sum_{n=1}^{\infty}\frac{H_n^{(2)}}{n^6} + 2\sum_{n=1}^{\infty}\frac{H_n}{n^7}$$

Accordingly we have the identity

$$6\varsigma(2)\varsigma(6) + 2\varsigma(3)\varsigma(5) - 2\sum_{n=1}^{\infty}\frac{H_n^{(3)}}{n^5} - 4\sum_{n=1}^{\infty}\frac{H_n^{(2)}}{n^6} - 6\sum_{n=1}^{\infty}\frac{H_n}{n^7} = \frac{1}{2}\varsigma^2(4) - 2\varsigma(2)\varsigma(6) + \sum_{n=1}^{\infty}\frac{H_n^{(2)}}{n^6} + 2\sum_{n=1}^{\infty}\frac{H_n}{n^7}$$

which is simplified to

$$5\sum_{n=1}^{\infty}\frac{H_n^{(2)}}{n^6} + 2\sum_{n=1}^{\infty}\frac{H_n^{(3)}}{n^5} = 8\varsigma(2)\varsigma(6) + 2\varsigma(3)\varsigma(5) - \frac{1}{2}\varsigma^2(4) - 8\sum_{n=1}^{\infty}\frac{H_n}{n^7}$$

Using (4.4.100jiii) this becomes

(4.4.100lvii) $$\qquad 5\sum_{n=1}^{\infty}\frac{H_n^{(2)}}{n^6} + 2\sum_{n=1}^{\infty}\frac{H_n^{(3)}}{n^5} = -\frac{21}{4}\varsigma(8) + 10\varsigma(3)\varsigma(5)$$

The above relation is recorded in [69, p.23].

It is clear that further differentiations will produce similar results for $\int_0^1 \frac{\log^r x\, Li_2(x)}{x(1+x)}dx$

and $\int_0^1 Li_p(x)\log^r x \frac{Li_2(x)}{x}dx$ involving higher order Euler sums.

We showed earlier that

$$\int_0^1 x^n \log^2 x \frac{Li_2(x)}{x}dx = -3\frac{\varsigma(2)}{n^3} - \frac{2H_n^{(3)} - 2\varsigma(3)}{n^2} + 4\frac{\varsigma(2) - H_n^{(2)}}{n^3} - 6\frac{H_n}{n^4}$$

A different integration of (4.4.100liii) results in

$$\int_0^1 x^n \log^2 x \frac{Li_2(x)}{x}dx = x^n \log^2 x\, Li_3(x)\Big|_0^1 - \int_0^1 Li_3(x)\left[nx^n \frac{\log^2 x}{x} + 2x^n \frac{\log x}{x}\right]dx$$

$$= -\int_0^1 Li_3(x)\left[nx^n \frac{\log^2 x}{x} + 2x^n \frac{\log x}{x}\right]dx$$

Summation then results in

$$\int_0^1 Li_p(x)\log^2 x \frac{Li_2(x)}{x}dx = -\int_0^1 Li_{p-1}(x)\log^2 x \frac{Li_3(x)}{x}dx - 2\int_0^1 Li_p(x)\log x \frac{Li_3(x)}{x}dx$$

and in particular

$$\int_0^1 Li_6(x)\log^2 x \frac{Li_2(x)}{x}dx = -\int_0^1 Li_5(x)\log^2 x \frac{Li_3(x)}{x}dx - 2\int_0^1 Li_6(x)\log x \frac{Li_3(x)}{x}dx$$

In fact, this is equivalent to integration by parts as noted below

$$\int_0^1 Li_6(x)\log^2 x \frac{Li_2(x)}{x}dx = Li_6(x)\log^2 x\, Li_3(x)\Big|_0^1 - \int_0^1 Li_3(x)\left[\log^2 x \frac{Li_5(x)}{x} + 2\log x \frac{Li_6(x)}{x}\right]dx$$

$$= -\int_0^1 Li_5(x)\log^2 x \frac{Li_3(x)}{x}dx - 2\int_0^1 Li_6(x)\log x \frac{Li_3(x)}{x}dx$$

A different integration by parts of (4.4.100kiii) produces

$$\int_0^1 x^n \log x \frac{Li_2(x)}{x}dx = x^n \log x\, Li_3(x)\Big|_0^1 - \int_0^1 Li_3(x)\left(nx^{n-1}\log x + x^{n-1}\right)dx$$

$$= -\int_0^1 Li_3(x)\left(nx^{n-1}\log x + x^{n-1}\right)dx$$

We have seen previously that

$$\int_0^1 x^n \frac{Li_3(x)}{x}dx = \frac{\varsigma(3)}{n} - \frac{\varsigma(2)}{n^2} + \frac{H_n}{n^3}$$

and hence we obtain

$$\int_0^1 x^n \log x \frac{Li_3(x)}{x}dx = -\frac{\varsigma(3)}{n^2} + \frac{\varsigma(2)}{n^3} - 3\frac{H_n}{n^4} - \frac{H_n^{(2)}}{n^3}$$

Completing the summation we get

(4.4.100m)

$$\int_0^1 Li_p(x)\log x \frac{Li_3(x)}{x}dx = -\varsigma(3)\varsigma(2+p)+\varsigma(2)\varsigma(3+p)-3\sum_{n=1}^{\infty}\frac{H_n}{n^{4+p}}-\sum_{n=1}^{\infty}\frac{H_n^{(2)}}{n^{3+p}}$$

Letting $p=4$ we have using integration by parts

$$\int_0^1 Li_4(x)\log x \frac{Li_3(x)}{x}dx = \frac{1}{2}\log x\left[Li_4(x)\right]^2\Big|_0^1 - \int_0^1 \frac{\left[Li_4(x)\right]^2}{x}dx = -\int_0^1 \frac{\left[Li_4(x)\right]^2}{x}dx$$

Freitas [69a] has evaluated the integral $\int_0^1 \frac{\left[Li_p(x)\right]^2}{x}dx$ and, in particular, we have

$$\int_0^1 \frac{\left[Li_4(x)\right]^2}{x}dx = 2\sum_{j=1}^{2}\varsigma(2j)\varsigma(9-2j) = 2\varsigma(2)\varsigma(7)+2\varsigma(4)\varsigma(5)-5\varsigma(9)$$

Therefore we conclude that

$$-\varsigma(3)\varsigma(6)+\varsigma(2)\varsigma(7)-3\sum_{n=1}^{\infty}\frac{H_n}{n^8}-\sum_{n=1}^{\infty}\frac{H_n^{(2)}}{n^7} = -2\varsigma(2)\varsigma(7)-2\varsigma(4)\varsigma(5)+5\varsigma(9)$$

From Flajolet and Salvy [69] we see that

$$\sum_{n=1}^{\infty}\frac{H_n}{n^q} = \left(1+\frac{q}{2}\right)\varsigma(q+1)-\frac{1}{2}\sum_{k=1}^{q-2}\varsigma(k+1)\varsigma(q-k)$$

and we therefore obtain

(4.4.100n) $$\sum_{n=1}^{\infty}\frac{H_n^{(2)}}{n^7} = -20\varsigma(9)+\frac{3}{2}\sum_{k=1}^{6}\varsigma(k+1)\varsigma(8-k)-\varsigma(3)\varsigma(6)+3\varsigma(2)\varsigma(7)$$

It may be noted from the recent paper by Borwein and Bradley, "Thirty-two Goldbach Variations" [30a] that there is a connection with the Witten $\varsigma$-function $W(r,s,t)$ defined for $r,s > 1/2$ by

$$W(r,s,t) = \sum_{n=1}^{\infty}\sum_{m=1}^{\infty}\frac{1}{n^r m^s (n+m)^t}$$

which also has the integral representation

(4.4.100o) $$W(r,s,t) = \frac{(-1)^{t-1}}{\Gamma(t)}\int_0^1 Li_r(x)Li_s(x)\frac{\log^{t-1}x}{x}dx$$

The $W(r,s,t)$ are also referred to as Tornheim double sums.

Following Lewin's procedure [100, p.20], integration by parts shows that

$$\int_0^t x^\alpha Li_p(x)\,dx = \frac{t^{\alpha+1}}{\alpha+1} Li_p(t) - \frac{1}{\alpha+1}\int_0^t x^{\alpha+1} \frac{Li_{p-1}(x)}{x}\,dx$$

$$= \frac{t^{\alpha+1}}{\alpha+1} Li_p(t) - \frac{1}{\alpha+1}\int_0^t x^\alpha Li_{p-1}(x)\,dx$$

$$= \frac{t^{\alpha+1}}{\alpha+1} Li_p(t) - \frac{t^{\alpha+1}}{(\alpha+1)^2} Li_{p-1}(t) + \frac{1}{(\alpha+1)^2}\int_0^t x^\alpha Li_{p-2}(x)\,dx$$

and hence we have

$$\int_0^t x^\alpha Li_p(x)\,dx = t^{\alpha+1}\sum_{k=0}^{p-2}(-1)^k \frac{Li_{p-k}(t)}{(\alpha+1)^{k+1}} + \frac{(-1)^p}{(\alpha+1)^{p-1}}\int_0^t x^{\alpha+1}\frac{\log(1-x)}{x}\,dx$$

$$= t^{\alpha+1}\sum_{k=0}^{p-2}(-1)^k \frac{Li_{p-k}(t)}{(\alpha+1)^{k+1}} + (-1)^p \frac{t^{\alpha+1}-1}{(\alpha+1)^p}\log(1-t) - \frac{(-1)^p}{(\alpha+1)^p}\int_0^t \frac{1-x^{\alpha+1}}{1-x}\,dx$$

Therefore, letting $t=1$ we get

(4.4.100p) $$\int_0^1 x^\alpha Li_p(x)\,dx = \sum_{k=0}^{p-2}(-1)^k \frac{\varsigma(p-k)}{(\alpha+1)^{k+1}} - \frac{(-1)^p}{(\alpha+1)^p}[\psi(\alpha+2)+\gamma]$$

With $t=1$ and $\alpha = n-1$ we have

$$\int_0^1 x^n \frac{Li_p(x)}{x}\,dx = \sum_{k=0}^{p-2}(-1)^k \frac{\varsigma(p-k)}{n^{k+1}} - (-1)^p \frac{H_n^{(1)}}{n^p}$$

Completing the summation of the above we see that

(4.4.100pi) $$\int_0^1 \frac{Li_p(x)}{x(1+x)}\,dx = -\sum_{k=0}^{p-2}(-1)^k \varsigma_a(k+1)\varsigma(p-k) - (-1)^p \sum_{n=1}^\infty (-1)^n \frac{H_n^{(1)}}{n^p}$$

(4.4.100pii) $$\int_0^1 \frac{Li_p(x)Li_q(x)}{x}\,dx = \sum_{k=0}^{p-2}(-1)^k \varsigma(p-k)\varsigma(q+k+1) - (-1)^p \sum_{n=1}^{\infty} \frac{H_n^{(1)}}{n^{p+q}}$$

With $q = p-1$ we obtain

$$\int_0^1 \frac{Li_p(x)Li_{p-1}(x)}{x}\,dx = \frac{1}{2}\varsigma^2(p)$$

and hence we get

(4.4.100piii) $$\frac{1}{2}\varsigma^2(p) = \sum_{k=0}^{p-2}(-1)^k \varsigma(p-k)\varsigma(p+k) - (-1)^p \sum_{n=1}^{\infty} \frac{H_n^{(1)}}{n^{2p-1}}$$

Differentiating (4.4.100p) gives us

(4.4.100piv) $$\int_0^1 x^\alpha \log^q x \, Li_p(x)\,dx =$$

$$\sum_{k=0}^{p-2} \frac{(-1)^{k+q}}{(\alpha+1)^{q+k+1}}(k+1)(k+2)\ldots(k+q+1)\varsigma(p-k) - \frac{(-1)^p}{(\alpha+1)^p}\psi^{(q)}(\alpha+2)$$

With $\alpha = n-1$ we obtain

(4.4.100pv) $$\int_0^1 x^n \log^q x \frac{Li_p(x)}{x}\,dx =$$

$$\sum_{k=0}^{p-2} \frac{(-1)^{k+q}}{n^{q+k+1}}(k+1)(k+2)\ldots(k+q+1)\varsigma(p-k) - \frac{(-1)^{p+q}q!}{n^p}\left[H_n^{(q+1)} - \varsigma(q+1)\right]$$

Upon making the familiar summation we may easily obtain an expression for
$$\int_0^1 \log^q x \frac{Li_p(x)Li_r(x)}{x}\,dx.$$

In (4.4.100p) let $\alpha = n + \beta - 1$

$$\int_0^1 x^\beta x^n \frac{Li_p(x)}{x}\,dx = \sum_{k=0}^{p-2}(-1)^k \frac{\varsigma(p-k)}{(\beta+n)^{k+1}} - \frac{(-1)^p}{(\beta+n)^p}\left[\psi(\beta+n+1) + \gamma\right]$$

and making another summation we get

(4.4.100pvi)

$$\int_0^1 x^\beta \frac{Li_p(x)Li_q(x)}{x} dx = \sum_{k=0}^{p-2}(-1)^k \varsigma(p-k)\sum_{n=1}^\infty \frac{1}{n^q(\beta+n)^{k+1}} - (-1)^p \sum_{n=1}^\infty \frac{[\psi(\beta+n+1)+\gamma]}{n^q(\beta+n)^p}$$

The substitution $y = e^{-x}$ in (4.4.93) gives us a form reminiscent of the Beta function

(4.4.100q) $$\int_0^1 \frac{y^{a-1}(1-y)^n}{\log^r y} dy = \frac{1}{(r-1)!}\sum_{k=0}^n \binom{n}{k}(-1)^k (a+k)^{r-1} \log(a+k)$$

and we have from (4.4.13)

$$g(x) = \frac{n!}{x(1+x)\ldots(n+x)} = \int_0^1 t^{x-1}(1-t)^n dt = \sum_{k=0}^n \binom{n}{k}\frac{(-1)^k}{k+x}$$

Integrating this with respect to $x$ we obtain

$$\int_1^a dx \int_0^1 t^{x-1}(1-t)^n dt = \sum_{k=0}^n \binom{n}{k}(-1)^k \log\frac{k+a}{k+1}$$

and we also have

$$\int_1^a dx \int_0^1 t^{x-1}(1-t)^n dt = \int_0^1 \frac{(t^{a-1}-1)(1-t)^n}{\log t} dt = \int_0^1 \frac{t^{a-1}(1-t)^n}{\log t} dt - \int_0^1 \frac{(1-t)^n}{\log t} dt$$

Employing (4.4.94a) this becomes

$$\int_0^1 \frac{t^{a-1}(1-t)^n}{\log t} dt = \sum_{k=0}^n \binom{n}{k}(-1)^k \log\frac{k+a}{k+1} + \sum_{k=0}^n \binom{n}{k}(-1)^{k+1}\log(1+k)$$

and thus we have

(4.4.100qi) $$\int_0^1 \frac{t^{a-1}(1-t)^n}{\log t} dt = \sum_{k=0}^n \binom{n}{k}(-1)^{k+1} \log(a+k)$$

Integration with respect to $a$ gives us

(4.4.100qii) $$\int_0^1 \frac{(t^{\beta-1}-t^{\beta-1})(1-t)^n}{\log^2 t} dt =$$

$$\sum_{k=0}^{n}\binom{n}{k}(-1)^{k+1}(\beta+k)\log(\beta+k)-\sum_{k=0}^{n}\binom{n}{k}(-1)^{k+1}(\alpha+k)\log(\alpha+k)-\beta+\alpha$$

It is clear that a further integration with respect to $a$ will result in a corresponding identity for $\int_{0}^{1}\frac{(t^{\beta-1}-t^{\beta-1})(1-t)^{n}}{\log^{3}t}dt$ and so on.

Completing the summation of (4.4.100qi) we get

$$\sum_{n=1}^{\infty}\frac{x^{n}}{n^{p}}\int_{0}^{1}\frac{t^{a-1}(1-t)^{n}}{\log t}dt=\sum_{n=1}^{\infty}\frac{x^{n}}{n^{p}}\sum_{k=0}^{n}\binom{n}{k}(-1)^{k+1}\log(a+k)$$

and hence we have

(4.4.100qiii)  $$\int_{0}^{1}\frac{t^{a-1}Li_{p}[x(1-t)]}{\log t}dt=\sum_{n=1}^{\infty}\frac{x^{n}}{n^{p}}\sum_{k=0}^{n}\binom{n}{k}(-1)^{k+1}\log(a+k)$$

Differentiation with respect to $a$ results in

(4.4.100qiv)  $$\int_{0}^{1}t^{a-1}Li_{p}[x(1-t)]dt=\sum_{n=1}^{\infty}\frac{x^{n}}{n^{p}}\sum_{k=0}^{n}\binom{n}{k}\frac{(-1)^{k+1}}{a+k}$$

(4.4.100qv)  $$\int_{0}^{1}t^{a-1}\log t\, Li_{p}[x(1-t)]dt=-\sum_{n=1}^{\infty}\frac{x^{n}}{n^{p}}\sum_{k=0}^{n}\binom{n}{k}\frac{(-1)^{k+1}}{(a+k)^{2}}$$

(4.4.100qvi)  $$\int_{0}^{1}t^{a-1}\log^{2}t\, Li_{p}[x(1-t)]dt=2\sum_{n=1}^{\infty}\frac{x^{n}}{n^{p}}\sum_{k=0}^{n}\binom{n}{k}\frac{(-1)^{k+1}}{(a+k)^{3}}$$

See also (4.4.24aaa).

The Wolfram Integrator is only able to produce a result for the dilogarithm which it does in terms of the hypergeometric functions

(4.4.100qvii)  $$\int_{0}^{t}x^{a-1}Li_{2}(1-x)dx=$$

$$\frac{1}{a(a+1)^{2}}t^{a}\left[t\,_{3}F_{2}(1,a+1,a+1;a+2;t)+(a+1)\left\{(a+1)Li_{2}(1-t)t\,_{2}F_{1}(1,a+1;a+2;t)\log t\right\}\right]$$

Letting $a = 1$ and using (4.4.123)

$$\sum_{k=0}^{n}\binom{n}{k}\frac{(-1)^k}{k+1} = \frac{1}{n+1}$$

results in

(4.4.100qviii) $\quad \int_0^1 Li_p[x(1-t)]dt = \sum_{n=1}^{\infty}\frac{1}{(n+1)n^p}x^n = \int_0^x Li_p(t)dt$

Similarly using (4.4.127)

$$\sum_{k=0}^{n}\binom{n}{k}\frac{(-1)^k}{(k+1)^2} = \frac{H_{n+1}}{n+1}$$

we obtain

(4.4.100qix) $\quad \int_0^1 \log t\, Li_p[x(1-t)]dt = -\sum_{n=1}^{\infty}\frac{H_{n+1}^{(1)}}{(n+1)n^p}x^n$

With (4.4.130)

$$\sum_{k=0}^{n}\binom{n}{k}\frac{(-1)^{k+1}}{(1+k)^3} = \frac{1}{2(n+1)}\left[\left(H_{n+1}^{(1)}\right)^2 + H_{n+1}^{(2)}\right]$$

we get

(4.4.100qx) $\quad \int_0^1 \log^2 t\, Li_p[x(1-t)]dt = \sum_{n=1}^{\infty}\frac{1}{(n+1)n^p}\left[\left(H_{n+1}^{(1)}\right)^2 + H_{n+1}^{(2)}\right]x^n$

We now refer back to (4.4.93) with $a = r = 1$

$$f(n) = \int_0^{\infty}\frac{e^{-x}(1-e^{-x})^n}{x}dx = \sum_{k=0}^{n}\binom{n}{k}(-1)^{k+1}\log(1+k)$$

Applying the integral identity

$$\frac{1}{n+1} = \int_0^{\infty} e^{-nt}e^{-t}dt$$

we have

$$\sum_{n=1}^{\infty}\frac{f(n)}{(n+1)^2}=\sum_{n=1}^{\infty}\int_0^{\infty}e^{-nt}e^{-t}dt\int_0^{\infty}e^{-nu}e^{-u}du\int_0^{\infty}\frac{e^{-x}(1-e^{-x})^n}{x}dx$$

The geometric series gives us

$$\sum_{n=1}^{\infty}\left[e^{-t}e^{-u}(1-e^{-x})\right]^n=\frac{e^{-(t+u)}(1-e^{-x})}{1-e^{-(t+u)}(1-e^{-x})}$$

and we may therefore write the integral as

$$I=\int_0^{\infty}\frac{e^{-x}}{x}dx\int_0^{\infty}e^{-u}du\int_0^{\infty}\frac{e^{-2t}\left[e^{-u}(1-e^{-x})\right]}{1-e^{-t}\left[e^{-u}(1-e^{-x})\right]}dt$$

Therefore, as before, we have

$$\int_0^{\infty}\frac{e^{-2t}\left[e^{-u}(1-e^{-x})\right]}{1-e^{-t}\left[e^{-u}(1-e^{-x})\right]}dt=-1-\frac{\log\left[1-e^{-u}(1-e^{-x})\right]}{\left[e^{-u}(1-e^{-x})\right]}$$

Substituting this in the above we obtain

$$I=-\int_0^{\infty}\frac{e^{-x}}{x}dx\int_0^{\infty}e^{-u}\left\{1+\frac{\log\left[1-e^{-u}(1-e^{-x})\right]}{\left[e^{-u}(1-e^{-x})\right]}\right\}du$$

$$=-\int_0^{\infty}\frac{e^{-x}}{x}dx\int_0^{\infty}\left\{e^{-u}+\frac{\log\left[1-e^{-u}(1-e^{-x})\right]}{(1-e^{-x})}\right\}du$$

We have seen from (4.4.57g) that

$$Li_{s+1}(t)=-\frac{1}{\Gamma(s)}\int_0^t u^{s-1}\log(1-te^{-u})du$$

and hence we have with $s=1$

$$I=-\int_0^{\infty}\frac{e^{-x}}{x}\left\{1-\frac{Li_2\left[1-e^{-x}\right]}{\left[1-e^{-x}\right]}\right\}dx$$

We therefore obtain

$$-\int_0^\infty \frac{e^{-x}}{x}\left\{1-\frac{Li_2\left[1-e^{-x}\right]}{\left[1-e^{-x}\right]}\right\}dx = \sum_{n=1}^\infty \frac{1}{(n+1)^2}\sum_{k=0}^n \binom{n}{k}(-1)^{k+1}\log(k+1)$$

and the substitution $y = e^{-x}$ results in

(4.4.100r) $$\int_0^1 \left\{1-\frac{Li_2(1-y)}{1-y}\right\}\frac{dy}{\log y} = \sum_{n=1}^\infty \frac{1}{(n+1)^2}\sum_{k=0}^n \binom{n}{k}(-1)^{k+1}\log(k+1)$$

A generalised version of this is presented in (4.4.100zi).

Completing the summation of (4.4.94a)

$$\int_0^1 \frac{(1-y)^n}{\log y} dy = \sum_{k=0}^n \binom{n}{k}(-1)^k \log(1+k)$$

we have

$$\sum_{n=1}^\infty \frac{1}{p^n}\int_0^1 \frac{(1-y)^n}{\log y} dy = \sum_{n=1}^\infty \frac{1}{p^n}\sum_{k=0}^n \binom{n}{k}(-1)^k \log(1+k)$$

The geometric series gives us for $|p| < 1$

$$\sum_{n=1}^\infty \frac{1}{p^n}\int_0^1 \frac{(1-y)^n}{\log y} dy = \int_0^1 \frac{(1-y)}{[p-(1-y)]\log y} dy$$

and hence we obtain

(4.4.100s) $$\int_0^1 \frac{1-y}{[1-t(1-y)]\log y} dy = \sum_{n=1}^\infty \frac{1}{p^n}\sum_{k=0}^n \binom{n}{k}(-1)^k \log(1+k)$$

If we let $p = 1/t$ in (4.4.100s) we rediscover the formula found by Guillera (see Sondow's paper [123a])

(4.4.100t) $$t^2\int_0^1 \frac{1-y}{[1-t(1-y)]\log y} dy = \sum_{n=1}^\infty t^{n+1}\sum_{k=0}^n \binom{n}{k}(-1)^k \log(1+k)$$

With $t = 1/2$ we get

(4.4.100u) $$\int_0^1 \frac{y-1}{(1+y)\log y}\,dy = \sum_{n=1}^\infty \frac{1}{2^n}\sum_{k=0}^n \binom{n}{k}(-1)^{k+1}\log(1+k) = \log\frac{\pi}{2}$$

where in the final part we have used (4.4.112).

We now revisit the Hasse identity (3.12) which is valid for all $s$ except $s=1$

$$\varsigma(s) = \frac{1}{s-1}\sum_{n=0}^\infty \frac{1}{n+1}\sum_{k=0}^n \binom{n}{k}\frac{(-1)^k}{(k+1)^{s-1}}$$

and with differentiation we get

$$\varsigma'(s) = -\frac{1}{(s-1)^2}\sum_{n=0}^\infty \frac{1}{n+1}\sum_{k=0}^n \binom{n}{k}\frac{(-1)^k}{(k+1)^{s-1}} - \frac{1}{s-1}\sum_{n=0}^\infty \frac{1}{n+1}\sum_{k=0}^n \binom{n}{k}\frac{(-1)^k \log(k+1)}{(k+1)^{s-1}}$$

$$= -\frac{1}{s-1}\left[\varsigma(s) + \sum_{n=0}^\infty \frac{1}{n+1}\sum_{k=0}^n \binom{n}{k}\frac{(-1)^k \log(k+1)}{(k+1)^{s-1}}\right]$$

Therefore we have

(4.4.100v) $$(s-1)\varsigma'(s) + \varsigma(s) = \sum_{n=0}^\infty \frac{1}{n+1}\sum_{k=0}^n \binom{n}{k}\frac{(-1)^{k+1}\log(k+1)}{(k+1)^{s-1}}$$

and reference to (4.4.92a) suggests that

(4.4.100w) $$\lim_{s\to 1}[(s-1)\varsigma'(s) + \varsigma(s)] = -\sum_{n=0}^\infty \frac{1}{n+1}\sum_{k=0}^n \binom{n}{k}(-1)^k \log(1+k) = \gamma$$

Using Riemann's functional equation

$$\varsigma(s) = 2(2\pi)^{s-1}\Gamma(1-s)\sin(\pi s/2)\varsigma(1-s)$$

it is easily seen that

$$\varsigma(0) = \lim_{s\to 0}\left[\frac{1}{\pi}\varsigma(1-s)\sin(\pi s/2)\right]$$

$$= \lim_{s\to 0}\left[\frac{s}{2}\varsigma(1-s)\frac{\sin(\pi s/2)}{\pi s/2}\right]$$

$$= \lim_{s \to 0}\left[\frac{s}{2}\varsigma(1-s)\right]$$

and using the fact that $\varsigma(0) = -\frac{1}{2}$ we have

$$-1 = \lim_{s \to 0}[s\varsigma(1-s)]$$

It should be noted that the above limit is negative because, as may be seen from (1.1) in Volume I, $\varsigma(x)$ is negative for $0 < x < 1$.

Alternatively, changing the limit we get

(4.4.100x) $\qquad 1 = \lim_{s \to 1}[(s-1)\varsigma(s)]$

We may therefore apply L'Hôpital's rule to $\lim_{s \to 1}\left[\frac{\varsigma(s)(s-1)-1}{s-1}\right]$ since both the numerator and the denominator approach zero as $s \to 1$. It is known from [135, p.271] and (4.4.99m) that $\lim_{s \to 1}\left[\varsigma(s) - \frac{1}{s-1}\right] = \gamma$ and hence we obtain

$$\lim_{s \to 1}[(s-1)\varsigma'(s) + \varsigma(s)] = \gamma$$

This accordingly gives us yet another proof of

$$\lim_{s \to 1}[(s-1)\varsigma'(s) + \varsigma(s)] = -\sum_{n=0}^{\infty}\frac{1}{n+1}\sum_{k=0}^{n}\binom{n}{k}(-1)^k \log(1+k) = \gamma$$

In this regard, see the limit given by (F.5b) in Appendix F.

$$\frac{\varsigma'(s)}{\varsigma(s)} = \frac{-1\{1/(s-1)^2\} + k + \ldots}{\{1/(s-1)\} + \gamma + k(s-1)\ldots} = -\frac{1}{s-1} + \gamma + \ldots$$

Upon differentiating (4.4.100v) we obtain

$$(s-1)\varsigma''(s) + 2\varsigma'(s) = \sum_{n=0}^{\infty}\frac{1}{n+1}\sum_{k=0}^{n}\binom{n}{k}\frac{(-1)^k \log^2(k+1)}{(k+1)^{s-1}}$$

and as $s \to 1$ we get

$$\lim_{s \to 1}[(s-1)\varsigma''(s) + 2\varsigma'(s)] = \sum_{n=0}^{\infty} \frac{1}{n+1} \sum_{k=0}^{n} \binom{n}{k} (-1)^k \log^2(k+1)$$

We may recall from (4.2.5) that

$$\lim_{s \to 1}\left[(1-2^{1-s})\varsigma(s)\right] = \log 2$$

Equivalently we have

$$\lim_{s \to 1}\left[2^{s-1}(1-2^{1-s})\varsigma(s)\right] = \lim_{s \to 1}\left[(2^{s-1}-1)\varsigma(s)\right] = \log 2 \lim_{s \to 1} 2^{s-1} = \log 2$$

and we may write this as

$$\lim_{s \to 1}\left[\left(\frac{2^{s-1}-1}{s-1}\right)(s-1)\varsigma(s)\right] = \log 2$$

Since by the definition of a derivative

$$\lim_{s \to 1}\left(\frac{2^{s-1}-1}{s-1}\right) = \log 2$$

we see again that $\lim_{s \to 1}[(s-1)\varsigma(s)] = 1$.

Integrating Guillera's formula (4.4.100t) we get (having first divided both sides by $t$)

$$\int_0^x t\,dt \int_0^1 \frac{1-y}{1-t(1-y)\log y}\,dy = \sum_{n=1}^{\infty} \frac{x^{n+1}}{n+1} \sum_{k=0}^{n} \binom{n}{k} (-1)^k \log(1+k)$$

Reversing the order of integration we have

$$\int_0^x t\,dt \int_0^1 \frac{1-y}{[1-t(1-y)]\log y}\,dy = \int_0^1 \frac{1-y}{\log y}\,dy \int_0^x \frac{t}{[1-t(1-y)]}\,dt$$

We see that

$$\int_0^x \frac{t}{[1-t(1-y)]}\,dt = \frac{1}{1-y}\left(\int_0^x \frac{1}{[1-t(1-y)]} - 1\right)dt$$

$$= -\frac{1}{(1-y)^2}\log[1-(1-y)x] - \frac{x}{1-y}$$

and accordingly we have

$$\int_0^1 \frac{1-y}{\log y}\,dy \int_0^x \frac{t}{[1-t(1-y)]}\,dt = -\int_0^1 \left(\frac{\log[1-(1-y)x]}{(1-y)\log y} + \frac{x}{\log y}\right)dy$$

With $x = 1$ we get

$$\int_0^1 t\,dt \int_0^1 \frac{1-y}{1-t(1-y)\log y}\,dy = \int_0^1 \left(\frac{1}{1-y} + \frac{1}{\log y}\right)dy$$

which we have previously shown in (4.4.99) is equal to $\gamma$.

Hence we rediscover (4.4.92a)

$$\sum_{n=1}^{\infty} \frac{1}{n+1} \sum_{k=0}^{n} \binom{n}{k} (-1)^k \log(1+k) = -\gamma$$

We also have

(4.4.100y) $$\int_0^1 \left(\frac{\log[1-(1-y)x]}{(1-y)\log y} + \frac{x}{\log y}\right)dy = \sum_{n=1}^{\infty} \frac{x^{n+1}}{n+1} \sum_{k=0}^{n} \binom{n}{k} (-1)^{k+1} \log(1+k)$$

Differentiating the above equation with respect to $x$ we obtain

(4.4.100z) $$\int_0^1 \left(\frac{-1}{[1-(1-y)x]\log y} + \frac{1}{\log y}\right)dy = \sum_{n=1}^{\infty} x^n \sum_{k=0}^{n} \binom{n}{k} (-1)^{k+1} \log(1+k)$$

With an eye on (4.4.100t) Sondow defines $f(t,s)$ for $t \in (-1,1)$ and $s \in \mathbf{C}$ by

$$f(t,s) = \sum_{n=0}^{\infty} t^{n+1} \sum_{k=0}^{n} \binom{n}{k} \frac{(-1)^k}{(k+1)^s}$$

and notes that $f(1/2, s) = \varsigma_a(s)$ by reference to the Hasse/Sondow identity (3.12).

Dividing (4.4.100y) by $x$ and then integrating results in

$$\int_o^t dx \int_0^1 \left( \frac{\log[1-(1-y)x]}{x(1-y)\log y} + \frac{1}{\log y} \right) dy = \sum_{n=1}^{\infty} \frac{t^{n+1}}{(n+1)^2} \sum_{k=0}^{n} \binom{n}{k} (-1)^{k+1} \log(1+k)$$

Then, reversing the order of integration we get

$$\int_o^t dx \int_0^1 \left( \frac{\log[1-(1-y)x]}{x(1-y)\log y} + \frac{1}{\log y} \right) dy = \int_o^1 \frac{dy}{\log y} \int_0^t \left( \frac{\log[1-(1-y)x]}{x(1-y)} + 1 \right) dx$$

We have

$$\int_0^t \left( \frac{\log[1-(1-y)x]}{x(1-y)} + 1 \right) dx = -\frac{Li_2[(1-y)x]}{1-y} + x \bigg|_0^t = t - \frac{Li_2[(1-y)t]}{1-y}$$

and hence we obtain

(4.4.100zi) $$\int_o^1 \left\{ t - \frac{Li_2[(1-y)t]}{1-y} \right\} \frac{dy}{\log y} = \sum_{n=1}^{\infty} \frac{t^{n+1}}{(n+1)^2} \sum_{k=0}^{n} \binom{n}{k} (-1)^{k+1} \log(1+k)$$

This is a generalised version of (4.4.100r) where we had $t=1$. See also (4.4.116c) for a different derivation.

Let us now divide (4.4.100zi) by $t$ and integrate. Since

$$\int_0^x \frac{Li_2[\alpha t]}{t} dt = \int_0^{\alpha x} \frac{Li_2[y]}{y} dy = Li_3[\alpha x]$$

we have

(4.4.100zii) $$\int_o^1 \left\{ t - \frac{Li_3[(1-y)t]}{1-y} \right\} \frac{dy}{\log y} = \sum_{n=1}^{\infty} \frac{t^{n+1}}{(n+1)^3} \sum_{k=0}^{n} \binom{n}{k} (-1)^{k+1} \log(1+k)$$

and it is clear that this generalises to

(4.4.100ziii) $$\int_o^1 \left\{ t - \frac{Li_p[(1-y)t]}{1-y} \right\} \frac{dy}{\log y} = \sum_{n=1}^{\infty} \frac{t^{n+1}}{(n+1)^p} \sum_{k=0}^{n} \binom{n}{k} (-1)^{k+1} \log(1+k)$$

# A LITTLE BIT OF $\log \frac{\pi}{2}$

In his 1997 paper, "A class of logarithmic integrals", [2a] Adamchik produced the rather astonishing integral

$$(4.4.101) \quad \int_0^1 \frac{x^{p-1}}{(1+x^n)^2} \log\log\left(\frac{1}{x}\right) dx = \frac{(n-p)(\log(2n)+\gamma)\left[\psi\left(\frac{p}{2n}\right) - \psi\left(\frac{n+p}{2n}\right)\right]}{2n^2}$$

$$-\frac{1}{2n}\left\{\gamma + \log(2n) - 2\log\frac{\Gamma\left(\frac{p}{2n}\right)}{\Gamma\left(\frac{n+p}{2n}\right)}\right\} + \frac{(n-p)}{2n^2}\left\{\varsigma'\left(1,\frac{p}{2n}\right) - \varsigma'\left(1,\frac{n+p}{2n}\right)\right\}$$

where $\varsigma'(s,z) = \frac{\partial}{\partial s}\varsigma(s,z)$ denotes the derivative of the Hurwitz zeta function and $\text{Re}(p) > 0$, $\text{Re}(n) > 0$. Adamchik's proof of a similar integral is shown in Appendix C in Volume VI of this series of papers.

Fortuitously, for $p = n$ this greatly simplifies to

$$(4.4.101a) \quad \int_0^1 \frac{x^{n-1}}{(1+x^n)^2} \log\log\left(\frac{1}{x}\right) dx = -\frac{1}{2n}\left[\gamma + \log\left(\frac{2n}{\pi}\right)\right]$$

In (4.4.101a) let $n = 1$ and hence we have

$$(4.4.101b) \quad \int_0^1 \frac{\log\log(1/t)}{(t+1)^2} dt = \frac{1}{2}\left(\log\frac{\pi}{2} - \gamma\right)$$

This integral is contained in G&R [74, p.565].

In Theorem 4.7 below we show that

$$(4.4.101c) \quad \log\frac{\pi}{2} = \int_0^1 \frac{(t-1)}{(t+1)\log t} dt$$

We have

$$\int_0^1 \frac{(t-1)}{(t+1)\log t} dt = \int_0^1 \frac{t(t-1)}{(t+1)} \left\{ \frac{-1/t}{\log(1/t)} \right\} dt$$

$$= \frac{t(t-1)}{(t+1)} \log\log(1/t) \Big|_0^1 - \int_0^1 \frac{t^2+2t-1}{(t+1)^2} \log\log(1/t) dt$$

$$= -\int_0^1 \frac{\left[(t+1)^2 - 2\right]}{(t+1)^2} \log\log(1/t) dt$$

$$= 2\int_0^1 \frac{\log\log(1/t)}{(t+1)^2} dt - \int_0^1 \log\log(1/t) dt$$

The gamma function is defined in (4.3.1) as

$$\Gamma(x) = \int_0^\infty y^{x-1} e^{-y} dy$$

and, with the substitution $y = -\log t$, this becomes

$$\Gamma(x) = \int_0^1 \left(\log(1/t)\right)^{x-1} dt$$

We therefore have

$$\Gamma'(x) = \int_0^1 \left(\log(1/t)\right)^{x-1} \log\log(1/t) dt$$

$$\Gamma'(1) = \int_0^1 \log\log(1/t) dt = -\gamma$$

and hence we have

(4.4.101d) $$\int_0^1 \frac{\log\log(1/t)}{(t+1)^2} dt = \frac{1}{2}\left(\log\frac{\pi}{2} - \gamma\right)$$

This is a particular case (with $p = n = 1$) of the general formula (4.4.100) derived by Adamchik in [2a] in 1997.

We have

(4.4.101e) $$\log\frac{\pi}{2} = \int_0^1 \frac{t(t-1)}{t+1}\frac{1}{t\log t}dt$$

and, using integration by parts, and the fact that $\int \frac{1}{t\log t}dt = \int \frac{1/t}{\log t}dt = \log\log t$ we obtain

(4.4.101f) $$= \frac{t(t-1)}{(t+1)}\log\log t \Big|_0^1 - \int_0^1 \frac{t}{t+1}\log\log t\, dt$$

(4.4.101g) $$= -\int_0^1 \frac{t}{t+1}\log\log t\, dt$$

Combining (4.4.101b) and (4.4.101g) we obtain

(4.4.101h) $$\gamma = 2\int_0^1 \frac{\log\log(1/t)}{(t+1)^2}dt - \int_0^1 \frac{t\log\log t}{t+1}dt$$

**(x) Theorem 4.7:**

(4.4.102) $$\log\frac{\pi}{2} = \sum_{n=1}^\infty (-1)^{n+1}\log\frac{n+1}{n} = \int_0^1 dx \int_0^1 \frac{y^x}{1+y}dy = \int_0^1 \frac{(y-1)}{(1+y)\log y}dy$$

**Proof:**

We have Wallis's infinite product for $\pi$ (which was published in Arithmetica Infinitorum in 1659).

(4.4.103) $$\frac{\pi}{2} = \frac{2}{1}\frac{2}{3}\frac{4}{3}\frac{4}{5}\frac{6}{5}\frac{6}{7}\frac{8}{7}\cdots$$

and taking the logarithm of (4.4.103) we obtain

(4.4.104) $$\log\frac{\pi}{2} = \sum_{n=1}^\infty (-1)^{n+1}\log\frac{n+1}{n}$$

We have (formally)

(4.4.105) $$S = \sum_{n=1}^\infty (-1)^{n+1}\log\frac{n+1}{n} = \sum_{n=1}^\infty (-1)^{n+1}\int_0^1 \frac{dx}{n+x}$$

Using the elementary integral

(4.4.106) $$\frac{1}{n+x} = \int_0^1 y^{n+x-1} dy$$

we deduce that

(4.4.107) $$S = \sum_{n=1}^{\infty} (-1)^{n+1} \int_0^1 dx \int_0^1 y^{n+x-1} dy$$

Applying the geometric series we obtain

(4.4.108) $$S = \int_0^1 dx \int_0^1 \frac{y^x}{1+y} dy$$

Interchanging the order of integration gives

(4.4.109) $$= \int_0^1 \frac{dy}{1+y} \int_0^1 y^x dx$$

(4.4.110) $$= \int_0^1 \frac{dy}{1+y} \left\{ \frac{y^x}{\log y} \Big|_0^1 \right\}$$

(4.4.111) $$= \int_0^1 \frac{(y-1)}{(1+y)\log y} dy$$

Accordingly, subject to a rigorous justification of the various intermediate steps, we have

(4.4.112) $$\log \frac{\pi}{2} = \int_0^1 \frac{(y-1)}{(1+y)\log y} dy$$

Rigour can be provided by initially considering the following finite sum (I have annotated the preceding equation numbers with an "a")

(4.4.105a) $$S_N = \sum_{n=1}^{N} (-1)^{n+1} \log \frac{n+1}{n} = \sum_{n=1}^{N} (-1)^{n+1} \int_0^1 \frac{dx}{n+x}$$

Using the elementary integral

(4.4.106a) $$\frac{1}{n+x} = \int_0^1 y^{n+x-1} dy$$

we deduce that

(4.4.107a) $$S_N = \sum_{n=1}^{N} (-1)^{n+1} \int_0^1 dx \int_0^1 y^{n+x-1} dy$$

Applying the geometric series we obtain

(4.4.108a) $$S_N = \int_0^1 dx \int_0^1 \frac{y^x \left[1 + (-y)^{N+1}\right]}{1+y} dy$$

Interchanging the order of integration gives

(4.4.109a) $$= \int_0^1 \frac{\left[1 + (-y)^{N+1}\right]}{1+y} dy \int_0^1 y^x dx$$

(4.4.110a) $$= \int_0^1 \frac{\left[1 + (-y)^{N+1}\right]}{1+y} dy \left\{ \frac{y^x}{\log y}\bigg|_0^1 \right\}$$

(4.4.111a) $$= \int_0^1 \frac{(y-1)\left[1 + (-y)^{N+1}\right]}{(1+y) \log y} dy$$

$$= \int_0^1 \frac{(y-1)}{(1+y) \log y} dy + (-1)^{N+1} \int_0^1 \frac{(y-1) y^{N+1}}{(1+y) \log y} dy$$

We now consider the remainder term

$$R_N = \int_0^1 \frac{y(y-1) y^N}{(1+y) \log y} dy$$

and note that

$$|R_N| = \left| \int_0^1 \frac{y(y-1) y^N}{(1+y) \log y} dy \right| \leq \int_0^1 \left| \frac{y(y-1) y^N}{(1+y) \log y} \right| dy \leq M \int_0^1 |y^N| dy = \frac{M}{N+1}$$

where $M = \sup_{y \in [0,1]} \left\{ \dfrac{y(y-1)}{(1+y)\log y} \right\}$

Using L'Hôpital's rule, it is easily seen that $f(y) = \dfrac{y(y-1)}{(1+y)\log y}$ is finite at $y = 0$ and $y = 1$ and hence is bounded throughout the interval $[0,1]$. Therefore, we have $\lim_{N \to \infty} R_N = 0$.

Accordingly we have as before

(4.4.112a) $\qquad \log \dfrac{\pi}{2} = \displaystyle\int_0^1 \dfrac{(y-1)}{(1+y)\log y} dy$

Three other proofs of this result are given in Sondow's paper [123a] "A faster product for $\pi$ and a new integral for $\log \dfrac{\pi}{2}$". This formula is also employed in Theorem 4.8. In June 2006, whilst looking for inspiration, I discovered that a derivation of this integral was also recorded in 1864 in Bertrand's treatise [22, Book 1, p.149], albeit with a trivial sign error being made in the last step of the proof. Bertrand showed that

(4.4.112ai) $\qquad \displaystyle\int_0^1 \dfrac{y^\alpha - y^\beta}{y(1+y^n)\log y} dy = \log\left[ \dfrac{\alpha}{\beta} \cdot \dfrac{\beta+n}{\alpha+n} \cdot \dfrac{\alpha+2n}{\beta+2n} \cdots \right]$

and (4.4.112a) follows by letting $\alpha = 1$, $\beta = 2$ and $n = 1$ and then applying the Wallis identity. Similarly we find that

(4.4.112aii) $\qquad \displaystyle\int_0^1 \dfrac{1-y^2}{y \log y} dy = \log\left[ \dfrac{1}{4} \cdot \dfrac{4+2}{1+2} \cdot \dfrac{1+4}{4+4} \cdots \right]$

and in connection with (4.4.24e) we also see that

(4.4.112aiii) $\qquad \displaystyle\int_0^1 \dfrac{y^{\alpha-1} - y^{\beta-1}}{(1+y)\log y} dy = \log\left[ \dfrac{\alpha}{\beta} \cdot \dfrac{\beta+1}{\alpha+1} \cdot \dfrac{\alpha+2}{\beta+2} \cdots \right]$

More generally, and using the same modus operandi, it is easily shown that

(4.4.112b) $\qquad \displaystyle\sum_{n=1}^\infty (-1)^{n+1} t^{n+1} \log \dfrac{n+1}{n} = \int_0^1 \dfrac{dy}{1+ty} \int_0^1 (ty)^x dx = \int_0^1 \dfrac{t(ty-1)}{(1+ty)\log ty} dy$

Letting $x = ty$ we may write this as

(4.4.112bi) $$\sum_{n=1}^{\infty}(-1)^{n+1}t^{n+1}\log\frac{n+1}{n}=\int_{0}^{t}\frac{x-1}{(1+x)\log x}dx$$

In 2006 Sondow and Hadjicostas [123ab] considered the generalised Euler constant function $\gamma(x)$ which they defined as (see also (E.43b) in Volume VI)

$$\gamma(t)=\sum_{n=1}^{\infty}t^{n-1}\left[\frac{1}{n}-\log\left(\frac{n+1}{n}\right)\right]$$

$$=-\frac{\log(1-t)}{t}-\frac{1}{t^2}\sum_{n=1}^{\infty}t^{n+1}\log\left(\frac{n+1}{n}\right)$$

Hence we have

$$\gamma(-t)=\frac{\log(1+t)}{t}-\frac{1}{t^2}\sum_{n=1}^{\infty}(-1)^{n+1}t^{n+1}\log\left(\frac{n+1}{n}\right)$$

(4.4.112bii) $$=\frac{\log(1+t)}{t}-\frac{1}{t^2}\int_{0}^{1}\frac{t(ty-1)}{(1+ty)\log ty}dy$$

and with $t=1$ we obtain

$$\gamma(-1)=\log 2-\int_{0}^{1}\frac{y-1}{(1+y)\log y}dy$$

Using (4.4.112b) this becomes

(4.4.112biii) $\gamma(-1)=\log 2-\log(\pi/2)=\log(4/\pi)$

where in the final part we used (4.4.102).

Differentiating (4.4.112bii) gives us

$$-\gamma'(-t)=\frac{1}{t(1+t)}-\frac{\log(1+t)}{t^2}+\frac{2}{t^3}\int_{0}^{1}\frac{t(ty-1)}{(1+ty)\log ty}dy$$

$$-\frac{1}{t^2}\int_{0}^{1}\frac{(1+ty)(2ty-1)\log(ty)-t(ty-1)[(1+ty)/t+y\log(ty)]}{[(1+ty)\log(ty)]^2}dy$$

With $t=1$ we get

$$-\gamma'(-1) = \frac{1}{2} - \log 2 + 2\int_0^1 \frac{y-1}{(1+y)\log y} dy - \int_0^1 \frac{y^2+2y-1}{(1+y)^2 \log y} dy + \int_0^1 \frac{y^2-1}{(1+y)^2 \log^2 y} dy$$

$$\gamma'(-1) = -\frac{1}{2} + \log 2 - 2\log\frac{\pi}{2} + \int_0^1 \frac{y^2+2y-1}{(1+y)^2 \log y} dy - \int_0^1 \frac{y^2-1}{(1+y)^2 \log^2 y} dy$$

Sondow and Hadjicostas [123ab] have shown that

(4.4.112biv) $\qquad \gamma'(-1) = \log\dfrac{2^{11/6} A^6}{\pi^{3/2} e}$

and

(4.4.112bv) $\qquad \gamma''(-1) = \log\dfrac{2^{10/3} A^{24}}{\pi^4 e^{13/4}} - \dfrac{7\varsigma(3)}{2\pi^2}$

where $A$ is the Glaisher-Kinkelin constant.

Adamchik [2a] proved that

$$\int_0^1 \frac{x^{p-1}}{1+x^n} \log\log\left(\frac{1}{x}\right) dx = \frac{\gamma+\log(2n)}{2n}\left[\psi\left(\frac{p}{2n}\right) - \psi\left(\frac{n+p}{2n}\right)\right] + \frac{1}{2n}\left[\varsigma'\left(1,\frac{p}{2n}\right) - \varsigma'\left(1,\frac{n+p}{2n}\right)\right]$$

where $\operatorname{Re}(p) > 0$, $\operatorname{Re}(n) > 0$ (a proof appears in (C.57)). He also showed that

$$\int_0^1 \frac{x^{p-1}}{(1+x^n)^2} \log\log\left(\frac{1}{x}\right) dx = \frac{(n-p)[\gamma+\log(2n)]}{2n^2}\left[\psi\left(\frac{p}{2n}\right) - \psi\left(\frac{n+p}{2n}\right)\right]$$

$$-\frac{1}{2n}\left[\gamma+\log(2n) - 2\log\frac{\Gamma\left(\frac{p}{2n}\right)}{\Gamma\left(\frac{n+p}{2n}\right)}\right] + \frac{n-p}{2n^2}\left[\varsigma'\left(1,\frac{p}{2n}\right) - \varsigma'\left(1,\frac{n+p}{2n}\right)\right]$$

and with $n=1$ these become

$$\int_0^1 \frac{x^{p-1}}{1+x} \log\log\left(\frac{1}{x}\right) dx = \frac{\gamma+\log 2}{2}\left[\psi\left(\frac{p}{2}\right) - \psi\left(\frac{1+p}{2}\right)\right] + \frac{1}{2}\left[\varsigma'\left(1,\frac{p}{2}\right) - \varsigma'\left(1,\frac{1+p}{2}\right)\right]$$

$$\int_0^1 \frac{x^{p-1}}{(1+x)^2} \log\log\left(\frac{1}{x}\right) dx = \frac{(1-p)[\gamma+\log 2]}{2}\left[\psi\left(\frac{p}{2}\right) - \psi\left(\frac{1+p}{2}\right)\right]$$

$$-\frac{1}{2}\left[\gamma + \log 2 - 2\log \frac{\Gamma\left(\frac{p}{2}\right)}{\Gamma\left(\frac{1+p}{2}\right)}\right] + \frac{1-p}{2}\left[\varsigma'\left(1, \frac{p}{2}\right) - \varsigma'\left(1, \frac{1+p}{2}\right)\right]$$

Using Adamchik's formulae may enable us to evaluate the integral $\int_0^1 \frac{y^2 - 1}{(1+y)^2 \log^2 y} dy$.

In (4.4.112niv) we will see that

$$\int_0^1 \frac{t^{u-1} x(1-t)}{[1 - x(1-t)] \log t} dt = \sum_{n=1}^{\infty} x^n \sum_{k=0}^{n} \binom{n}{k} (-1)^k \log(k + u)$$

and with $u = 1$ we get

$$\int_0^1 \frac{x(1-t)}{[1 - x(1-t)] \log t} dt = \sum_{n=1}^{\infty} x^n \sum_{k=0}^{n} \binom{n}{k} (-1)^k \log(k + 1)$$

Sondow and Hadjicostas [123ab] have also shown that for $|x| < 1$

$$\frac{x}{1-x} \gamma\left(\frac{-x}{1-x}\right) + \log(1-x) = \sum_{n=1}^{\infty} x^n \sum_{k=0}^{n} \binom{n}{k} (-1)^k \log(k + 1)$$

and hence we obtain

$$\frac{x}{1-x} \gamma\left(\frac{-x}{1-x}\right) + \log(1-x) = \int_0^1 \frac{x(1-t)}{[1 - x(1-t)] \log t} dt$$

In [123a] Sondow and Guillera noted that

$$\frac{1}{x} \frac{\partial}{\partial s} Li_s\left(\frac{x}{x-1}\right)\bigg|_{s=0} = \int_0^1 \frac{x(1-t)}{[1 - x(1-t)] \log t} dt$$

and hence we have

$$\frac{x}{1-x} \gamma\left(\frac{-x}{1-x}\right) + \log(1-x) = \frac{1}{x} \frac{\partial}{\partial s} Li_s\left(\frac{-x}{1-x}\right)\bigg|_{s=0}$$

Letting $x = 1/2$ results in

$$\gamma(-1) - \log 2 = 2 \frac{\partial}{\partial s} Li_s(-1)\bigg|_{s=0}$$

and since $Li_s(-1) = -\varsigma_a(s)$ we have

$$\gamma(-1) - \log 2 = -2\varsigma_a'(0)$$

We have

$$\varsigma_a'(s) = (1 - 2^{1-s})\varsigma'(s) + 2^{1-s}\varsigma(s)\log 2$$

and thus

$$\varsigma_a'(0) = -\varsigma'(0) - \log 2 = \frac{1}{2}\log(2\pi) - \log 2$$

This then implies that

$$\gamma(-1) - \log 2 = -\log(2\pi) + 2\log 2$$

and hence we have

$$\gamma(-1) = \log\left(\frac{4}{\pi}\right)$$

as previously determined by Sondow and Hadjicostas [123ab].

Letting $t = \frac{-x}{1-x}$ gives us

$$-t\gamma(t) - \log(1-t) = \frac{t-1}{t}\frac{\partial}{\partial s}Li_s(t)\bigg|_{s=0}$$

or alternatively

$$t^2\gamma(t) + t\log(1-t) = (1-t)\frac{\partial}{\partial s}Li_s(t)\bigg|_{s=0}$$

in agreement with the result obtained by Sondow and Hadjicostas [123ab].

Let us now differentiate the above identity

$$t^2\gamma'(t) + 2t\gamma(t) - \frac{t}{1-t} + \log(1-t) = (1-t)\frac{\partial}{\partial t}\frac{\partial}{\partial s}Li_s(t) - \frac{\partial}{\partial s}Li_s(t)\bigg|_{s=0}$$

We have

$$\frac{\partial}{\partial t}\frac{\partial}{\partial s}Li_s(t) = \frac{\partial}{\partial s}\frac{\partial}{\partial t}Li_s(t) = \frac{1}{t}\frac{\partial}{\partial s}Li_{s-1}(t)$$

but further consideration will have to await another day.

□

The following integral is well-known (see for example [126a, p.173]): a proof was also given in (4.4.24c).

(4.4.112c) $$\int_0^1 \frac{y^b - y^a}{\log y}dy = \log\frac{b+1}{a+1}$$

and hence we have

(4.4.112d) $$\int_0^1 \frac{y-1}{\log y}dy = \log 2$$

Therefore, using (4.4.112) we have

$$\log \pi = \int_0^1 \frac{(y-1)}{(1+y)\log y}dy + \log 2$$

$$= \int_0^1 \frac{(y-1)}{(1+y)\log y}dy + \int_0^1 \frac{y-1}{\log y}dy$$

(4.4.112e) $$= \int_0^1 \frac{(y+2)(y-1)}{(y+1)\log y}dy$$

Integrating (4.4.112c) with respect to $b$ we obtain

$$\int_0^1 \left[\frac{y^u - 1}{\log^2 y} - \frac{uy^a}{\log y}\right]dy = (u+1)\log(u+1) - u - u\log(a+1)$$

By using the binomial theorem we have (as mentioned by Adamchik in [2a])

$$\lambda \int_0^z \frac{y^{p-1}}{\lambda + y^n} dy = \sum_{k=0}^{\infty} \frac{(-1)^k}{\lambda^k} \int_0^z y^{nk+p-1} dy \quad \text{(provided } \left|\frac{y^n}{\lambda}\right| < 1\text{)}$$

$$= \sum_{k=0}^{\infty} \frac{(-1)^k z^{nk+p}}{\lambda^k (nk+p)}$$

With $z = n = \lambda = 1$ and $p - 1 = x$ we have

$$\int_0^1 \frac{y^x}{1+y} dy = \sum_{k=0}^{\infty} \frac{(-1)^k}{(k+1+x)}$$

Integrating the above with respect to $x$ we have

$$\int_0^t dx \int_0^1 \frac{y^x}{1+y} dy = \sum_{k=0}^{\infty} \int_0^t \frac{(-1)^k}{(k+1+x)} dx$$

$$= \sum_{k=0}^{\infty} (-1)^k \log \frac{k+1+t}{k+1}$$

$$= \sum_{k=1}^{\infty} (-1)^{k+1} \log \frac{k+t}{k}$$

The next identity is reported as an exercise in Whittaker & Watson [135, p.262]

(4.4.112f) $$\int_0^1 \frac{y^x}{1+y} dy = \frac{1}{2}\left[\psi\left(\frac{x+2}{2}\right) - \psi\left(\frac{x+1}{2}\right)\right]$$

and a proof is contained in (C.56) of Volume VI of this series of papers.

By definition we have

$$\frac{1}{2}\psi\left(\frac{x+2}{2}\right) = \frac{d}{dx} \log \Gamma\left(\frac{x+2}{2}\right)$$

and hence for the integral in (4.4.108) we have

$$S = \int_0^1 \frac{1}{2}\left[\psi\left(\frac{x+2}{2}\right) - \psi\left(\frac{x+1}{2}\right)\right]dx = \log\Gamma\left(\frac{x+2}{2}\right) - \log\Gamma\left(\frac{x+1}{2}\right)\Big|_0^1$$

$$= \log\Gamma(3/2) + \log\Gamma(1/2)$$

Using the functional equation for $\Gamma(x)$ we have $\Gamma(3/2) = \frac{1}{2}\Gamma(1/2) = \frac{\sqrt{\pi}}{2}$ and hence we have a further proof that $S = \log\frac{\pi}{2}$.

(4.4.112fa) $$\frac{1}{2}\int_0^1\left[\psi\left(\frac{x+2}{2}\right) - \psi\left(\frac{x+1}{2}\right)\right]dx = \log\frac{\pi}{2}$$

As mentioned previously, it transpires that Sondow's integral (4.4.112) is not new: in addition to Bertrand's proof, it is reported in Whittaker & Watson [135, p.262] that Kummer discovered the following identity (for $\alpha > 0$, $\beta > 0$)

(4.4.112g) $$\int_0^1 \frac{t^{\alpha-1} - t^{\beta-1}}{(1+t)\log t}dt = \log\frac{\Gamma\left(\frac{1+\alpha}{2}\right)\Gamma\left(\frac{\beta}{2}\right)}{\Gamma\left(\frac{1+\beta}{2}\right)\Gamma\left(\frac{\alpha}{2}\right)}$$

(this integral is also contained in G&R [74, p.541]). Therefore, for $\alpha = 2$ and $\beta = 1$ we have

$$\int_0^1 \frac{t-1}{(1+t)\log t}dt = \log\frac{\Gamma(3/2)\Gamma(1/2)}{\Gamma(1)\Gamma(1)} = \log\frac{\pi}{2}$$

and this particular integral is also contained in G&R [74, p.540]. Reference should also be made to (3.86h).

From (4.4.112aiii) we see that

$$\int_0^1 \frac{y^{\alpha-1} - y^{\beta-1}}{(1+y)\log y}dy = \log\left[\frac{\alpha}{\beta}\cdot\frac{\beta+1}{\alpha+1}\cdot\frac{\alpha+2}{\beta+2}\cdots\right]$$

and hence we obtain

(4.4.112ga) $$\frac{\Gamma\left(\frac{1+\alpha}{2}\right)\Gamma\left(\frac{\beta}{2}\right)}{\Gamma\left(\frac{1+\beta}{2}\right)\Gamma\left(\frac{\alpha}{2}\right)} = \frac{\alpha}{\beta} \cdot \frac{\beta+1}{\alpha+1} \cdot \frac{\alpha+2}{\beta+2} \cdots$$

Due to the explosive growth in mathematics, especially during the last two centuries, it is very difficult to believe, with any certitude, that one has discovered a new result: accordingly, I am no longer convinced that many of the few results, which I initially thought were new, are actually new (it's simply that I haven't seen them before!).

Let us now use the identity (4.4.112f)

$$\int_0^1 \frac{y^{x-1}}{1+y} dy = \frac{1}{2}\left[\psi\left(\frac{x+1}{2}\right) - \psi\left(\frac{x}{2}\right)\right]$$

Integrating again we have

$$\int_\beta^\alpha dx \int_0^1 \frac{y^{x-1}}{1+y} dy = \frac{1}{2}\int_\beta^\alpha \left[\psi\left(\frac{x+1}{2}\right) - \psi\left(\frac{x}{2}\right)\right] dx$$

Changing the order of integration we have

$$\int_\beta^\alpha dx \int_0^1 \frac{y^{x-1}}{1+y} dy = \int_0^1 \frac{dy}{1+y} \int_\beta^\alpha y^{x-1} dx$$

$$= \int_0^1 \frac{y^{\alpha-1} - y^{\beta-1}}{(1+y)\log y} dy$$

We have by direct integration

$$\frac{1}{2}\int_\beta^\alpha \left[\psi\left(\frac{x+1}{2}\right) - \psi\left(\frac{x}{2}\right)\right] dx = \log \frac{\Gamma\left(\frac{1+\alpha}{2}\right)\Gamma\left(\frac{\beta}{2}\right)}{\Gamma\left(\frac{1+\beta}{2}\right)\Gamma\left(\frac{\alpha}{2}\right)}$$

and we have therefore proved Kummer's identity (4.4.112g)

$$\int_0^1 \frac{t^{\alpha-1} - t^{\beta-1}}{(1+t)\log t} dt = \log \frac{\Gamma\left(\frac{1+\alpha}{2}\right)\Gamma\left(\frac{\beta}{2}\right)}{\Gamma\left(\frac{1+\beta}{2}\right)\Gamma\left(\frac{\alpha}{2}\right)}$$

In (4.4.24e) we showed that

$$(4.4.112\text{gb}) \qquad \int_0^1 \frac{(t^{\alpha-1} - t^{\beta-1})}{(1+t)\log t} dt = \sum_{n=0}^{\infty} \frac{1}{2^{n+1}} \sum_{k=0}^{n} \binom{n}{k} (-1)^k \log \frac{k+\alpha}{k+\beta}$$

and hence we obtain

$$(4.4.112\text{h}) \qquad \sum_{n=0}^{\infty} \frac{1}{2^{n+1}} \sum_{k=0}^{n} \binom{n}{k} (-1)^k \log \frac{k+\alpha}{k+\beta} = \log \frac{\Gamma\left(\frac{1+\alpha}{2}\right)\Gamma\left(\frac{\beta}{2}\right)}{\Gamma\left(\frac{1+\beta}{2}\right)\Gamma\left(\frac{\alpha}{2}\right)}$$

Therefore we have

$$\sum_{n=0}^{\infty} \frac{1}{2^{n+1}} \sum_{k=0}^{n} \binom{n}{k} (-1)^k \log \frac{k+2}{k+1} = \log \frac{\Gamma\left(\frac{3}{2}\right)\Gamma\left(\frac{1}{2}\right)}{\Gamma(1)\Gamma(1)} = \log \frac{\pi}{2}$$

In (4.4.102) we showed that

$$\log \frac{\pi}{2} = \sum_{n=1}^{\infty} (-1)^{n+1} \log \frac{n+1}{n}$$

Euler's work [90, p.244] on the transformation of series shows that

$$(4.4.112\text{i}) \qquad \sum_{n=1}^{\infty} (-1)^{n+1} a_n = \sum_{n=0}^{\infty} \frac{1}{2^{n+1}} \sum_{k=0}^{n} \binom{n}{k} (-1)^k a_{k+1}$$

is valid for any convergent series of complex numbers, and Sondow [123a] showed that this series acceleration technique results in

$$(4.4.112\text{j}) \qquad \log \frac{\pi}{2} = \sum_{n=1}^{\infty} (-1)^{n+1} \log \frac{n+1}{n} = \sum_{n=0}^{\infty} \frac{1}{2^{n+1}} \sum_{k=0}^{n} \binom{n}{k} (-1)^k \log \frac{k+2}{k+1}$$

Since from (4.4.113)

$$\log \frac{\pi}{2} = \sum_{n=1}^{\infty} \frac{1}{2^n} \sum_{k=0}^{n} \binom{n}{k} (-1)^{k+1} \log(1+k) = -2 \sum_{n=0}^{\infty} \frac{1}{2^{n+1}} \sum_{k=0}^{n} \binom{n}{k} (-1)^k \log(1+k)$$

we easily see that

$$(4.4.112k) \quad \sum_{n=0}^{\infty} \frac{1}{2^{n+1}} \sum_{k=0}^{n} \binom{n}{k} (-1)^k \log(k+2) = -\sum_{n=0}^{\infty} \frac{1}{2^{n+1}} \sum_{k=0}^{n} \binom{n}{k} (-1)^k \log(k+1) = \frac{1}{2} \log \frac{\pi}{2}$$

In view of the symmetry involved in (4.4.112h) we have

$$\sum_{n=0}^{\infty} \frac{1}{2^{n+1}} \sum_{k=0}^{n} \binom{n}{k} (-1)^k \log(k+x) = \log \frac{\Gamma\left(\frac{1+x}{2}\right)}{\Gamma\left(\frac{x}{2}\right)} + c$$

The constant $c$ may be determined from (4.4.112k) or from (4.4.113) by letting $x=1$ and, using our old faithful identity (6.60) $\Gamma\left(\frac{1}{2}\right) = \sqrt{\pi}$, we get

$$(4.4.112l) \quad \sum_{n=0}^{\infty} \frac{1}{2^{n+1}} \sum_{k=0}^{n} \binom{n}{k} (-1)^k \log(k+x) = \log \frac{\Gamma\left(\frac{1+x}{2}\right)}{\Gamma\left(\frac{x}{2}\right)} + \log \sqrt{2}$$

Differentiating equation (4.4.112l) we get

$$(4.4.112m) \quad \sum_{n=0}^{\infty} \frac{1}{2^{n+1}} \sum_{k=0}^{n} \binom{n}{k} \frac{(-1)^k}{(k+x)} = \frac{1}{2}\left[\psi\left(\frac{1+x}{2}\right) - \psi\left(\frac{x}{2}\right)\right]$$

and for $s \geq 2$ we have

$$(4.4.112n) \quad (-1)^{s+1}(s-1)! \sum_{n=0}^{\infty} \frac{1}{2^{n+1}} \sum_{k=0}^{n} \binom{n}{k} \frac{(-1)^k}{(k+x)^s} = \frac{1}{2^s}\left[\psi^{(s-1)}\left(\frac{1+x}{2}\right) - \psi^{(s-1)}\left(\frac{x}{2}\right)\right]$$

Referring to (4.4.79) we see a connection with the Hurwitz-Lerch zeta function

$$(4.4.112\text{ni}) \quad (-1)^{s+1}(s-1)! \Phi(-1, s, x) = \frac{1}{2^s}\left[\psi^{(s-1)}\left(\frac{1+x}{2}\right) - \psi^{(s-1)}\left(\frac{x}{2}\right)\right]$$

Letting $x=1$ in (4.4.112j) we obtain

$$\sum_{n=0}^{\infty} \frac{1}{2^{n+1}} \sum_{k=0}^{n} \binom{n}{k} \frac{(-1)^k}{(k+1)} = \frac{1}{2}\left[\psi(1) - \psi\left(\frac{1}{2}\right)\right]$$

Since $\psi(1) = -\gamma$ and $\psi\left(\dfrac{1}{2}\right) = -\gamma - 2\log 2$ we have

$$\sum_{n=0}^{\infty} \frac{1}{2^{n+1}} \sum_{k=0}^{n} \binom{n}{k} \frac{(-1)^k}{k+1} = \log 2$$

or equivalently, using (3.11), we may express this as

$$\varsigma_a(1) = \log 2$$

It is well known that [1, 6.4.4] for $k \geq 1$

(4.4.112nii) $\qquad \psi^{(k)}(1/2) = (-1)^{k+1} k! (2^{k+1} - 1) \varsigma(k+1)$

and hence for $s \geq 2$ we have

$$\psi^{(s-1)}(1/2) = (-1)^s (s-1)! (2^s - 1) \varsigma(s)$$

We also have

$$\psi^{(s-1)}(1) = (-1)^s (s-1)! \varsigma(s)$$

Therefore, with $x = 1$ in (4.4.112n) we obtain the Hasse/Sondow identity (3.11).

With $x = 2$ in (4.4.112n) we obtain

$$(-1)^{s+1}(s-1)! \sum_{n=0}^{\infty} \frac{1}{2^{n+1}} \sum_{k=0}^{n} \binom{n}{k} \frac{(-1)^k}{(k+2)^s} = \frac{1}{2^s} \left[ \psi^{(s-1)}(3/2) - \psi^{(s-1)}(1) \right]$$

We have from [126, p.22]

$$\psi^{(s-1)}(3/2) = \psi^{(s-1)}(1 + 1/2) = \psi^{(s-1)}(1/2) + (-1)^{s-1}(s-1)! 2^s$$

and hence

$$\psi^{(s-1)}(3/2) - \psi^{(s-1)}(1) = (-1)^s (s-1)!(2^s - 2)\varsigma(s) + (-1)^{s-1}(s-1)! 2^s$$

giving us

(4.4.112niii) $\qquad \displaystyle\sum_{n=0}^{\infty} \frac{1}{2^{n+1}} \sum_{k=0}^{n} \binom{n}{k} \frac{(-1)^k}{(k+2)^s} = 1 - \varsigma_a(s)$

In this connection, see also (4.4.112z(i)).

We have seen in (4.4.12) of Volume II(b) that

$$\int_0^1 t^{x-1}(1-yt)^n \, dt = \sum_{k=0}^{n} \binom{n}{k}(-1)^k \frac{y^k}{k+x}$$

and hence

$$\int_0^1 t^{x-1}(1-t)^n \, dt = \sum_{k=0}^{n} \binom{n}{k}\frac{(-1)^k}{k+x} = B(x, n+1)$$

Integration results in

$$\int_0^1 \frac{(t^{u-1}-1)(1-t)^n}{\log t} \, dt = \sum_{k=0}^{n} \binom{n}{k}(-1)^k \log \frac{k+u}{k+1}$$

and having regard to (4.4.9a) we then see that

$$\int_0^1 \frac{t^{u-1}(1-t)^n}{\log t} \, dt = \sum_{k=0}^{n} \binom{n}{k}(-1)^k \log(k+u)$$

Completing the summation we get

$$\sum_{n=1}^{\infty} x^n \int_0^1 \frac{t^{u-1}(1-t)^n}{\log t} \, dt = \sum_{n=1}^{\infty} x^n \sum_{k=0}^{n} \binom{n}{k}(-1)^k \log(k+u)$$

and hence we have for $|x(1-t)| < 1$

(4.4.112niv) $$\int_0^1 \frac{t^{u-1}x(1-t)}{[1-x(1-t)]\log t} \, dt = \sum_{n=1}^{\infty} x^n \sum_{k=0}^{n} \binom{n}{k}(-1)^k \log(k+u)$$

With $x = 1/2$ we get

$$\int_0^1 \frac{t^{u-1}(1-t)}{(1+t)\log t} \, dt = \sum_{n=1}^{\infty} \frac{1}{2^n} \sum_{k=0}^{n} \binom{n}{k}(-1)^k \log(k+u)$$

which corresponds with (4.4.100u) in the case where $u = 1$.

Integrating (4.4.112niv) gives us

(4.4.112nv) $$-\int_0^1 \frac{t^{u-1}}{\log t}\left\{w + \frac{\log[1-w(1-t)]}{1-t}\right\} dt = \sum_{n=1}^{\infty} \frac{w^{n+1}}{n+1}\sum_{k=0}^{n}\binom{n}{k}(-1)^k \log(k+u)$$

Letting $w = 1$ we get

(4.4.112nvi) $$-\int_0^1 t^{u-1}\left[\frac{1}{\log t} + \frac{1}{1-t}\right] dt = \sum_{n=1}^{\infty} \frac{1}{n+1}\sum_{k=0}^{n}\binom{n}{k}(-1)^k \log(k+u)$$

$$= \psi(u) - \log u$$

and with $u = 1$ we regain Sondow's formula (4.4.92a)

$$\gamma = \int_0^1 \left[\frac{1}{\log t} + \frac{1}{1-t}\right] dt = \sum_{n=1}^{\infty} \frac{1}{n+1}\sum_{k=0}^{n}\binom{n}{k}(-1)^{k+1} \log(k+1)$$

Differentiating (4.4.112nvi) results in for $p \geq 2$

(4.4.112nvii) $$\int_0^1 t^{u-1} \log^{p-1} t \left[\frac{1}{\log t} + \frac{1}{1-t}\right] dt = (-1)^{p+1}(p-2)!\sum_{n=1}^{\infty} \frac{1}{n+1}\sum_{k=0}^{n}\binom{n}{k}\frac{(-1)^{k+1}}{(k+u)^{p-1}}$$

or alternatively for $p \geq 1$

$$\int_0^1 t^{u-1} \log^p (1/t) \left[\frac{1}{\log t} + \frac{1}{1-t}\right] dt = \Gamma(p)\sum_{n=1}^{\infty} \frac{1}{n+1}\sum_{k=0}^{n}\binom{n}{k}\frac{(-1)^k}{(k+u)^p}$$

We recall the Hasse formula (3.12a) from Volume I

$$\varsigma(p,u) = \sum_{n=0}^{\infty} \frac{1}{(n+u)^p} = \frac{1}{p-1}\sum_{n=0}^{\infty} \frac{1}{n+1}\sum_{k=0}^{n}\binom{n}{k}\frac{(-1)^k}{(k+u)^{p-1}}$$

$$(p-1)\varsigma(p,u) - \frac{1}{u^{p-1}} = \sum_{n=1}^{\infty} \frac{1}{n+1}\sum_{k=0}^{n}\binom{n}{k}\frac{(-1)^k}{(k+u)^{p-1}}$$

and we therefore have

(4.4.112nviii) $$\int_0^1 t^{u-1} \log^{p-1} t \left[\frac{1}{\log t} + \frac{1}{1-t}\right] dt = (-1)^{p+1}(p-1)!\left(\varsigma(p,u) - \frac{1}{(p-1)u^{p-1}}\right)$$

and with $u = 1$

(4.4.112nix) $$\int_0^1 \log^{p-1} t \left[\frac{1}{\log t} + \frac{1}{1-t}\right] dt = (-1)^{p+1}(p-1)!\left(\varsigma(p) - \frac{1}{p-1}\right)$$

Hence in the limit as $p \to 1$ we have the well-known result

$$\int_0^1 \left[\frac{1}{\log t} + \frac{1}{1-t}\right] dt = \lim_{p \to 1}\left(\varsigma(p) - \frac{1}{p-1}\right) = \gamma$$

Differentiation results in

$$\int_0^1 t^{u-1} \log^p(1/t) \log\log(1/t) \left[\frac{1}{\log t} + \frac{1}{1-t}\right] dt$$

$$= -\Gamma(p) \sum_{n=1}^{\infty} \frac{1}{n+1} \sum_{k=0}^{n} \binom{n}{k} \frac{(-1)^k \log(k+u)}{(k+u)^p} + \Gamma'(p) \sum_{n=1}^{\infty} \frac{1}{n+1} \sum_{k=0}^{n} \binom{n}{k} \frac{(-1)^k}{(k+u)^p}$$

Substituting $t = e^{-x}$ in (4.4.112nviii) gives us

$$\int_0^1 t^{u-1} \log^{p-1} t \left[\frac{1}{\log t} + \frac{1}{1-t}\right] dt = (-1)^{p+1} \int_0^{\infty} \frac{x^{p-1} e^{-px}}{1 - e^{-x}} dx - (-1)^{p+1} \int_0^{\infty} x^{p-2} e^{-px} dx$$

$$= (-1)^{p+1}(p-1)!\left(\varsigma(p,u) - \frac{1}{(p-1)u^{p-1}}\right)$$

$\square$

From (4.4.112l) we deduce that

$$\frac{1}{2}\int_0^1 \frac{t^{u-1}(1-t)}{(1+t)\log t} dt = \log \frac{\Gamma\left(\frac{1+u}{2}\right)}{\Gamma\left(\frac{u}{2}\right)} + \frac{1}{2}\log 2$$

and integration results in

$$\frac{1}{2}\int_a^z du \int_0^1 \frac{t^{u-1}(1-t)}{(1+t)\log t} dt = \int_a^z \log \frac{\Gamma\left(\frac{1+u}{2}\right)}{\Gamma\left(\frac{u}{2}\right)} du + \frac{1}{2}(z-a)\log 2$$

Therefore we obtain for $a, z > 0$

$$\int_0^1 \frac{(t^{z-1} - t^{a-1})(1-t)}{(1+t)\log^2 t} dt = 2\int_a^z \log \frac{\Gamma\left(\frac{1+u}{2}\right)}{\Gamma\left(\frac{u}{2}\right)} du + (z-a)\log 2$$

and the integral on the right-hand side may be evaluated using (4.4.112p).

□

As noted by Coffey [45h] in a recent paper, we may write for $m \geq 2$

$$\frac{d^{m-2}}{dx^{m-2}} \sum_{k=0}^n \binom{n}{k} \frac{(-1)^k y^k}{k+x} = \frac{d^{m-2}}{dx^{m-2}} \sum_{k=0}^n \binom{n}{k} (-1)^k y^k \int_0^\infty e^{-(k+x)t} dt$$

$$= (-1)^m \sum_{k=0}^n \binom{n}{k} (-1)^k \int_0^\infty t^{m-2} y^k e^{-(k+x)t} dt$$

$$= (-1)^m \int_0^\infty t^{m-2} e^{-xt} (1 - ye^{-t})^n dt$$

where the interchange of differentiation and integration is justified by the absolute convergence of the above integral.

Hence we have

(4.4.112nx)
$$(m-2)! \sum_{k=0}^n \binom{n}{k} \frac{(-1)^k y^k}{(k+x)^{m-1}} = \int_0^\infty t^{m-2} e^{-xt} (1 - ye^{-t})^n dt$$

and we have the summation

$$(m-2)! \sum_{n=0}^\infty u^n \sum_{k=0}^n \binom{n}{k} \frac{(-1)^k y^k}{(k+x)^{m-1}} = \sum_{n=0}^\infty u^n \int_0^\infty t^{m-2} e^{-xt} (1 - ye^{-t})^n dt$$

$$= \int_0^\infty t^{m-2} e^{-xt} \sum_{n=0}^\infty u^n (1 - ye^{-t})^n dt$$

$$= \int_0^\infty \frac{t^{m-2} e^{-xt}}{1-u(1-ye^{-t})} dt$$

We therefore have

$$(m-2)! \sum_{n=0}^\infty u^n \sum_{k=0}^n \binom{n}{k} \frac{(-1)^k y^k}{(k+x)^{m-1}} = \int_0^\infty \frac{t^{m-2} e^{-(x-1)t}}{e^t - u(e^t - y)} dt$$

Integration with respect to $u$ results in

$$(m-2)! \sum_{n=0}^\infty \frac{v^{n+1}}{n+1} \sum_{k=0}^n \binom{n}{k} \frac{(-1)^k y^k}{(k+x)^{m-1}} = \int_0^v du \int_0^\infty \frac{t^{m-2} e^{-(x-1)t}}{e^t - u(e^t - y)} dt$$

$$= -\int_0^\infty \frac{t^{m-2} e^{-(x-1)t} \left( \log\left[ e^t - v(e^t - y) \right] - \log\left[ e^t \right] \right)}{e^t - y} dt$$

$$= \int_0^\infty \frac{t^{m-1} e^{-(x-1)t}}{e^t - y} dt - \int_0^\infty \frac{t^{m-2} e^{-(x-1)t} \log\left[ e^t(1-v) + vy \right]}{e^t - y} dt$$

With $v = 1$ we get

(4.4.112ny)
$$(m-2)! \sum_{n=0}^\infty \frac{1}{n+1} \sum_{k=0}^n \binom{n}{k} \frac{(-1)^k y^k}{(k+x)^{m-1}} = \int_0^\infty \frac{t^{m-1} e^{-(x-1)t}}{e^t - y} dt - \log y \int_0^\infty \frac{t^{m-2} e^{-(x-1)t}}{e^t - y} dt$$

When $v = y = 1$ this becomes

$$(m-2)! \sum_{n=0}^\infty \frac{1}{n+1} \sum_{k=0}^n \binom{n}{k} \frac{(-1)^k}{(k+x)^{m-1}} = \int_0^\infty \frac{t^{m-1} e^{-(x-1)t}}{e^t - 1} dt$$

and with $x = 1$ we have

$$(m-2)! \sum_{n=0}^\infty \frac{1}{n+1} \sum_{k=0}^n \binom{n}{k} \frac{(-1)^k}{(k+1)^{m-1}} = \int_0^\infty \frac{t^{m-1}}{e^t - 1} dt$$

We have seen in (4.4.85) that

$$\frac{\log y}{y} Li_{s-1}(y) + \frac{s-1}{y} Li_s(y) = \sum_{n=0}^\infty \frac{1}{n+1} \sum_{k=0}^n \binom{n}{k} \frac{(-1)^k y^k}{(k+1)^{s-1}}$$

and we saw in (4.4.25) that

$$Li_s(y) = \frac{y}{\Gamma(s)} \int_0^\infty \frac{t^{s-1}}{e^t - y} dt$$

Thus we have

$$\frac{\log y}{y} Li_{s-1}(y) + \frac{s-1}{y} Li_s(y) = \frac{\log y}{y} \frac{y}{\Gamma(s-1)} \int_0^\infty \frac{t^{s-2}}{e^t - y} dt + \frac{s-1}{y} \frac{y}{\Gamma(s)} \int_0^\infty \frac{t^{s-1}}{e^t - y} dt$$

$$= \frac{\log y}{\Gamma(s-1)} \int_0^\infty \frac{t^{s-2}}{e^t - y} dt + \frac{1}{\Gamma(s-1)} \int_0^\infty \frac{t^{s-1}}{e^t - y} dt$$

but this does not agree with (4.4.112ny) because of the negative sign.

From (4.4.25) we have

$$\Phi(y, s, x) = \sum_{n=0}^\infty \frac{y^n}{(n+x)^s} = \frac{1}{\Gamma(s)} \int_0^\infty \frac{t^{s-1} e^{-(x-1)t}}{e^t - y} dt$$

and thus

$$\Phi(y, s, x+1) = \frac{1}{\Gamma(s)} \int_0^\infty \frac{t^{s-1} e^{-xt}}{e^t - y} dt$$

and we may note that

$$\Phi(y, s, x+1) = \frac{1}{y} \left[ \Phi(y, s, x) - \frac{1}{x^s} \right]$$

Referring to (4.4.112ny) (subject to the negative sign error)

$$(m-2)! \sum_{n=0}^\infty \frac{1}{n+1} \sum_{k=0}^n \binom{n}{k} \frac{(-1)^k y^k}{(k+x)^{m-1}} = \int_0^\infty \frac{t^{m-1} e^{-(x-1)t}}{e^t - y} dt - \log y \int_0^\infty \frac{t^{m-2} e^{-(x-1)t}}{e^t - y} dt$$

$$= \Gamma(m) \Phi(y, m, x) - \Gamma(m-1) \log y \, \Phi(y, m-1, x)$$

This may be written as

$$\sum_{n=0}^{\infty}\frac{1}{n+1}\sum_{k=0}^{n}\binom{n}{k}\frac{(-1)^k y^k}{(k+x)^{m-1}} = (m-1)\Phi(y,m,x) - \log y\,\Phi(y,m-1,x)$$

We also have for the polylogarithm

$$\frac{Li_m(y)}{y} = \Phi(y,m,1)$$

With $m=2$ in (4.4.112nx) we have

$$\sum_{n=0}^{\infty} u^n \sum_{k=0}^{n}\binom{n}{k}\frac{(-1)^k y^k}{k+x} = \int_0^{\infty}\frac{e^{-(x-1)t}}{e^t - u(e^t - y)}\,dt$$

and integration with respect to $x$ gives us

$$\sum_{n=0}^{\infty} u^n \sum_{k=0}^{n}\binom{n}{k}(-1)^k y^k \log\frac{k+b}{k+a} = \int_0^{\infty}\frac{e^{-(a-1)t} - e^{-(b-1)t}}{t\left[e^t - u(e^t - y)\right]}\,dt$$

Letting $u=1/2$ results in

$$\sum_{n=0}^{\infty}\frac{1}{2^{n+1}}\sum_{k=0}^{n}\binom{n}{k}(-1)^k y^k \log\frac{k+b}{k+a} = \int_0^{\infty}\frac{e^{-(a-1)t} - e^{-(b-1)t}}{t(e^t + y)}\,dt$$

and looking ahead we see that this is equivalent to (4.4.117ki) when $y=1$.

We have

$$\sum_{n=0}^{\infty}\frac{1}{2^{n+1}}\sum_{k=0}^{n}\binom{n}{k}(-1)^k y^k \log\frac{k+b}{k+a} = \int_0^{\infty}\frac{t^{-1}e^{-(a-1)t}}{e^t + y}\,dt - \int_0^{\infty}\frac{t^{-1}e^{-(b-1)t}}{e^t + y}\,dt$$

and this suggests that

$$\sum_{n=0}^{\infty}\frac{1}{2^{n+1}}\sum_{k=0}^{n}\binom{n}{k}(-1)^k y^k \log\frac{k+b}{k+a} = \lim_{s\to 0}\Gamma(s)\left[\Phi(-y,s,a) - \Phi(-y,s,b)\right]$$

## AN ALTERNATIVE DERIVATION OF THE GLAISHER-KINKELIN CONSTANTS

This part was written a couple of years ago (long before I strengthened my acquaintance with the Barnes multiple gamma function).

**Theorem 4.7(a):**

$$\log A = \frac{1}{12} - \varsigma'(-1)$$

where $A$ is the Glaisher-Kinkelin constant

$$A = \lim_{n \to \infty} \left[ \sum_{k=1}^{n} k \log k - \left( \frac{n^2}{2} + \frac{n}{2} + \frac{1}{12} \right) \log n + \frac{n^2}{4} \right]$$

**Proof:**

Whilst the result is well-known, the following presentation appears to be new.

Since $\int \log(k+x) dx = (k+x)\log(k+x) - x$, upon integrating (4.4.112l) over the interval $[a,b]$ we get

$$\sum_{n=0}^{\infty} \frac{1}{2^{n+1}} \sum_{k=0}^{n} \binom{n}{k} (-1)^k (k+b) \log(k+b) - \sum_{n=0}^{\infty} \frac{1}{2^{n+1}} \sum_{k=0}^{n} \binom{n}{k} (-1)^k (k+a) \log(k+a)$$

$$-\sum_{n=0}^{\infty} \frac{1}{2^{n+1}} \sum_{k=0}^{n} \binom{n}{k} (-1)^k (b-a) = \int_a^b \log \frac{\Gamma\left(\frac{1+x}{2}\right)}{\Gamma\left(\frac{x}{2}\right)} dx + \log \sqrt{2}(b-a)$$

Since $\sum_{k=0}^{n} \binom{n}{k} (-1)^k = 0 \; \forall n \geq 1$ and is equal to 1 for $n=0$, this may be written as

(4.4.112o) $\sum_{n=0}^{\infty} \frac{1}{2^{n+1}} \sum_{k=0}^{n} \binom{n}{k} (-1)^k (k+b) \log(k+b) - \sum_{n=0}^{\infty} \frac{1}{2^{n+1}} \sum_{k=0}^{n} \binom{n}{k} (-1)^k (k+a) \log(k+a)$

$$= \int_a^b \log \frac{\Gamma\left(\frac{1+x}{2}\right)}{\Gamma\left(\frac{x}{2}\right)} dx + \frac{1}{2}(b-a) + \log \sqrt{2}(b-a)$$

Obviously we have

$$\int_a^b \log \frac{\Gamma\left(\frac{1+x}{2}\right)}{\Gamma\left(\frac{x}{2}\right)} dx = \int_a^b \log \Gamma\left(\frac{1+x}{2}\right) dx - \int_a^b \log \Gamma\left(\frac{x}{2}\right) dx$$

And using the substitution $\frac{1+x}{2} = 1 + y$ we get $\int_a^b \log \Gamma\left(\frac{1+x}{2}\right) dx = 2 \int_{(a-1)/2}^{(b-1)/2} \log \Gamma(1+y) dy$

Similarly we have $\int_a^b \log \Gamma\left(\frac{x}{2}\right) dx = 2 \int_{(a-2)/2}^{(b-2)/2} \log \Gamma(1+y) dy$.

Accordingly we have

$$2 \int_{(a-1)/2}^{(b-1)/2} \log \Gamma(1+y) dy = 2 \int_0^{(b-1)/2} \log \Gamma(1+y) dy - 2 \int_0^{(a-1)/2} \log \Gamma(1+y) dy$$

and hence we obtain

$$\int_a^b \log \frac{\Gamma\left(\frac{1+x}{2}\right)}{\Gamma\left(\frac{x}{2}\right)} dx = 2 \left\{ \int_0^{(b-1)/2} \log \Gamma(1+y) dy - \int_0^{(a-1)/2} \log \Gamma(1+y) dy \right\}$$

$$- 2 \left\{ \int_0^{(b-2)/2} \log \Gamma(1+y) dy - \int_0^{(a-2)/2} \log \Gamma(1+y) dy \right\}$$

With $a = 1$ and $b = 2$ this nicely simplifies to

$$\int_1^2 \log \frac{\Gamma\left(\frac{1+x}{2}\right)}{\Gamma\left(\frac{x}{2}\right)} dx = 2 \left\{ \int_0^{1/2} \log \Gamma(1+y) dy + \int_0^{-1/2} \log \Gamma(1+y) dy \right\}$$

$$= 2 \int_0^{1/2} \log \frac{\Gamma(1+y)}{\Gamma(1-y)} dy$$

In order to evaluate this integral, fortuitously we find that Srivastava and Choi [126, p.32] report Alexeiewsky's formula

(4.4.112p) $$\int_0^z \log\Gamma(t+1)\,dt = \frac{1}{2}[\log(2\pi)-1]z - \frac{z^2}{2} + z\log\Gamma(z+1) - \log G(z+1)$$

where $G(z)$ is the Barnes double gamma function which we shall also refer to later in (6.55).

$$G(z+1) = (2\pi)^{z/2} \exp\left[-\frac{1}{2}(z^2 + \gamma z^2 + z)\right]\prod_{n=1}^{\infty}\left(1+\frac{z}{n}\right)^n e^{-z+z^2/2n}$$

For $z = 1/2$ we get

$$\int_0^{1/2} \log\Gamma(1+y)\,dy = \frac{1}{4}[\log(2\pi)-1] - \frac{1}{8} + \frac{1}{2}\log\Gamma\left(\frac{3}{2}\right) - \log G\left(\frac{3}{2}\right)$$

From [126, p.26] we have

$$G(z+1) = \Gamma(z)G(z)$$

and hence

$$G(3/2) = \Gamma(1/2)G(1/2)$$

The value of $G(1/2)$ was originally determined by Barnes [17aa] in 1899 as

(4.4.112q) $$G(1/2) = A^{-\frac{3}{2}} \pi^{-\frac{1}{4}} e^{\frac{1}{8}} 2^{\frac{1}{24}}$$

where $A$ is the Glaisher-Kinkelin constant which we later refer to in (4.4.225) and (6.83)

$$A = \lim_{n\to\infty}\left[\sum_{k=1}^{n} k\log k - \left(\frac{n^2}{2} + \frac{n}{2} + \frac{1}{12}\right)\log n + \frac{n^2}{4}\right]$$

Hence we have

$$\log G(3/2) = \log\Gamma(1/2) + \log G(1/2)$$

$$\log G(1/2) = -\frac{3}{2}\log A - \frac{1}{4}\log\pi + \frac{1}{8} + \frac{1}{24}\log 2$$

$$\log G(3/2) = -\frac{3}{2}\log A + \frac{1}{4}\log\pi + \frac{1}{8} + \frac{1}{24}\log 2$$

Therefore, as reported in [126, p.35], we have

(4.4.112r) $$\int_0^{1/2} \log \Gamma(x+1)\,dx = -\frac{1}{2} - \frac{7}{24}\log 2 + \frac{1}{4}\log \pi + \frac{3}{2}\log A$$

For $z = -1/2$ we get [126, p.216]

(4.4.112s) $$\int_0^{-1/2} \log \Gamma(1+x)\,dx = -\frac{1}{4}[\log(2\pi) - 1] - \frac{1}{8} - \frac{1}{2}\log \Gamma\left(\frac{1}{2}\right) - \log G\left(\frac{1}{2}\right)$$

$$= -\frac{7}{24}\log 2 - \frac{1}{4}\log \pi + \frac{3}{2}\log A$$

Accordingly we have

(4.4.112t) $$\int_0^{1/2} \log \Gamma(1+x)\,dx + \int_0^{-1/2} \log \Gamma(1+x)\,dx = -\frac{1}{2} - \frac{7}{12}\log 2 + 3\log A$$

More generally we have from (4.4.112p) (and also refer to (6.65))

(4.4.112u)

$$\int_0^z \log \Gamma(1+x)\,dx + \int_0^{-z} \log \Gamma(1+x)\,dx = -z^2 + z\log \frac{\Gamma(1+z)}{\Gamma(1-z)} - \log\left[G(1+z)G(1-z)\right]$$

Hence, using (4.4.112o) and (4.4.112u) we have

(4.4.112v)

$$\sum_{n=0}^{\infty} \frac{1}{2^{n+1}} \sum_{k=0}^{n} \binom{n}{k}(-1)^k (k+2)\log(k+2) - \sum_{n=0}^{\infty} \frac{1}{2^{n+1}} \sum_{k=0}^{n} \binom{n}{k}(-1)^k (k+1)\log(k+1)$$

$$= 2\left[-\frac{1}{2} - \frac{7}{12}\log 2 + 3\log A\right] + \frac{1}{2} + \frac{1}{2}\log 2 = -\frac{1}{2} - \frac{2}{3}\log 2 + 6\log A$$

In (4.4.83) we showed that

$$\sum_{n=0}^{\infty} \frac{1}{2^{n+1}} \sum_{k=0}^{n} \binom{n}{k}\frac{(-1)^k}{(k+x)^s} = \sum_{n=0}^{\infty} \frac{(-1)^n}{(n+x)^s}$$

and for $x = 1$ we have

(4.4.112w) $$\sum_{n=0}^{\infty}\frac{1}{2^{n+1}}\sum_{k=0}^{n}\binom{n}{k}\frac{(-1)^k}{(k+1)^s} = \sum_{n=0}^{\infty}\frac{(-1)^n}{(n+1)^s} = \sum_{n=1}^{\infty}\frac{(-1)^{n+1}}{n^s} = \varsigma_a(s)$$

which is the Hasse/Sondow identity (3.11) and this is valid on the whole complex plane except for $s=1$.

We therefore have

(4.4.112x) $$\sum_{n=0}^{\infty}\frac{1}{2^{n+1}}\sum_{k=0}^{n}\binom{n}{k}\frac{(-1)^k \log(k+1)}{(k+1)^s} = -\varsigma_a'(s)$$

and hence letting $s=0$ we get

$$\sum_{n=0}^{\infty}\frac{1}{2^{n+1}}\sum_{k=0}^{n}\binom{n}{k}(-1)^k \log(k+1) = -\varsigma_a'(0)$$

We note from the following that the above coincides with (4.4.113). Since

$$\varsigma_a'(s) = (1-2^{1-s})\varsigma'(s) + 2^{1-s}\varsigma(s)\log 2$$

we have

$$\varsigma_a'(0) = -\varsigma'(0) + 2\varsigma(0)\log 2 = -\frac{1}{2}\log(2\pi) - \log 2 = \frac{1}{2}\log\frac{2}{\pi}$$

where we have used (F.2) and (F.6).

We recall from (4.3.132) in Volume II(a) that $\lim_{a\to 0}[\varsigma'(-p,a)-\varsigma'(-p)]=0$ and this may be employed in connection with the series $\sum_{n=0}^{\infty}\frac{1}{2^{n+1}}\sum_{k=0}^{n}\binom{n}{k}(-1)^k k^{p+1}\log k$.

We also see from (4.4.112x) that

(4.4.112y) $$\sum_{n=0}^{\infty}\frac{1}{2^{n+1}}\sum_{k=0}^{n}\binom{n}{k}(-1)^k (k+1)\log(k+1) = -\varsigma_a'(-1)$$

Therefore we have

(4.4.112z) $$\sum_{n=0}^{\infty}\frac{1}{2^{n+1}}\sum_{k=0}^{n}\binom{n}{k}(-1)^k k \log(k+1) = -\frac{1}{2}\log\frac{\pi}{2} - \varsigma_a'(-1)$$

We see that

$$\varsigma_a(s) = \frac{1}{1^s} - \frac{1}{2^s} + \ldots + \frac{(-1)^N}{(N-1)^s} + \left\{ \frac{(-1)^{N+1}}{N^s} + \frac{(-1)^{N+2}}{(N+1)^s} + \ldots \right\}$$

$$= \frac{1}{1^s} - \frac{1}{2^s} + \ldots + \frac{(-1)^N}{(N-1)^s} + (-1)^{N+1} \sum_{n=0}^{\infty} \frac{(-1)^n}{(n+N)^s}$$

With $N = 2$ we have

(4.4.112z(i)) $$\sum_{n=0}^{\infty} \frac{1}{2^{n+1}} \sum_{k=0}^{n} \binom{n}{k} \frac{(-1)^k}{(k+2)^s} = \sum_{n=0}^{\infty} \frac{(-1)^n}{(n+2)^s} = 1 - \varsigma_a(s)$$

and therefore upon differentiating the above and letting $s = -1$ we obtain

$$\sum_{n=0}^{\infty} \frac{1}{2^{n+1}} \sum_{k=0}^{n} \binom{n}{k} (-1)^k (k+2) \log(k+2) = \varsigma_a'(-1)$$

This is equivalent to

$$\sum_{n=0}^{\infty} \frac{1}{2^{n+1}} \sum_{k=0}^{n} \binom{n}{k} (-1)^k k \log(k+2) = \varsigma_a'(-1) - \log \frac{\pi}{2}$$

Accordingly we have

(4.4.112za)

$$\sum_{n=0}^{\infty} \frac{1}{2^{n+1}} \sum_{k=0}^{n} \binom{n}{k} (-1)^k (k+2) \log(k+2) - \sum_{n=0}^{\infty} \frac{1}{2^{n+1}} \sum_{k=0}^{n} \binom{n}{k} (-1)^k (k+1) \log(k+1) = 2\varsigma_a'(-1)$$

Since $\varsigma_a(s) = (1 - 2^{1-s}) \varsigma(s)$ and, assuming that this relation is valid for $s$ less than 0, we have

$$\varsigma_a'(s) = (1 - 2^{1-s}) \varsigma'(s) + 2^{1-s} \varsigma(s) \log 2$$

$$\varsigma_a'(-1) = -3\varsigma'(-1) + 4\varsigma(-1) \log 2$$

and since by (3.11b) $\varsigma(-1) = -\frac{1}{12}$ we have

$$\varsigma'_a(-1) = -3\varsigma'(-1) - \frac{1}{3}\log 2$$

(a more rigorous proof using Hardy's functional equation for the alternating zeta function is contained in Appendix F of Volume VI).

We have thus shown that

(4.4.112zb)

$$\sum_{n=0}^{\infty}\frac{1}{2^{n+1}}\sum_{k=0}^{n}\binom{n}{k}(-1)^k(k+2)\log(k+2) - \sum_{n=0}^{\infty}\frac{1}{2^{n+1}}\sum_{k=0}^{n}\binom{n}{k}(-1)^k(k+1)\log(k+1) = -6\varsigma'(-1) - \frac{2}{3}\log 2$$

and comparing (4.4.112v) and (4.4.112zb) we have

$$-\frac{1}{2} - \frac{2}{3}\log 2 + 6\log A = -6\varsigma'(-1) - \frac{2}{3}\log 2$$

and hence we have rediscovered the well-known relationship between $\varsigma'(-1)$ and $\log A$, namely

(4.4.112zc) $$\log A = \frac{1}{12} - \varsigma'(-1)$$

This relationship was found by Voros [133b] in 1987 and was extended by Adamchik's formula [4] for the generalised Glaisher-Kinkelin constants which was derived in his 1998 paper "PolyGamma functions of negative order"

$$\log A_n = \frac{B_{n+1}H_n}{n+1} - \varsigma'(-n) \quad \text{(where } A_1 = A\text{)}$$

We have

$$\sum_{n=0}^{\infty}\frac{1}{2^{n+1}}\sum_{k=0}^{n}\binom{n}{k}(-1)^k(k+2)\log(k+2)$$

$$= \sum_{n=0}^{\infty}\frac{1}{2^{n+1}}\sum_{k=0}^{n}\binom{n}{k}(-1)^k k \log(k+2) + 2\sum_{n=0}^{\infty}\frac{1}{2^{n+1}}\sum_{k=0}^{n}\binom{n}{k}(-1)^k \log(k+2)$$

$$= \sum_{n=0}^{\infty}\frac{1}{2^{n+1}}\sum_{k=0}^{n}\binom{n}{k}(-1)^k(k+2)\log(k+2) + \log\frac{\pi}{2}$$

and

$$\sum_{n=0}^{\infty}\frac{1}{2^{n+1}}\sum_{k=0}^{n}\binom{n}{k}(-1)^k(k+1)\log(k+1)$$

$$= \sum_{n=0}^{\infty} \frac{1}{2^{n+1}} \sum_{k=0}^{n} \binom{n}{k} (-1)^k k \log(k+1) + \sum_{n=0}^{\infty} \frac{1}{2^{n+1}} \sum_{k=0}^{n} \binom{n}{k} (-1)^k \log(k+1)$$

$$= \sum_{n=0}^{\infty} \frac{1}{2^{n+1}} \sum_{k=0}^{n} \binom{n}{k} (-1)^k (k+1) \log(k+2) - \frac{1}{2} \log \frac{\pi}{2}$$

We therefore get

(4.4.112zd)

$$\sum_{n=0}^{\infty} \frac{1}{2^{n+1}} \sum_{k=0}^{n} \binom{n}{k} (-1)^k (k+2) \log(k+2) - \sum_{n=0}^{\infty} \frac{1}{2^{n+1}} \sum_{k=0}^{n} \binom{n}{k} (-1)^k (k+1) \log(k+1)$$

$$= \sum_{n=0}^{\infty} \frac{1}{2^{n+1}} \sum_{k=0}^{n} \binom{n}{k} (-1)^k k \log \frac{(k+2)}{(k+1)} + \frac{3}{2} \log \frac{\pi}{2} = -6\varsigma'(-1) - \frac{2}{3} \log 2$$

Applying Euler's transformation of series (4.4.112i) we get

(4.4.112zdi)

$$\sum_{n=1}^{\infty} (-1)^{n+1} (n-1) \log \frac{n}{n-1} = -6\varsigma'(-1) - \frac{2}{3} \log 2 - \frac{3}{2} \log \frac{\pi}{2}$$

We have
$$\sum_{n=1}^{\infty} (-1)^{n+1} (n-1) \log \frac{n}{n-1} = \sum_{n=2}^{\infty} (-1)^{n+1} (n-1) \log \frac{n}{n-1} = \sum_{n=2}^{\infty} (-1)^{n+1} n \log \frac{n}{n-1} + \sum_{n=2}^{\infty} (-1)^{n+1} \log \frac{n-1}{n}$$

As in Sondow's paper [123aa] we have using the Wallis formula

$$\sum_{n=2}^{\infty} (-1)^{n+1} \log \frac{n-1}{n} = \log \prod_{n=2}^{\infty} \left( \frac{n-1}{n} \right)^{(-1)^{n+1}} = \log \left[ \frac{2}{1} \frac{2}{3} \frac{4}{3} \frac{4}{5} \frac{6}{5} \frac{6}{7} \cdots \right] = \log \frac{\pi}{2}$$

and hence
$$\sum_{n=1}^{\infty} (-1)^{n+1} (n-1) \log \frac{n}{n-1} = \sum_{n=2}^{\infty} (-1)^{n+1} n \log \frac{n}{n-1} + \log \frac{\pi}{2}$$

This then results in

$$\sum_{n=2}^{\infty}(-1)^{n+1}n\log\frac{n}{n-1}=-6\varsigma'(-1)+\frac{4}{3}\log 2-2\log\pi$$

In Sondow's paper [121] there is a reference to a formula derived by Hardy in 1922 which is valid for $\mathrm{Re}(s) > -1$ and $s \neq 1$

$$2\varsigma(s)\left(1-2^{1-s}\right)-1=\sum_{n=1}^{\infty}(-1)^{n+1}\left[\frac{1}{n^s}-\frac{1}{(n+1)^s}\right]$$

and upon differentiation we obtain

$$2\varsigma(s)2^{1-s}\log 2+2\varsigma'(s)\left(1-2^{1-s}\right)=\sum_{n=1}^{\infty}(-1)^{n+1}\left[\frac{\log(n+1)}{(n+1)^s}-\frac{\log n}{n^s}\right]$$

and, with $s=-1$ (assuming that it was valid), the above expression bears a close resemblance to (4.4.112zdi).

$$-6\varsigma'(-1)-\frac{2}{3}\log 2=\sum_{n=1}^{\infty}(-1)^{n+1}\left[(n+1)\log(n+1)-n\log n\right]$$

As an aside, Yung and Sondow have shown how Hardy's formula can be used to derive the Wallis identity (see the Wallis formula on the Mathworld website).

**Theorem 4.7(b):**

$$\log B = -\varsigma'(-2)$$

**Proof:**

We now wish to take the above procedure one step further. From (4.4.112o) with $a=1$ we have

$$\sum_{n=0}^{\infty}\frac{1}{2^{n+1}}\sum_{k=0}^{n}\binom{n}{k}(-1)^k(k+b)\log(k+b)-\sum_{n=0}^{\infty}\frac{1}{2^{n+1}}\sum_{k=0}^{n}\binom{n}{k}(-1)^k(k+1)\log(k+1)$$

$$=\int_{1}^{b}\log\frac{\Gamma\left(\frac{1+x}{2}\right)}{\Gamma\left(\frac{x}{2}\right)}dx+\frac{1}{2}(b-1)+(b-1)\log\sqrt{2}$$

We now wish to integrate the above equation with respect to $b$ over the interval $[1,2]$. Since

$$\int (k+x)\log(k+x)dx = \frac{1}{2}(k+x)^2 \log(k+x) - \frac{1}{2}(k+x)^2$$

for the left-hand side we get

$$\sum_{n=0}^{\infty}\frac{1}{2^{n+1}}\sum_{k=0}^{n}\binom{n}{k}(-1)^k \int_{1}^{2}(k+b)\log(k+b)\,db - \sum_{n=0}^{\infty}\frac{1}{2^{n+1}}\sum_{k=0}^{n}\binom{n}{k}(-1)^k(k+1)\log(k+1) =$$

$$\frac{1}{2}\sum_{n=0}^{\infty}\frac{1}{2^{n+1}}\sum_{k=0}^{n}\binom{n}{k}(-1)^k(k+2)^2\log(k+2) - \frac{1}{2}\sum_{n=0}^{\infty}\frac{1}{2^{n+1}}\sum_{k=0}^{n}\binom{n}{k}(-1)^k(k+1)^2\log(k+1)$$

$$-\frac{1}{2}\sum_{n=0}^{\infty}\frac{1}{2^{n+1}}\sum_{k=0}^{n}\binom{n}{k}(-1)^k\left[(k+2)^2 - (k+1)^2\right] - \sum_{n=0}^{\infty}\frac{1}{2^{n+1}}\sum_{k=0}^{n}\binom{n}{k}(-1)^k(k+1)\log(k+1)$$

Using (4.4.112x) this is equal to

$$\frac{1}{2}\sum_{n=0}^{\infty}\frac{1}{2^{n+1}}\sum_{k=0}^{n}\binom{n}{k}(-1)^k(k+2)^2\log(k+2) - \frac{1}{2}\sum_{n=0}^{\infty}\frac{1}{2^{n+1}}\sum_{k=0}^{n}\binom{n}{k}(-1)^k(k+1)^2\log(k+1) - \frac{1}{2} + \varsigma'_a(-1)$$

since

$$\sum_{n=0}^{\infty}\frac{1}{2^{n+1}}\sum_{k=0}^{n}\binom{n}{k}(-1)^k\left[(k+2)^2 - (k+1)^2\right] = 2\sum_{n=0}^{\infty}\frac{1}{2^{n+1}}\sum_{k=0}^{n}\binom{n}{k}(-1)^k k + 3\sum_{n=0}^{\infty}\frac{1}{2^{n+1}}\sum_{k=0}^{n}\binom{n}{k}(-1)^k = 1$$

as can be easily determined as follows. In deriving (3.11b) we showed that

$$\sum_{n=0}^{\infty}\frac{1}{2^{n+1}}\sum_{k=0}^{n}\binom{n}{k}(-1)^k(k+1) = \sum_{n=0}^{\infty}\frac{1}{2^{n+1}}\sum_{k=0}^{n}\binom{n}{k}(-1)^k k + \sum_{n=0}^{\infty}\frac{1}{2^{n+1}}\sum_{k=0}^{n}\binom{n}{k}(-1)^k = \frac{1}{4}$$

$$\sum_{n=0}^{\infty}\frac{1}{2^{n+1}}\sum_{k=0}^{n}\binom{n}{k}(-1)^k = \frac{1}{2}$$

and hence we have

$$\sum_{n=0}^{\infty}\frac{1}{2^{n+1}}\sum_{k=0}^{n}\binom{n}{k}(-1)^k k = -\frac{1}{4}$$

Differentiation of (4.4.83) gives us

$$\sum_{n=0}^{\infty}\frac{1}{2^{n+1}}\sum_{k=0}^{n}\binom{n}{k}(-1)^{k}(k+1)^{2}\log(k+1)=-\varsigma_{a}'(-2)$$

and differentiation of (4.4.112z(i)) results in

$$\sum_{n=0}^{\infty}\frac{1}{2^{n+1}}\sum_{k=0}^{n}\binom{n}{k}(-1)^{k}(k+2)^{2}\log(k+2)=\varsigma_{a}'(-2)$$

Accordingly, the left-hand side becomes equal to $\varsigma_{a}'(-2)-\frac{1}{2}+\varsigma_{a}'(-1)$.

In passing we note from (4.4.112z (i)) and (4.4.84) that

(4.4.112ze) $$\sum_{n=0}^{\infty}\frac{1}{2^{n+1}}\sum_{k=0}^{n}\binom{n}{k}\frac{(-1)^{k}}{(k+2)^{s}}=1-\varsigma_{a}(s)=1-\sum_{n=0}^{\infty}\frac{1}{2^{n+1}}\sum_{k=0}^{n}\binom{n}{k}\frac{(-1)^{k}}{(k+1)^{s}}$$

and therefore with differentiation we see that

(4.4.112zf) $$\sum_{n=0}^{\infty}\frac{1}{2^{n+1}}\sum_{k=0}^{n}\binom{n}{k}\frac{(-1)^{k}}{(k+2)^{s}}\log(k+2)=\varsigma_{a}'(s)=-\sum_{n=0}^{\infty}\frac{1}{2^{n+1}}\sum_{k=0}^{n}\binom{n}{k}\frac{(-1)^{k}}{(k+1)^{s}}\log(k+1)$$

We also have for part of the right-hand side

$$f(b)=\int_{1}^{b}\log\frac{\Gamma\left(\frac{1+x}{2}\right)}{\Gamma\left(\frac{x}{2}\right)}dx=2\int_{0}^{(b-1)/2}\log\Gamma(1+y)\,dy-2\int_{0}^{(b-2)/2}\log\Gamma(1+y)\,dy+2\int_{0}^{-1/2}\log\Gamma(1+y)\,dy$$

We have from (4.4.112p)

$$\int_{0}^{(b-1)/2}\log\Gamma(1+y)\,dy=\frac{1}{2}[\log(2\pi)-1]\left(\frac{b-1}{2}\right)-\frac{1}{2}\left(\frac{b-1}{2}\right)^{2}+\left(\frac{b-1}{2}\right)\log\Gamma\left(\frac{b+1}{2}\right)-\log G\left(\frac{b+1}{2}\right)$$

and

$$\int_{0}^{(b-2)/2}\log\Gamma(1+y)\,dy=\frac{1}{2}[\log(2\pi)-1]\left(\frac{b-2}{2}\right)-\frac{1}{2}\left(\frac{b-2}{2}\right)^{2}+\left(\frac{b-2}{2}\right)\log\Gamma\left(\frac{b}{2}\right)-\log G\left(\frac{b}{2}\right)$$

It is easily shown that

$$\int_1^2 \left[\frac{1}{2}[\log(2\pi)-1]\left(\frac{b-1}{2}\right)-\frac{1}{2}\left(\frac{b-1}{2}\right)^2\right]-\left[\frac{1}{2}[\log(2\pi)-1]\left(\frac{b-2}{2}\right)-\frac{1}{2}\left(\frac{b-2}{2}\right)^2\right]db = \frac{1}{4}[\log(2\pi)-1]$$

We now need to determine the following integral

$$\int_1^2 \log G\left(\frac{b+1}{2}\right)db = 2\int_0^{1/2} \log G(1+x)dx$$

and in [126, p.39] we are kindly provided with the answer

(4.4.112zg) $$\int_0^{1/2} \log G(1+x)dx = \frac{1}{24}+\frac{1}{24}\log 2+\frac{1}{16}\log\pi-\frac{1}{4}\log A-\frac{7}{4}\log B$$

Similarly we obtain

$$\int_1^2 \log G\left(\frac{b}{2}\right)db = -2\int_0^{-1/2} \log G(1+x)dx$$

and we have from [126, p.216]

(4.4.112zh) $$\int_0^{-1/2} \log G(1+x)dx = -\frac{1}{24}+\frac{1}{12}\log 2+\frac{1}{16}\log\pi+\frac{1}{4}\log A-\frac{7}{4}\log B$$

(where, as we shall independently show at the end of this proof, $\log B = -\varsigma'(-2) = \frac{\varsigma(3)}{4\pi^2}$).
The constant $B = A_2$ was first considered by Choi and Srivastava [45aa] in 1999.

Hence we obtain

(4.4.112zi) $$\int_0^{1/2} \log G(1+t)dt + \int_0^{-1/2} \log G(1+t)dt = \frac{1}{8}\log 2+\frac{1}{8}\log\pi-\frac{7}{2}\log B$$

Therefore we have

$$\int_1^2 \log G\left(\frac{b}{2}\right)db - \int_1^2 \log G\left(\frac{b+1}{2}\right)db = -2\left[\int_0^{1/2} \log G(1+x)dx + \int_0^{-1/2} \log G(1+x)dx\right]$$

$$= -\frac{1}{4}\log 2 - \frac{1}{4}\log \pi + 7\log B$$

$$= -\frac{1}{4}\log(2\pi) + 7\log B$$

We note in passing that by substituting $x \to -x$ in the second integral that

(4.4.112zj) $\quad \int_0^{1/2} \log G(1+x)\,dx + \int_0^{-1/2} \log G(1+x)\,dx = \int_0^{1/2} \log \frac{G(1+x)}{G(1-x)}\,dx$

The last integral then reminds us of (6.69b)

$$\log \frac{G(1+x)}{G(1-x)} = x\log(2\pi) - \int_0^x \pi t \cot \pi t \, dt$$

and hence we have

$$\int_0^{1/2} \log \frac{G(1+x)}{G(1-x)}\,dx = \frac{1}{2}\log(2\pi) - \int_0^{1/2} dx \int_0^x \pi t \cot \pi t \, dt$$

The following integrals also need to be determined

$$\int_1^2 \left(\frac{b-1}{2}\right)\log\Gamma\left(\frac{b+1}{2}\right)db = 2\int_0^{1/2} x\log\Gamma(1+x)\,dx$$

$$\int_1^2 \left(\frac{b-2}{2}\right)\log\Gamma\left(\frac{b}{2}\right)db = -2\int_0^{-1/2} x\log\Gamma(1+x)\,dx$$

Therefore we get

$$\int_1^2 \left(\frac{b-1}{2}\right)\log\Gamma\left(\frac{b+1}{2}\right)db - \int_1^2 \left(\frac{b-2}{2}\right)\log\Gamma\left(\frac{b}{2}\right)db = 2\left[\int_0^{1/2} x\log\Gamma(1+x)\,dx + \int_0^{-1/2} x\log\Gamma(1+x)\,dx\right]$$

We have from (4.4.112p)

(4.4.112zk) $\quad \int_0^z \log\Gamma(1+t)\,dt = \frac{1}{2}[\log(2\pi)-1]z - \frac{z^2}{2} + z\log\Gamma(1+z) - \log G(1+z)$

and integration by parts readily gives us

$$\int_0^x t\log\Gamma(1+t)\,dt = t\left[\frac{1}{2}[\log(2\pi)-1]t - \frac{t^2}{2} + t\log\Gamma(1+t) - \log G(1+t)\right]\Big|_0^x$$

$$- \int_0^x \left[\frac{1}{2}[\log(2\pi)-1]t - \frac{t^2}{2} + t\log\Gamma(1+t) - \log G(1+t)\right]dt$$

Hence we obtain

(4.4.112zl)

$$2\int_0^x t\log\Gamma(1+t)\,dt = \frac{1}{4}[\log(2\pi)-1]t^2 - \frac{1}{3}t^3 + t^2\log\Gamma(1+t) - t\log G(1+t) + \int_0^x \log G(1+t)\,dt$$

and therefore

$$2\left[\int_0^{1/2} t\log\Gamma(1+t)\,dt + \int_0^{-1/2} t\log\Gamma(1+t)\,dt\right]$$

$$= \frac{1}{8}[\log(2\pi)-1] + \frac{1}{4}\log\Gamma(3/2)\Gamma(1/2) - \frac{1}{2}\log\frac{G(3/2)}{G(1/2)} + \frac{1}{8}\log 2 + \frac{1}{8}\log\pi - \frac{7}{2}\log B$$

$$= -\frac{1}{8} + \frac{1}{4}\log\pi - \frac{7}{2}\log B$$

The left and right-hand sides are therefore equal to

$$\varsigma_a'(-2) - \frac{1}{2} + \varsigma_a'(-1) =$$

$$\frac{1}{2}[\log(2\pi)-1] + 2\left\{-\frac{1}{8} + \frac{1}{4}\log\pi - \frac{7}{2}\log B\right\} + 2\left\{-\frac{1}{4}\log(2\pi) + 7\log B\right\}$$

$$+ 2\left\{-\frac{7}{24}\log 2 - \frac{1}{4}\log\pi + \frac{3}{2}\log A\right\} + \frac{1}{4}(1+\log 2)$$

$$= -\frac{1}{2} + 7\log B + 3\log A - \frac{7}{12}\log 2$$

We showed previously at (4.4.112za) that

$$\varsigma'_a(-1) = -3\varsigma'(-1) - \frac{1}{3}\log 2$$

and similarly we have

$$\varsigma'_a(-2) = -7\varsigma'(-2) + 4\varsigma(-2)\log 2$$

This then becomes

$$\varsigma'_a(-2) = -7\varsigma'(-2)$$

since $\varsigma(-2n) = 0$ per (3.11c)

Using (4.4.112zc), we then have an alternative proof of the well-known formula

(4.4.112zm) $\qquad -\varsigma'(-2) = \log B$

With considerable labour, it may be possible to extend the above analysis to the other Glaisher-Kinkelin constants. Some further applications of the Barnes function are given in Volume II(a) and Volume V.

We also note from [75aa] that

(4.4.112zn) $\qquad 2\log A - \dfrac{1}{4} = \sum_{n=0}^{\infty} \dfrac{1}{n+1} \sum_{k=0}^{n} \binom{n}{k} (-1)^{k+1}(k+1)^2 \log(k+1)$

The following analysis is partly based on [6c].

Lemma: If $a$ and $b$ are positive real numbers, then

(4.4.112zo) $\qquad \sum_{n=1}^{\infty} \left[ \dfrac{1}{an} - \dfrac{1}{an+b} \right] = \dfrac{1}{b} + \dfrac{\gamma + \psi(b/a)}{a}$

Proof: We have

$$\sum_{n=1}^{\infty} \left[ \dfrac{1}{an} - \dfrac{1}{an+b} \right] = \lim_{s \to 1+} \sum_{n=1}^{\infty} \left[ \dfrac{1}{(an)^s} - \dfrac{1}{(an+b)^s} \right] = \lim_{s \to 1+} \dfrac{1}{a^s} \sum_{n=1}^{\infty} \left[ \dfrac{1}{n^s} - \dfrac{1}{(n+b/a)^s} \right]$$

$$= \lim_{s \to 1+} \frac{1}{a} \sum_{n=1}^{\infty} \left[ \frac{1}{n^s} - \frac{1}{(n+b/a)^s} \right]$$

$$= \frac{1}{a} \lim_{s \to 1+} \left[ \varsigma(s) - \varsigma\left(s, \frac{b}{a}\right) + \left(\frac{a}{b}\right)^s \right]$$

$$= \frac{1}{b} + \frac{1}{a} \lim_{s \to 1+} \left[ \left\{ \varsigma(s) - \frac{1}{s-1} \right\} - \left\{ \varsigma\left(s, \frac{b}{a}\right) - \frac{1}{s-1} \right\} \right]$$

$$= \frac{1}{b} + \frac{1}{a} \left[ \gamma + \frac{\Gamma'\left(\frac{b}{a}\right)}{\Gamma\left(\frac{b}{a}\right)} \right] = \frac{1}{b} + \frac{\gamma + \psi(b/a)}{a}$$

From (E.12a) we have the well-known formula for $x > 0$

$$\sum_{n=1}^{\infty} \left[ \frac{x}{n} - \log\left(1 + \frac{x}{n}\right) \right] = \log \Gamma(x) + \gamma x + \log x$$

and letting $x = 1/a$ we obtain

$$\sum_{n=1}^{\infty} \left[ \frac{1}{an} - \log \frac{an+1}{an} \right] = \log \Gamma\left(\frac{1}{a}\right) + \frac{\gamma}{a} - \log a$$

Lemma: If $a$ and $b$ are positive real numbers, then

(4.4.112zp) $$\sum_{n=1}^{\infty} \left[ \frac{1}{an+b} - \log \frac{an+b+1}{an+b} \right] = \log \Gamma\left(\frac{b+1}{a}\right) - \log \Gamma\left(\frac{b}{a}\right) - \frac{\psi(b/a)}{a}$$

Proof: We write the $n$th term of the series as

$$\frac{1}{an+b} - \log \frac{an+b+1}{an+b} = \frac{1}{an+b} - \frac{1}{an} - \frac{b}{an} + \log\left(1 + \frac{b}{an}\right) - \frac{b+1}{an} - \log\left(1 + \frac{b+1}{an}\right)$$

and the result follows directly from the above.

Letting $a = 1$ and $b = x$ we have

$$\sum_{n=1}^{\infty} \left[ \frac{1}{n+x} - \log \frac{n+x+1}{n+x} \right] = \log \Gamma(x+1) - \log \Gamma(x) - \psi(x)$$

and integration results in

$$\int_1^y \sum_{n=1}^{\infty}\left[\frac{1}{n+x} - \log\frac{n+x+1}{n+x}\right]dx = \int_1^y \left[\log\Gamma(x+1) - \log\Gamma(x)\right]dx - \log\Gamma(y)$$

Since

$$\int_1^y \left[\log\Gamma(x+1) - \log\Gamma(x)\right]dx = \int_1^y \log x\, dx = y\log y - y + 1$$

we do not have to use Alexeiewsky's formula (4.4.112p)

We then obtain

$$\sum_{n=1}^{\infty}\left[(n+y+1)\log\frac{n+y+1}{n+y} - (n+2)\log\frac{n+2}{n+1}\right] = y\log y - y + 1 - \log\Gamma(y)$$

We now refer back to (4.4.112l)

$$\sum_{n=0}^{\infty}\frac{1}{2^{n+1}}\sum_{k=0}^{n}\binom{n}{k}(-1)^k \log(k+x) = \log\frac{\Gamma\left(\frac{1+x}{2}\right)}{\Gamma\left(\frac{x}{2}\right)} + \log\sqrt{2}$$

and note that Guillera and Sondow [75aa] have recently shown that

(4.4.112zq) $$-\frac{\partial\Phi}{\partial s}(-1, s, x)\bigg|_{s=0} = \log\frac{\Gamma\left(\frac{1+x}{2}\right)}{\Gamma\left(\frac{x}{2}\right)} + \log\sqrt{2}$$

We refer to the following formula in (6.69c)

$$Cl_2(2\pi t) = \sum_{n=1}^{\infty}\frac{\sin 2\pi nt}{n^2} = -2\pi t\log\left[\frac{\sin \pi t}{\pi}\right] - 2\pi \log\frac{G(1+t)}{G(1-t)}$$

and integration produces

$$\frac{1}{4\pi^2}\sum_{n=1}^{\infty}\frac{\cos(2\pi nt) - 1}{n^2} = \int_0^t t\log[\sin \pi t]dt - \frac{t^2}{2}\log \pi + \int_0^t \log G(1+t)dt + \int_0^{-t} \log G(1+t)dt$$

Letting $t = 1/2$ we get

$$\int\limits_0^{1/2} \log G(1+t)dt + \int\limits_0^{-1/2} \log G(1+t)dt = \frac{1}{8}\log 2 + \frac{1}{8}\log \pi - \frac{7}{8}\frac{\varsigma(3)}{\pi^2}$$

and this is in agreement with (4.4.112zi) where we showed that

$$\int\limits_0^{1/2} \log G(1+t)dt + \int\limits_0^{-1/2} \log G(1+t)dt = \frac{1}{8}\log 2 + \frac{1}{8}\log \pi - \frac{7}{2}\log B$$

The above is a particular case of the more general formula shown by Choi, Srivastava and Adamchik in [45].

(4.4.112zs)
$$\int\limits_0^u \log G(1+t)dt + \int\limits_0^{-u} \log G(1+t)dt = \frac{1}{4\pi^2}\sum_{n=1}^{\infty}\frac{\cos 2n\pi u - 1}{n^3} - \int\limits_0^u t\log\sin \pi t\, dt + \frac{u^2}{2}\log \pi$$

We now recall the Hasse/Sondow formula (3.11) again

$$\varsigma_a(s) = \sum_{n=0}^{\infty}\frac{1}{2^{n+1}}\sum_{k=0}^{n}\binom{n}{k}\frac{(-1)^k}{(k+1)^s}$$

and immediately deduce from (4.4.112m) that with $x=1$

(4.4.112zt) $\qquad \varsigma_a(s) = \frac{(-1)^{s+1}}{2^s(s-1)!}\left[\psi^{(s-1)}\left(\frac{1+x}{2}\right) - \psi^{(s-1)}\left(\frac{x}{2}\right)\right]_{x=1}$

Letting $x \to x+1$ in (4.4.112m) we get

(4.4.112zu) $\qquad \sum_{n=0}^{\infty}\frac{1}{2^{n+1}}\sum_{k=0}^{n}\binom{n}{k}\frac{(-1)^k}{(k+x+1)} = \frac{1}{2}\left[\psi\left(\frac{2+x}{2}\right) - \psi\left(\frac{1+x}{2}\right)\right]$

Employing the binomial theorem we deduce

$$\sum_{n=0}^{\infty}\frac{1}{2^{n+1}}\sum_{k=0}^{n}\binom{n}{k}\frac{(-1)^k}{(k+x+1)} = \sum_{n=0}^{\infty}\frac{1}{2^{n+1}}\sum_{k=0}^{n}\binom{n}{k}\frac{(-1)^k}{1+k}\sum_{m=0}^{\infty}\frac{(-1)^k x^m}{(1+k)^m}$$

The coefficient of $x^m$ is easily seen to be

$$(-1)^m \sum_{n=0}^{\infty} \frac{1}{2^{n+1}} \sum_{k=0}^{n} \binom{n}{k} \frac{(-1)^k}{(1+k)^{m+1}} = (-1)^m \varsigma_a(m+1)$$

and we therefore have

(4.4.112zv) $$\sum_{n=0}^{\infty} \frac{1}{2^{n+1}} \sum_{k=0}^{n} \binom{n}{k} \frac{(-1)^k}{(k+x+1)} = \sum_{n=0}^{\infty} (-1)^n \varsigma_a(n+1) x^n$$

With $x = 0$ we get the known result

$$\sum_{n=0}^{\infty} \frac{1}{2^{n+1}} \sum_{k=0}^{n} \binom{n}{k} \frac{(-1)^k}{k+1} = \log 2$$

Alternatively we see that

(4.4.112zw) $$\frac{1}{2}\left[\psi\left(\frac{2+x}{2}\right) - \psi\left(\frac{x+1}{2}\right)\right] = \sum_{n=0}^{\infty} (-1)^n \varsigma_a(n+1) x^n$$

In 1961 Beumer [23a] showed that G&R [74, p.897]

(4.4.112zy) $$\frac{1}{2}\left[\psi\left(\frac{x+2}{2}\right) - \psi\left(\frac{x+1}{2}\right)\right] = \log 2 - \sum_{n=1}^{\infty} (-1)^{n+1}\left(1 - \frac{1}{2^n}\right) \varsigma(n+1) x^n$$

and it is easily seen that the above two formulae are equivalent.

Integrating (4.4.112zw) we get

$$\frac{1}{2} \int_0^t \left[\psi\left(\frac{x+1}{2}\right) - \psi\left(\frac{2+x}{2}\right)\right] dx = \sum_{n=0}^{\infty} (-1)^n \frac{\varsigma_a(n+1)}{n+1} t^{n+1}$$

We have seen in (4.4.112fa) that

$$\frac{1}{2} \int_0^1 \left[\psi\left(\frac{x+1}{2}\right) - \psi\left(\frac{2+x}{2}\right)\right] dx = \log \frac{\pi}{2}$$

and we therefore have

(4.4.112zz) $$\log \frac{\pi}{2} = \sum_{n=0}^{\infty} (-1)^n \frac{\varsigma_a(n+1)}{n+1}$$

**(xi) Theorem 4.8:**

(4.4.113)
$$\log\frac{\pi}{2} = \sum_{n=1}^{\infty}\frac{1}{2^n}\sum_{k=0}^{n}\binom{n}{k}(-1)^{k+1}\log(1+k)$$

$$\int_0^1\left\{u+\frac{\log[1-u(1-y)]}{(1-y)}\right\}\frac{dy}{\log y} = \sum_{n=1}^{\infty}\frac{u^{n+1}}{n+1}\sum_{k=0}^{n}\binom{n}{k}(-1)^{k+1}\log(1+k)$$

**Proof:**

Using Anglesio's formula (4.4.94) we have

(4.4.114)
$$\sum_{n=1}^{\infty}\frac{1}{2^n}\int_0^{\infty}\frac{e^{-x}(1-e^{-x})^n}{x}\,dx = \sum_{n=1}^{\infty}\frac{1}{2^n}\sum_{k=0}^{n}\binom{n}{k}(-1)^{k+1}\log(1+k)$$

and, completing the summation of the left-hand side, we have

(4.4.115)
$$\sum_{n=1}^{\infty}\frac{1}{2^n}\int_0^{\infty}\frac{e^{-x}(1-e^{-x})^n}{x}\,dx = \int_0^{\infty}\frac{e^{-x}(1-e^{-x})}{x(1+e^{-x})}\,dx$$

The substitution $y = e^{-x}$ gives us

(4.4.116)
$$\int_0^{\infty}\frac{e^{-x}(1-e^{-x})}{x(1+e^{-x})}\,dx = -\int_0^1\frac{(1-y)}{(1+y)\log y}\,dy$$

$$= \log(\pi/2)$$

where, in the final step, we used the result (4.4.112) from the preceding theorem. Accordingly we have

$$\log\frac{\pi}{2} = \sum_{n=1}^{\infty}\frac{1}{2^n}\sum_{k=0}^{n}\binom{n}{k}(-1)^{k+1}\log(1+k)$$

We may generalise the above result as follows. Using (4.4.94) we have

$$\sum_{n=1}^{\infty}t^n\int_0^{\infty}\frac{e^{-x}(1-e^{-x})^n}{x}\,dx = \sum_{n=1}^{\infty}t^n\sum_{k=0}^{n}\binom{n}{k}(-1)^{k+1}\log(1+k)$$

and, completing the summation of the left-hand side, we have

$$\sum_{n=1}^{\infty} t^n \int_0^{\infty} \frac{e^{-x}(1-e^{-x})^n}{x} dx = \int_0^{\infty} \frac{e^{-x} t(1-e^{-x})}{x\left[1-t(1+e^{-x})\right]} dx$$

The substitution $y = e^{-x}$ gives us

$$\int_0^{\infty} \frac{e^{-x} t(1-e^{-x})}{x\left[1-t(1+e^{-x})\right]} dx = -\int_0^1 \frac{t(1-y)}{[1-t(1-y)]\log y} dy$$

and accordingly we have

(4.4.116a) $\quad -\int_0^1 \frac{t(1-y)}{[1-t(1-y)]\log y} dy = \sum_{n=1}^{\infty} t^n \sum_{k=0}^{n} \binom{n}{k} (-1)^{k+1} \log(1+k)$

Integrating the above with respect to $t$ results in

$$-\int_0^u dt \int_0^1 \frac{t(1-y)}{[1-t(1-y)]\log y} dy = -\int_0^1 \frac{dy}{\log y} \int_0^u \frac{t(1-y)}{[1-t(1-y)]} dt$$

We have

$$-\int_0^u \frac{t(1-y)}{[1-t(1-y)]} dt = t + \frac{\log[1-t(1-y)]}{(1-y)}\bigg|_0^u = u + \frac{\log[1-u(1-y)]}{(1-y)}$$

and therefore we get

(4.4.116b) $\quad \int_0^1 \left\{u + \frac{\log[1-u(1-y)]}{(1-y)}\right\} \frac{dy}{\log y} = \sum_{n=1}^{\infty} \frac{u^{n+1}}{n+1} \sum_{k=0}^{n} \binom{n}{k} (-1)^{k+1} \log(1+k)$

Letting $u = 1$ we obtain

$$\int_0^1 \left\{1 + \frac{\log y}{(1-y)}\right\} \frac{dy}{\log y} = \sum_{n=1}^{\infty} \frac{1}{n+1} \sum_{k=0}^{n} \binom{n}{k} (-1)^{k+1} \log(1+k)$$

and referring to (4.4.92a) we see that the integral is equal to $\gamma$.

We note that

$$\int \frac{\log(1-ax)}{x} dx = -Li_2(ax)$$

and dividing (4.4.116b) by $u$ and integrating results in another derivation of (4.4.100zi)

(4.4.116c) $\int_0^1 \left\{ u - \frac{Li_2[u(1-y)]}{(1-y)} \right\} \frac{dy}{\log y} = \sum_{n=1}^{\infty} \frac{u^{n+1}}{(n+1)^2} \sum_{k=0}^{n} \binom{n}{k} (-1)^{k+1} \log(1+k)$

Sondow [123a] and Coffey [45c] have recently obtained alternative proofs of (4.4.113) using the Hasse/Sondow formula (3.11)

$$\varsigma_a(s) = \sum_{n=0}^{\infty} \frac{1}{2^{n+1}} \sum_{k=0}^{n} \binom{n}{k} \frac{(-1)^k}{(k+1)^s}$$

Differentiating with respect to $s$ we obtain

$$\varsigma_a'(s) = -\sum_{n=0}^{\infty} \frac{1}{2^{n+1}} \sum_{k=0}^{n} \binom{n}{k} \frac{(-1)^k \log(k+1)}{(k+1)^s}$$

and we have $\varsigma_a'(s) = (1-2^{1-s})\varsigma'(s) + 2^{1-s}\varsigma(s)\log 2$

Therefore we have

$$\varsigma'(s) = -\frac{\log 2}{(2^{s-1}-1)}\varsigma(s) - \frac{1}{(1-2^{1-s})} \sum_{n=0}^{\infty} \frac{1}{2^{n+1}} \sum_{k=0}^{n} \binom{n}{k} \frac{(-1)^k \log(k+1)}{(k+1)^s}$$

or

$$\frac{\varsigma'(s)}{\varsigma(s)} = -\frac{\log 2}{(2^{s-1}-1)} - \frac{1}{\varsigma(s)(1-2^{1-s})} \sum_{n=0}^{\infty} \frac{1}{2^{n+1}} \sum_{k=0}^{n} \binom{n}{k} \frac{(-1)^k \log(k+1)}{(k+1)^s}$$

Letting $s = 0$, and noting from (3.11a) that $\varsigma(0) = -\frac{1}{2}$ and from (F.5c) that

$$\frac{\varsigma'(0)}{\varsigma(0)} = \log(2\pi)$$

we get as before

$$\log \frac{\pi}{2} = \sum_{n=0}^{\infty} \frac{1}{2^n} \sum_{k=0}^{n} \binom{n}{k} (-1)^{k+1} \log(k+1)$$

We now recall the identity involving the alternating Hurwitz-Lerch zeta function which we proved at (4.4.24a)

$$\sum_{n=0}^{\infty}\frac{1}{2^{n+1}}\sum_{k=0}^{n}\binom{n}{k}\frac{(-1)^k}{(k+x)^s}=\sum_{k=0}^{\infty}\frac{(-1)^k}{(k+x)^s}$$

and the well-known result which was derived at (C.61) in Appendix C

$$\sum_{k=1}^{\infty}\frac{(-1)^k \log k}{k}=\log 2\left[\gamma-\frac{\log 2}{2}\right]$$

We therefore have

(4.4.116d) $$\varsigma'_a(1)=\sum_{k=1}^{\infty}\frac{(-1)^k \log k}{k}=\log 2\left[\gamma-\frac{\log 2}{2}\right]$$

Differentiating with respect to $s$ we obtain

(4.4.116e) $$\sum_{n=0}^{\infty}\frac{1}{2^{n+1}}\sum_{k=0}^{n}\binom{n}{k}\frac{(-1)^k \log(k+x)}{(k+x)^s}=\sum_{k=0}^{\infty}\frac{(-1)^k \log(k+x)}{(k+x)^s}$$

and, with $x=1$, this becomes

(4.4.116f) $$\sum_{n=0}^{\infty}\frac{1}{2^{n+1}}\sum_{k=0}^{n}\binom{n}{k}\frac{(-1)^k \log(k+1)}{(k+1)^s}=\sum_{k=1}^{\infty}\frac{(-1)^k \log k}{k^s}$$

Hence we obtain for $s=1$

(4.4.116g) $$\sum_{n=0}^{\infty}\frac{1}{2^{n+1}}\sum_{k=0}^{n}\binom{n}{k}\frac{(-1)^k \log(k+1)}{(k+1)}=\log 2\left[\frac{\log 2}{2}-\gamma\right]$$

This concurs with proposition 8(a) in Coffey's paper [45c].

The left-hand side of (4.4.116) may be written as

(4.4.117) $$\int_0^{\infty}\frac{\tanh(x/2)}{xe^x}dx=\log(\pi/2)$$

This integral is not new: a more general form is contained in G&R [74, p.381] which states that

(4.4.117a) $$\int_0^{\infty}\frac{\tanh x}{xe^{\beta x}}dx=\log\frac{\beta}{4}+2\log\frac{\Gamma\left(\frac{\beta}{4}\right)}{\Gamma\left(\frac{\beta}{4}+\frac{1}{2}\right)}\quad\text{for Re}(\beta)>0$$

and letting $x = y/2$ and $\beta = 2$ we rediscover (4.4.117).

Differentiating (4.4.117a) with respect to $\beta$ we get

(4.4.117b) $$-\int_0^\infty \frac{\tanh x}{e^{\beta x}}dx = \frac{1}{\beta} + \frac{1}{2}\psi\left(\frac{\beta}{4}\right) - \frac{1}{2}\psi\left(\frac{\beta}{4} + \frac{1}{2}\right)$$

We now refer to Anglesio's [9] formula (4.4.93)

$$\int_0^\infty \frac{e^{-ax}(1-e^{-x})^n}{x^r} dx = \frac{(-1)^r}{(r-1)!}\sum_{k=0}^n \binom{n}{k}(-1)^k (a+k)^{r-1}\log(a+k)$$

where $a \geq 0$ and $1 \leq r \leq n$ (except for $a = 0, r = 1$).

We now let $r = 2$ to obtain

(4.4.117c) $$\int_0^\infty \frac{e^{-ax}(1-e^{-x})^n}{x^2} dx = \sum_{k=0}^n \binom{n}{k}(-1)^k (a+k)\log(a+k)$$

Making the substitution $y = e^{-x}$ (4.4.117c) becomes

(4.4.117d) $$\int_0^1 \frac{y^{a-1}(1-y)^n}{\log^2 y} dy = \sum_{k=0}^n \binom{n}{k}(-1)^k (a+k)\log(a+k)$$

Completing the summation from $n = 2$ we obtain

$$\sum_{n=2}^\infty \frac{1}{2^{n+1}}\sum_{k=0}^n \binom{n}{k}(-1)^k (k+a)\log(k+a) = \sum_{n=2}^\infty \frac{1}{2^{n+1}}\int_0^1 \frac{y^{a-1}(1-y)^n}{\log^2 y}dy$$

$$= \frac{1}{4}\int_0^1 \frac{(1-y)^2 y^{a-1}}{(1+y)\log^2 y}dy$$

Adding two more terms we get

$$\sum_{n=0}^\infty \frac{1}{2^{n+1}}\sum_{k=0}^n \binom{n}{k}(-1)^k (k+a)\log(k+a) =$$

$$\frac{1}{2}\sum_{k=0}^{0}\binom{0}{k}(-1)^k(k+a)\log(k+a)+\frac{1}{4}\sum_{k=0}^{1}\binom{1}{k}(-1)^k(k+a)\log(k+a)+\int_0^1\frac{(1-y)^2 y^{a-1}}{(1+y)\log^2 y}dy$$

and accordingly we have

$$\sum_{n=0}^{\infty}\frac{1}{2^{n+1}}\sum_{k=0}^{n}\binom{n}{k}(-1)^k(k+a)\log(k+a)=\frac{3}{4}\log a-\frac{1}{4}(1+a)\log(1+a)+\int_0^1\frac{(1-y)^2 y^{a-1}}{(1+y)\log^2 y}dy$$

With reference to (4.4.112za) where $a=1$ and $b=2$ we therefore have

$$\sum_{n=0}^{\infty}\frac{1}{2^{n+1}}\sum_{k=0}^{n}\binom{n}{k}(-1)^k(k+2)\log(k+2)-\sum_{n=0}^{\infty}\frac{1}{2^{n+1}}\sum_{k=0}^{n}\binom{n}{k}(-1)^k(k+1)\log(k+1)=2\varsigma_a'(-1)$$

(4.4.117e) $$2\varsigma_a'(-1)=\frac{1}{4}\log 2-\frac{3}{4}\log 3-\int_0^1\frac{(1-y)^3}{(1+y)\log^2 y}dy$$

Integrating (4.4.112n) we obtain

$$\sum_{n=0}^{\infty}\frac{1}{2^{n+1}}\sum_{k=0}^{n}\binom{n}{k}(-1)^k\log(k+u+1)-\sum_{n=0}^{\infty}\frac{1}{2^{n+1}}\sum_{k=0}^{n}\binom{n}{k}(-1)^k\log(k+1)=$$

$$u\log 2-\sum_{n=1}^{\infty}(-1)^{n+1}\frac{\varsigma_a(n+1)}{n+1}u^{n+1}$$

and with $u=1$ we have

$$\sum_{n=0}^{\infty}\frac{1}{2^{n+1}}\sum_{k=0}^{n}\binom{n}{k}(-1)^k\log(k+2)-\sum_{n=0}^{\infty}\frac{1}{2^{n+1}}\sum_{k=0}^{n}\binom{n}{k}(-1)^k\log(k+1)=\log 2-\sum_{n=1}^{\infty}(-1)^{n+1}\frac{\varsigma_a(n+1)}{n+1}$$

By (4.4.112k) the left-hand side is equal to $\log\frac{\pi}{2}$ and hence we get

(4.4.117f) $$\log\frac{\pi}{4}=\sum_{n=1}^{\infty}(-1)^n\frac{\varsigma_a(n+1)}{n+1}$$

and this concurs with (4.4.112zz).

We now refer back to Beumer's formula [23a]

(4.4.117g) $$\frac{1}{2}\left[\psi\left(\frac{x+2}{2}\right)-\psi\left(\frac{x+1}{2}\right)\right]=\log 2-\sum_{n=1}^{\infty}(-1)^{n+1}\left(1-\frac{1}{2^n}\right)\varsigma(n+1)x^n$$

and upon integration we obtain for $t \leq 1$

$$\frac{1}{2}\int_0^t\left[\psi\left(\frac{x+2}{2}\right)-\psi\left(\frac{x+1}{2}\right)\right]dx = t\log 2-\sum_{n=1}^{\infty}(-1)^{n+1}\left(1-\frac{1}{2^n}\right)\frac{\varsigma(n+1)}{n+1}t^{n+1}$$

$$= t\log 2-\sum_{n=2}^{\infty}(-1)^n\left(1-\frac{1}{2^{n-1}}\right)\frac{\varsigma(n)}{n}t^n$$

Whittaker & Watson [135, p.276] show that

(4.4.117h) $$\sum_{n=2}^{\infty}\frac{\varsigma(n,a)}{n}t^n = \log\Gamma(a-t)-\log\Gamma(a)+t\psi(a) \qquad (|t|<|a|)$$

and with $a=1$ we have

$$\sum_{n=2}^{\infty}\frac{\varsigma(n)}{n}t^n = \log\Gamma(1-t)-\log\Gamma(1)+t\psi(1) \qquad (|t|<1)$$

$$= \log\Gamma(1-t)-\gamma t$$

Letting $t \to -t$ we have

(4.4.117i) $$\sum_{n=2}^{\infty}(-1)^n\frac{\varsigma(n)}{n}t^n = \log\Gamma(1+t)+\gamma t$$

This Maclaurin series is derived in a different way in (E.22n).

Letting $t \to t/2$ in (4.4.117i) we have

(4.4.117j) $$\sum_{n=2}^{\infty}(-1)^n\frac{\varsigma(n)}{n2^n}t^n = \log\Gamma(1+t/2)+\frac{1}{2}\gamma t$$

Hence we obtain

$$\frac{1}{2}\int_0^t\left[\psi\left(\frac{x+2}{2}\right)-\psi\left(\frac{x+1}{2}\right)\right]dx = t\log 2 - \log\Gamma(1+t) + 2\log\Gamma(1+t/2)$$

and accordingly we have as in (4.4.112fa)

(4.4.117k) $$\frac{1}{2}\int_0^1\left[\psi\left(\frac{x+2}{2}\right)-\psi\left(\frac{x+1}{2}\right)\right]dx = \log 2 + 2\log\Gamma(3/2) = \log\frac{\pi}{2}$$

Referring to (4.4.117i) we see upon integration that

$$\sum_{n=2}^{\infty}(-1)^n\frac{\varsigma(n)}{n(n+1)}u^{n+1} = \int_0^u\log\Gamma(1+t)dt + \frac{1}{2}\gamma u^2$$

and making use of (4.4.112r)

$$\int_0^{1/2}\log\Gamma(x+1)dx = -\frac{1}{2} - \frac{7}{24}\log 2 + \frac{1}{4}\log\pi + \frac{3}{2}\log A$$

we obtain

$$\sum_{n=2}^{\infty}(-1)^n\frac{\varsigma(n)}{n(n+1)2^{n+1}} = -\frac{1}{2} - \frac{7}{24} + \frac{1}{4}\log\pi + \frac{3}{2}\log A + \frac{1}{8}\gamma$$

Differentiating (4.4.112ga) with respect to $\beta$

(4.4.117ki) $$\int_0^1\frac{(t^{\beta-1}-t^{\alpha-1})}{(1+t)\log t}dt = \sum_{n=0}^{\infty}\frac{1}{2^{n+1}}\sum_{k=0}^n\binom{n}{k}(-1)^k\log\frac{k+\beta}{k+\alpha}$$

we obtain G&R [74, p.896]

(4.4.117l) $$\int_0^1\frac{t^{\beta-1}}{1+t}dt = \sum_{n=0}^{\infty}\frac{1}{2^{n+1}}\sum_{k=0}^n\binom{n}{k}\frac{(-1)^k}{k+\beta} = \frac{1}{2}\left[\psi\left(\frac{1+\beta}{2}\right)-\psi\left(\frac{\beta}{2}\right)\right]$$

From (4.4.83) we have

$$\sum_{n=0}^{\infty}\frac{1}{2^{n+1}}\sum_{k=0}^n\binom{n}{k}\frac{(-1)^k}{(k+x)^s} = \sum_{n=0}^{\infty}\frac{(-1)^n}{(n+x)^s}$$

and hence we get

(4.4.117m) $$\int_0^1\frac{t^{\beta-1}}{1+t}dt = \sum_{n=0}^{\infty}\frac{(-1)^n}{n+\beta}$$

We have

$$\int_0^\infty \frac{t^{\beta-1}}{1+t}dt = \int_0^1 \frac{t^{\beta-1}}{1+t}dt + \int_1^\infty \frac{t^{\beta-1}}{1+t}dt$$

and making the substitution $t = 1/u$ we have

$$\int_1^\infty \frac{t^{\beta-1}}{1+t}dt = \int_0^1 \frac{u^{-\beta}}{1+u}du = \int_0^1 \frac{u^{\alpha-1}}{1+u}du$$

where $\alpha = -\beta + 1$ and, if $0 < \beta < 1$, then we also have $0 < \alpha < 1$. Therefore we obtain

$$\int_0^\infty \frac{t^{\beta-1}}{1+t}dt = \int_0^1 \frac{t^{\beta-1}}{1+t}dt + \int_0^1 \frac{t^{\alpha-1}}{1+t}dt$$

$$= \sum_{n=0}^\infty \frac{(-1)^n}{n+\beta} + \sum_{n=0}^\infty \frac{(-1)^n}{n+\alpha}$$

$$= \sum_{n=0}^\infty \frac{(-1)^n}{n+\beta} + \sum_{n=0}^\infty \frac{(-1)^n}{n+1-\beta}$$

$$= \frac{1}{\beta} + \sum_{n=1}^\infty \frac{(-1)^n}{n+\beta} + \sum_{n=0}^\infty \frac{(-1)^n}{n+1-\beta}$$

$$= \frac{1}{\beta} + \sum_{n=1}^\infty \frac{(-1)^n}{n+\beta} + \sum_{m=1}^\infty \frac{(-1)^{m-1}}{m-\beta}$$

Therefore we have

(4.4.117n) $$\int_0^\infty \frac{t^{\beta-1}}{1+t}dt = \frac{1}{\beta} + 2\beta \sum_{n=1}^\infty \frac{(-1)^{n-1}}{n^2 - \beta^2}$$

Using, for example, Fourier series we can show that (see (6.134) in Volume V)

(4.4.117o) $$\frac{\pi\beta}{\sin \pi\beta} = 1 + 2\beta^2 \sum_{n=1}^\infty \frac{(-1)^{n-1}}{n^2 - \beta^2}$$

Hence we obtain another proof of Euler's reflection formula for the gamma function (a further proof is given in Appendix C of Volume VI).

(4.4.117p) $$\int_0^\infty \frac{t^{\beta-1}}{1+t} dt = \Gamma(\beta)\Gamma(1-\beta) = \frac{\pi}{\sin \pi\beta}$$

We also note that from (4.4.117l) that

(4.4.117q) $$\frac{1}{2}\left[\psi\left(\frac{1+\beta}{2}\right) - \psi\left(\frac{\beta}{2}\right)\right] = \sum_{n=0}^\infty \frac{(-1)^n}{n+\beta}$$

and this appears in G&R [74, p.897].

Let us now reconsider the integral $I = \int_0^\infty \frac{t^{\beta-1}}{1+t} dt$ and, upon making the substitution $t = \alpha x$, we obtain

$$I = \alpha \int_0^\infty \frac{\alpha^{\beta-1} x^{\beta-1}}{1+\alpha x} dx$$

Differentiating with respect to $\alpha$ we get

$$\frac{\partial I}{\partial \alpha} = \int_0^\infty \frac{(1+\alpha x)\beta \alpha^{\beta-1} x^{\beta-1} - \alpha^{\beta-1} x^\beta}{(1+\alpha x)^2} dx$$

and letting $\alpha = 1$ we have

$$\beta \int_0^\infty \frac{x^{\beta-1}}{1+x} dx = \int_0^\infty \frac{x^\beta}{(1+x)^2} dx = 2\int_0^\infty \frac{u^{2\beta+1}}{(1+u^2)^2} du$$

which concurs with [25, p.193].

We now recall the Hasse/Sondow identity (3.11)

(4.4.117r) $$\varsigma_a(s) = \sum_{n=0}^\infty \frac{1}{2^{n+1}} \sum_{k=0}^n \binom{n}{k} \frac{(-1)^k}{(k+1)^s}$$

and differentiating (4.4.117b) we obtain

(4.4.117s) $$\varsigma_a'(s) = -\sum_{n=0}^\infty \frac{1}{2^{n+1}} \sum_{k=0}^n \binom{n}{k} \frac{(-1)^k \log(k+1)}{(k+1)^s}$$

Hence we have

$$(4.4.117t) \qquad \varsigma_a'(-1) = -\sum_{n=0}^{\infty} \frac{1}{2^{n+1}} \sum_{k=0}^{n} \binom{n}{k} (-1)^k (k+1) \log(k+1)$$

Donal F. Connon
Elmhurst
Dundle Road
Matfield
Kent TN12 7HD
dconnon@btopenworld.com